\documentclass[a4paper,12pt,reqno]{amsart}
\usepackage{amssymb}
\usepackage[dvips]{graphicx}
\usepackage{hyperref}

\input{poincare.def}

\title{Limiting behavior of trajectories of complex\\
polynomial vector fields}
\author{S. Ivashkovich}
\date{\today}

\address{
Universit\'e de Lille-1, UFR de Math\'ematiques, 59655 Villeneuve
d'Ascq, France} \email{ivachkov@math.univ-lille1.fr}
\address{IAPMM Nat. Acad. Sci. Ukraine
Lviv, Naukova 3b, 79601 Ukraine}

\begin{document}
\begin{abstract}
We prove that every trajectory  of a polynomial vector field on the
complex projective plane accumulates to the singular locus of the
vector field. This statement represents a holomorphic version of the
Poincar\'e-Bendixson theorem and solves the complex analytic
counterpart of Hilbert's $16$th  problem. The main result can be
also reformulated as the nonexistence of {\slsf "exceptional
minimals"} of holomorphic foliations on $\pp^2$ and, in particular,
implies the nonexistence of real analytic Levi flat hypersurfaces in
the complex projective plane. Finally, we describe (in the first
approximation) the way a {\slsf minimal} complex trajectory
approaches the singular locus of the vector field.
\end{abstract}

\maketitle

%\footnote{rty}
\bigskip\noindent{\bf Section 1. Introduction. \rm \longpoints pp.\ 2--10}

{\bgroup \baselineskip=8.5pt\parindent=0pt\small

\medskip
\centerline{\vbox{\hsize=5.0truein  \noindent 1.1.\ Statement of the
main result. 1.2.\ BLM-trichotomy. 1.3.\ The role of the holonomy
group. 1.4.\ Reduction to nef models. 1.5.\ Minimal sets in
projective spaces. 1.6. \ Levi flat hypersurfaces. 1.7. \
Exceptional minimals of Briot-Bouquet foliations. 1.8. Exceptional
minimals and Levi flats in Hirzebruch surfaces. 1.9. \ Approach of a
leaf to the singular locus. 1.10.\ Notes. 1.11.\ Acknowledgement. }}

\medskip\noindent{\bf Section 2. Pseudoconvexity of Poincar\'e domains.
\rm \longpoints pp.\ 10--21}

\smallskip \centerline{
\vbox{\hsize=5.0truein \noindent 2.1.\ Minimal sets of holomorphic
foliations. 2.2.\ Poincar\'e domains. 2.3.\ Docquier-Grauert
criterion and Fujita theorems. 2.4.\ Local pseudoconvexity of
Poincar\'e domains. 2.5.\ Nef models of holomorphic foliations.
2.6.\ Obstructions to the existence of universal covering Poincar\'e
domains. 2.7.\ A Rothstein type extension theorem.}}

\bigskip\noindent{\bf Section 3. Holomorphic representations of the
fundamental group. \rm \longpoints pp.\ 21--30}

\smallskip
\centerline{\vbox{\hsize=5.0truein  \noindent  3.1.\ A germ of the
holomorphic representation. 3.2.\ Expansion of the holomorphic
representation. 3.3.\ Universal covering Poincar\'e domains of
hyperbolic foliations. 3.4.\ Imbedding of the expanded Poincar\'e
domain into $\cc^2$. 3.5.\ Universal covering Poincar\'e domains of
parabolic foliations. }}

\bigskip\noindent{\bf Section 4. Proofs of the main results.
\rm \longpoints pp.\ 30--40}

\smallskip
\centerline{\vbox{\hsize=5.0truein  \noindent 4.1.\ BLM-trichotomy.
4.2.\ Proof of Theorem \ref{p2-inv-set}. 4.3.\ Limiting behavior of
leaves with hyperbolic holonomy. 4.4.\ Minimal leaves in the product
of projective lines. 4.5.\ Briot-Bouquet foliations. 4.6.\ Levi
problem in Hirzebruch surfaces. 4.7.\ Exceptional minimals and Levi
flats in Hirzebruch surfaces. 4.8.\ Minimal sets in projective
spaces. }}

\break
\bigskip\noindent{\bf Section 5. Psudoconvexity vs. rational curves.
\rm \longpoints pp.\ 40--46}

\smallskip
\centerline{\vbox{\hsize=5.0truein  \noindent 5.1.\ Analytic
objects. 5.2.\ Extension of analytic objects. 5.3.\ Pseudoconvexity
of the universal covering Poincar\'e domains. 5.4.\ Proof of Corollary
\ref{approach3}.}}

\bigskip\noindent{\bf References. \rm \longpoints \quad pp.\ 46--49}
}

%\tableofcontents
\newsect[sectINT]{Introduction}

\newprg[prgINT.proof]{Statement of the main result}
In this paper we study the limiting behavior of trajectories of
polynomial vector fields. Let $\v = P(x,y)\d/\d x + Q(x,y)\d/\d y$
be a complex vector field on $\cc^2$ with polynomial coefficients.
Trajectories of $\v$ define a holomorphic foliation $\call$, which
naturally extends onto the complex projective plane $\pp^2$. Vice
versa, every holomorphic foliation on $\pp^2$ is defined as the set
of trajectories of a polynomial vector field starting from an
appropriately chosen affine chart. In what follows we shall not
distinguish between trajectories of polynomial vector fields and
leaves of holomorphic foliations. Denote by $\Sing\call$  the
singular locus of a  holomorphic foliation $\call$ on a compact
complex surface $X$, \ie the set where the corresponding vector
field vanishes. This is a non empty finite subset if $X=\pp^2$. For
a point $m\not\in \Sing\call$ the leaf $\call_m$ through $m$ is, by
definition, the leaf of the smooth foliation
$\call^{\reg}\deff\call|_{X^{\reg}}$, where $X^{\reg}\deff
X\setminus\Sing\call$. If $m\in \Sing\call$ then leaves through $z$
are not defined, \ie a stationary point is not considered as being a
trajectory.

\smallskip
The main goal of this paper is to prove the following:

\begin{thm}
\label{p2-inv-set} Let  $\call$ be a holomorphic foliation on
$\pp^2$ and let $\call_m$ be any of its leaves. Then
\begin{equation}
\eqqno(intersect) \overline{\call_m} \cap \Sing\call \not=
\emptyset.
\end{equation}
\end{thm}

\smallskip The {\slsf limiting set} of the leaf $\call_m$
is defined as
\begin{equation}
\eqqno(lim-set) \lim \call_m \deff \bigcap_{K\comp \call_m}
\overline{\call_m\setminus K},
\end{equation}
where $K$ runs over all compact subsets of $\call_m$. Theorem
\ref{p2-inv-set} can be restated also in a following way: {\it the
limiting set of a leaf of a holomorphic foliation on $\pp^2$
intersects the singularity set of the foliation.} Recall that the
Poincar\'e - Bendixson theorem, in the form it was originally proved
by Poincar\'e in \cite{P}, states the following: {\it Let $\v$ be a
polynomial vector field on $\rr\pp^2$ (or, on $\sss^2$) and let
$\gamma$ be its trajectory. Then, either $\gamma$ is a periodic
trajectory (an orbit), or for each of the limiting sets
$\lim^{\pm}\gamma$ the following holds: either $\lim^{\pm}\gamma$ is
an orbit, or $\lim^{\pm}\gamma\cap \Sing \v\not=\emptyset$.} In its
turn the second part of the 16-th Hilbert's problem asks: {\it If
the number of limiting orbits of a polynomial vector field on
$\rr\pp^2$ of degree $n$ is bounded by some number, depending only
on $n$?} See \cite{H}.

\smallskip
The answer is still unknown. Theorem \ref{p2-inv-set} should be
viewed as an answer to a complex analytic counterpart of Hilbert's
problem: {\it the limiting set of a complex trajectory {\slsf always
accumulates} to the  singular locus of the vector field.} It was
known already for a long time that in the complex case there exist
no "orbits", \ie algebraic invariant curves in $X^{\reg}$, this
readily follows from the Camacho-Sad formula, see Appendix in
\cite{CS}. The problem was: if there could exist some leaves with
{\slsf massive} limiting sets away from singularities? This question
is mentioned in \cite{CLS1} and was explicitly posed in \cite{Ca} on
ICM-1990. Theorem \ref{p2-inv-set} tells that such sets {\slsf do
not exist}.

\newprg[INT.blm]{BLM trichotomy} Our point of depart for the proof
of Theorem \ref{p2-inv-set} will be the following trichotomy due to
Bonatti-Langevin-Moussu, see \cite{BLM}. A non-empty subset
$I\subset X$ is called {\slsf $\call$-invariant}, if it is not
contained in $\Sing\call$ and if for every point $m\in
I\setminus\Sing\call$ the leaf $\call_m$ of $\call$ through $m$ is
entirely contained in $I$. For example, the closure of each leaf is
a closed invariant set.  Limiting set of every leaf is also a closed
invariant set, unless it is entirely contained in $\Sing\call$. But
singular points $m\in\Sing\call$ are not considered as invariant
sets. A closed invariant set $\calm$ is called {\slsf minimal} if it
doesn't contains any proper closed invariant subset. Every closed
invariant set contains a (non unique in general) minimal subset.
Finally, every minimal set is the closure
of any of its leaves. If the leaf $\call_m$ is such that
$\overline{\call_m}$ is minimal we call it a {\slsf minimal leaf}.

\begin{rema} \rm
\label{exc-min} {\slsf "Exceptional minimals"} is a common name for
minimal sets of holomorphic foliations, which do not intersect the
singularity set.
\end{rema}
Let $(X,\call)$ be a foliated pair, where $X$ is a projective
surface, and let $\calm =\overline{\call_m}$ be an exceptional
minimal for $\call$ in $X$. Recall that every foliation on a complex
projective surface is defined by a {\slsf global} meromorphic form
$\Omega$. For such a form the divisors of poles $\Pole (\Omega )$
and zeroes $\Zero (\Omega )$ are correctly defined and both of them
contain $\Sing\call$.

\smallskip According to the result from \cite{BLM} if $\calm =
\overline{\call_m}$ is a leaf of a holomorphic foliation $\call$
such that $\overline{\call_m}$ doesn't intersects at least one of
$\Pole (\Omega)$ or $\Zero (\Omega )$ then one has only the three
following possibilities:

\begin{itemize}

\smallskip\item the defining meromorphic from $\Omega$ is (algebraically)
closed;

\smallskip\item $\call_m$ is a compact leaf;

\smallskip\item there exists a point $n\in \overline{\call_m}$ with
$\overline{\call_n}=\overline{\call_m}$ such that the leaf $\call_n$
contains a loop with hyperbolic holonomy.
\end{itemize}

\noindent First two possibilities do produce exceptional minimals,
see examples in Section \ref{sect.PROOF}. But in $\pp^2$ they are
obviously not possible, see Subsection \ref{prgPROOF.blm}. Therefore
we must exclude the third case. Recall that a loop $\gamma \in
\pi_1(\call_m,m)$ is said to have a hyperbolic holonomy if the
derivative of its holonomic representative is by modulus less then
one.

\newprg[prgINT.abel]{The role of the holonomy group}

For the sake of simplicity we suppose at most instances along this
discussion that $X = \pp^2$. Let $\call_m$ be a minimal leaf and let
$D$ be a Poincar\'e disc through $m$, \ie an image of $\Delta\subset
\cc$ under a holomorphic imbedding into $X$, such that it is
transversal to the leaves of $\call$. Let $\calp_D$ be the
Poincar\'e domain of $\call$ over $D$, \ie $\calp_D = \bigcup_{z\in
D}\call_z$ - the union of leaves cutting $D$. This is an open subset
of $X$, if $\Sing\call\not= \emptyset$ (ex. $X=\pp^2$) then this set
is also a proper subset of $X$. Would it so happen that $\calp_D$ is
pseudoconvex (in $\pp^2$ that means Stein) and $\overline{\call_m}$
doesn't intersects $\Sing\call$ we would get a contradiction.
Indeed, from one side $\overline{\call_m}$ should be contained in
$\calp_D$ by minimality, from another side - this contradicts to the
maximum principle. Therefore we are left with the option when
$\calp_D$ is not pseudoconvex {\slsf for every} $D$. Corollary
\ref{cor-ps2} now states the following:

\smallskip{\slsf  If $X=\pp^2$ and $\call$ is a holomorphic foliation,
which contains an exceptional minimal $\calm = \overline{\call_m}$,
then for a given Poincar\'e disc $D\ni m$ the Poincar\'e domain
$\calp_D$ can be non pseudoconvex only when all points of
$\Sing\call$ are the {\slsf isolated boundary points} for $\calp_D$.
Moreover, in that case $\calp_D=\pp^2\setminus \Sing\call$.}

\smallskip The only use that we shall make from this first observation
is that $\calp_D$ is simply connected in this case and contains a
lot of rational curves. Examples with $\calp_D=\pp^2\setminus
\Sing\call$ do exist, see \cite{LR}. But up to this moment we didn't
use the restrictive properties of the holonomy group. Now, up to
taking another point $n\in\overline{\call_m}$ and a leaf $\call_n$
(with the same closure) we can suppose that $\hln (\call_m, m)$
contains a hyperbolic (\ie contracting) element $\alpha$. Therefore
one can choose a complex coordinate $t$ in $D$ such that this
element will become to be a multiplication by some $\alpha \in\cc^*$
with $|\alpha|<1$ (but the whole group $\hln (\call_m, m)$ might be
even non abelian).

\begin{rema} \rm
\label{cerveau} It is an appropriate moment to mention the
remarkable dichotomy of Cerveau, who proved that either the  {\slsf
whole} holonomy group $\hln (\call_m,m)$ is abelian or, our
exceptional minimal $\calm =\overline{\call_m}$ is an
$\call$-invariant Levi flat hypersurface! See \cite{C}.
Unfortunately we are not able to make {\slsf any} use from $\calm$
being a Levi flat hypersurface and our proof uses only the existence
of {\slsf some} loop in $\pi_1(\call_m,m)$ with contractible
holonomy. Levi flats will be ruled out together will exceptional
minimals.
\end{rema}

Following Ilyashenko \cite{Iy1}, we consider
the universal covering Poincar\'e domain (a ''skew cylinder`` in the
terminology of \cite{Iy1}): $\tilde\calp_D \deff \bigcup_{z\in D}
\tilde\call_z$ - the union of universal coverings of leaves cutting
$D$. The natural topology on $\tilde\calp_D$ might be non Hausdorff.
The obvious reason is the possible presence of vanishing cycles. But
the more deeper reason is the presence of some special
$\call$-invariant rational curves in $X$. After bringing $(X,\call)$
to a {\slsf nef} model $(Y,\calf)$ we get reed from this problem. In
the case when $\tilde\calp_D$ is Hausdorff the natural universal
covering maps $\tilde p_z:\tilde\call_z\to \call_z$ glue together to
a locally biholomorphic foliated projection $\tilde p :
\tilde\calp_D\to \calp_D\subset X$. This aloud to consider the pair
$(\tilde\calp_D,\tilde p)$ as a Riemann domain over $X$.

\smallskip One should  take care at this point, because the nef model
can have singularities other than $\Sing\call$, \ie the surface $X$
itself can become singular. Also one should check that nothing
changes in our initial data. All this is not really difficult and
therefore we continue the explanation of our proof assuming that
$(X, \call)$ is {\slsf already nef}.

\smallskip The key observation of this paper can be described
as follows: using $\alpha$ we {\slsf expand} the universal covering
Poincar\'e domain $\tilde\calp_D$ to a foliated domain
$\tilde\calp_{\cc}$ over the whole complex plane $\cc$. Let us
formulate the precise statement. By a foliated domain we mean a
triple $(W,\pi, S)$, where $W$ is a connected complex surface, $S$
is a complex curve and $\pi :W \to S$ is a holomorphic submersion
with connected fibers. For example the universal covering
$\tilde\calp_D$ of a Poincar\'e domain $\calp_D$ is foliated over
$D$, the corresponding map is denoted as $\tilde\pi
:\tilde\calp_D\to D$, see Section \ref{sect.PPD} for more details.
In Section \ref{sectAB} we prove the following:

\begin{thm}
\label{hol-rep} Let $\call_m$ be a leaf of a holomorphic foliation
on a compact complex surface such that $\hln (\call_m, m)$
contains a hyperbolic element $\alpha$. Suppose furthermore, that
$\tilde\call_D$ is Hausdorff and Rothstein. Then there exists a foliated
domain $\tilde\calp_{\cc}$ over $\cc$ such that:

\smallskip\sli $\tilde\calp_D$ is a foliated subdomain of
$\hat\calp_{\cc}$, \ie $\tilde\pi$ extends to $\tilde\calp_{\cc}$,
and, moreover, $\tilde p:\tilde\calp_D\to X$

\quad also extends to $\tilde\calp_{\cc}$;

\smallskip\slii $\tilde\calp_{\cc}$ is periodic, \ie $\alpha$
lifts to a foliated, $\tilde p$-invariant biholomorphism
$\tilde\alpha$ of $\tilde\calp_{\cc}$.
\end{thm}

For the notion of {\slsf rothsteiness} we refer to Subsection
\ref{prgPPD.roth}. For the moment let us say that every Stein
manifold (or, normal space) is Rothstein.

\newprg[prgINT.nef]{Reduction to nef models} Having in our disposal
the {\slsf expanded} Poincar\'e domain we are in the position to
deploy the powerful ''Noncommutative Mori theory`` of McQuillan, see
\cite{McQ1}, \cite{McQ2}. One of the main results of this theory
tells that, after reducing $(X,\call)$ to the {\slsf nef model}, we
find ourselves  under the following alternative:

\begin{itemize}
 \item either all leaves of $\call$ are parabolic, \ie covered by $\cc$;

\smallskip\item or, the set $\{ \text{parabolic leaves}\}\cup \Sing\call$ is a
proper algebraic subset $\cala$ of $X$, and, moreover, the hyperbolic
distance along the (hyperbolic!) leaves in $X\setminus\cala$ is
continuous.
\end{itemize}

This is not a very precise statement. One should exclude the case of
rational quasi-fibrations and one should also take care about the
items already mentioned above (\ie of the fact that the nef model may
have cyclic singularities). But all this doesn't really matter
and therefore we continue the explanation of our proof assuming that
$(X, \call)$ is already a nef model, \ie that dichotomy of McQuillan
actually takes place already for our $\call$.

\smallskip\noindent{\slsf 4. Parabolic case.} In parabolic case we
have the following possibilities for a leaf with hyperbolic
holonomy, see Theorem \ref{parab-case}:

\begin{itemize}

\smallskip\item $\call$ is a rational quasi-fibration (not excluded
by McQuillan's alternative).

\smallskip\item $\call_m$ is compact, \ie is a torus or, a projective
line.

\smallskip\item $\call_m$ is biholomorphic to $\cc^*$ and is a
{\slsf locally closed  analytic subset} of $X\setminus \Sing\call$.
\end{itemize}

\begin{rema} \rm
\label{parab-rema} When stating these possibilities we do not
suppose that $\overline{\call_m}$ is neither exceptional nor even
that it is minimal.
\end{rema}

In all these cases it is easy to determine what happens with our
$\calm=\overline{\call_m}$ when $X=\pp^2$. First case is trivial. In
the case of $\pp^2$ the second cases doesn't happens. Third case
obviously leads to a contradiction in the case when minimality and
exceptionality are additionally assumed.

\smallskip\noindent{\slsf 5. Hyperbolic case.}
In hyperbolic case we prove, see Theorem \ref{prop-disc-hyp}, that:

\begin{itemize}
\smallskip\item $\tilde p:\tilde\calp_{\cc}\to \calp_D$ is a regular
covering.
\end{itemize}

If, for example, $X$ is simply connected (ex. $\pp^2$ is such) this
last fact quickly leads to a contradiction. Indeed, if $X=\pp^2$ we
have $\calp_D=\pp^2\setminus \Sing\call$, as it was noticed, \ie
$\calp_D$ is also simply connected. Therefore
$\tilde\calp_{\cc}=\calp_D$ by monodromy. But $\tilde\calp_{\cc}$ is
foliated over $\cc$ by a submersive holomorphic map $\tilde
\pi:\tilde\calp_{\cc}\to \cc$. This map then extends through
$\Sing\call$ onto $\pp^2$ and, therefore, should be constant.
Contradiction.

\smallskip In reality the last argument is a bit more technical and
makes more use of rational curves in $\pp^2\setminus \Sing\call$
then of simple connectivity, see Subsection \ref{prgNEF.p2} for more
details.

\newprg[prgINT.pn]{Minimal sets in projective spaces}
Starting from a polynomial vector field on $\cc^2$ one can
compactify it to a holomorphic foliation as well on $\pp^2$ as on
$\pp^1\times \pp^1$, for example. Therefore let us state the
following:

\begin{corol}
\label{p1p1-inv-set} Let $\call$ be a holomorphic foliation on
$\pp^1\times\pp^1$ and let $\call_m$ be its minimal leaf. Then:

\smallskip\sli either $\overline{\call_m}\cap \Sing\call \not=
\emptyset$;

\smallskip\slii or, $\call$ is a rational fibration and $\call_m$
is its fiber.
\end{corol}

The fact that the closure of a leaf of a codimension one holomorphic
foliation on $\pp^n$, for $n\ge 3$, must intersect $\Sing\call$ was
proved in \cite{L}. In fact, the principal result of \cite{L} reads
as follows: {\sl the singular set $\Sing\call$ of a holomorphic
codimension one foliation on $\pp^n$, $n\ge 2$, has at least one
irreducible component of codimension two.} Since the complement to
$\overline{\call_m}$ in $\pp^n$ is Stein by Fujita's theorem, we see
readily that $\Sing\call\cap \overline{\call_m}\not= \emptyset$
provided $n\ge 3$. In dimension two, however, the Stein domain
$\pp^2\setminus \overline{\call_m}$ may well contain codimension two
analytic sets $ = $ finite  sets of points. Therefore the Theorem
\ref{p2-inv-set} could be true only because the closures of leaves
of holomorphic foliations on $\pp^n$ posed more restrictive
properties with respect to the singular locus of the foliation then
just to ``intersect it''. Indeed, if one looks on the statement of
Theorem \ref{p2-inv-set} from the point of view of higher dimensions
then it states that: $\overline{\call_m}$ {\sl {\slsf contains} at
least one irreducible component of $\Sing\call$ of codimension two.}
This should be true in general. In this paper we provide a partial
step in this direction:

\begin{corol}
\label{pn-inv-set} Let $\call$ be a codimension one holomorphic
foliation on the complex projective space $\pp^n$, $n\ge 3$, and let
$\call_m$ be any of its leaves. Then $\overline{\call_m}$ intersects
at least one irreducible component of $\Sing\call$ of codimension
two by a closed set of positive $(2n-4)$-dimensional Hausdorff
measure.
\end{corol}

The proof follows from Theorem \ref{p2-inv-set} by taking generic
sections and is given in Subsection \ref{prgPROOF.pn}.

\newprg[INT.levi]{Levi flat hypersurfaces}

A special case (of a special interest) of the ''exceptional minimals
problem`` is the question of (non)existence of Levi flat
hypersurfaces in certain complex manifolds (like $\pp^2$). Recall
that a real hypersurface $M$ in a complex manifold $X$ is called
{\slsf Levi flat} if $M$ is locally foliated by complex
hypersurfaces. Equivalently, if $M$ locally divides  $X$ onto
pseudoconvex parts. In this paper we consider only {\slsf real
analytic} $M$-s, if the opposite is not explicitly stated, and, one
more point: our hypersurfaces will be always {\slsf compact}.

\smallskip To start with let us remark that given a real analytic
Levi flat hypersurface $M$ in a complex manifold $X$ there exists a
neighborhood $U$ of $M$ and a holomorphic foliation $\call$ on $U$
by complex hypersurfaces which extends the Levi foliation of $M$.
\begin{rema} \rm
\label{inv-logic1} If $M$ appeared as a minimal set of a holomorphic
foliation, the existence of $\call$ doesn't comes into a question.
But it is also true (and easy) that a Levi flat real analytic $M$
itself induces a holomorphic foliation, but {\slsf only on a
neighborhood} of it, in general.
\end{rema}

\smallskip\noindent{\slsf 1. A Levi-type problem.} Recall that the Levi
problem on a (not necessarily compact) complex manifold $X$ consists
in finding the necessary and sufficient conditions on a relatively
compact subdomain $D\comp X$ to be Stein. The first step in our
approach can be described as follows:

\begin{itemize}
\item Either all components of $X\setminus M$
are "convex enough" (ex. Stein), or the foliated pair $(U,\call)$ is
degenerate in a neighborhood of $M$.
\end{itemize}

\smallskip We are not going to make this statement more precise:
a substantial amount of classical results on the Levi problem show
that the "convex" case is quite typical in this setting. Therefore
the failure of $X\setminus M$ to be "convex enough" is the {\slsf
first raison d'\^etre} for a Levi flat hypersurface in a compact
complex manifold. Remark that the example of Grauert, see
\cite{Na2}, was viewed as an example of pseudoconvex manifold which
doesn't carries non-constant holomorphic functions. In this example
one has a Levi flat hypersurface $M$ in a complex torus $\ttt^2$  of
dimension two such that $D=\ttt^2\setminus M$ is the said
pseudoconvex manifold (and the corresponding foliation is clearly
''degenerate``).

\smallskip\noindent{\slsf 2. BLM-trichotomy.} On the other side, if
the components of $X\setminus M$ are ''convex enough`` (like in the
case of $\pp^2$), then the Levi foliation extends to a holomorphic
foliation $\call$ on the whole of $X$.

\begin{itemize}
 \item In that case $M$ should contain an ''exceptional minimal``
 $\calm$ of $\call$.
\end{itemize}

This is the  {\slsf second raison d'\^etre} for a Levi flat
hypersurface in a compact complex manifold. Remark that examples of
{\slsf Levi scrolls} of Ohsawa, see \cite{Oh3},  are of that kind.
Then one should check BLM-trichotomy for $\calm$. In our setting
this simply means that either $\calm = M$ and then we apply our
machine directly to $M$, or $\overline{\call_m}$ is a proper closed
invariant subset of $M$ and we derive the global behavior of $M$
from that of $\calm = \overline{\call_m}$. In this way one easily
obtains the following corollaries.

\begin{corol}
\label{proj} Complex projective plane $\cc\pp^2$ doesn't contains
any real analytic Levi flat hypersurface.
\end{corol}

\begin{corol}
\label{prod-proj} Let $M$ be a real analytic Levi flat hypersurface
in $\pp^{1}\times \pp^{1}$. Then $M=\gamma\times\pp^{1}$ (or $\pp^1
\times\gamma$), where $\gamma$ is a closed, real analytic curve  in
$\pp^1$.
\end{corol}

Proofs are given in Remarks \ref{proof-proj} and
\ref{proof-prod-proj}.

\newprg[prgINT.bb]{Exceptional minimals of Briot-Bouquet foliations}

Let us emphasis that neither $\pp^2$ nor $\pp^1\times\pp^1$ are
always natural manifolds to carry a holomorphic foliation, which
comes from algebraic differential equation. Consider, for example,
the polynomial equation of Briot and Bouquet, studied for the first
time in \cite{BB}. It is an equation of the form
\begin{equation}
\eqqno(bb) F(z,z')=0,
\end{equation}
where $F$ is an irreducible polynomial of two complex variables, non
constant both in $z$ and $z'$. This class of equations includes the
Riccati equation, equation for elliptic curves and so on. In
\cite{Dn},  after the geometrization of the problem, it was showed
that \eqqref(bb) naturally raises to a holomorphic foliation on the
product $\Sigma_g\times \pp^1$, where $\Sigma_g$ is the compact
Riemann surface of genus $g\ge 1$. The affine part $\Sigma^0_g$ of
$\Sigma_g$, in an appropriate affine coordinates $p,q$, is given  by
the equation $F(p,q)=0$. It is natural to call holomorphic
foliations on $\Sigma_g\times \pp^1$ - the {\slsf Briot-Bouquet
foliations}.  We prove the following:

\begin{corol}
\label{bb-min}
The only exceptional minimals of Briot-Bouquet
foliations are fibers $E_1\deff \{\pt\}\times \pp^1$ or $E_2 \deff
\Sigma_g\times \{pt\}$.
\end{corol}

The proof is given in Subsection \ref{prgPROOF.bb}.

\newprg[prgINT.hirz]{Exceptional minimals and Levi flats in
Hirzebruch surfaces}

In Subsections \ref{prgPROOF.hirz1} and \ref{prgPROOF.hirz2} we
classify the exceptional minimals and Levi flats in Hirzebruch
surfaces:

\begin{corol}
\label{inv-set-hirz} Let $\calm$ be a minimal set of a holomorphic
foliation on the Hirzebruch surface $H_k$, $k\ge 1$. Then:

\smallskip\sli either $\calm\cap \Sing\call \not= \emptyset$;

\smallskip\slii or, $\call$ is the rational fibration $\call_{s} =
\pi^{-1}(s)$ and $\calm$ is one of its fibers.
\end{corol}

\begin{corol}
\label{levi-fl-hirz} Let $M$ be a real analytic Levi flat
hypersurface in Hirzebruch surface $H_k$. Then $M=\pi^{-1}(\gamma)$,
where $\gamma$ is a real analytic imbedded loop in $\pp^1$.
\end{corol}

\newprg[prgINT.approach]{Approach of a leaf to the singular locus}

The classical  Poincar\'e - Bendixson theory apart of the
description of the limiting behavior of a trajectory of a vector
field on the real plane, which stays away from the singular locus of
the vector field (it accumulates to an orbit), describes also the
way a trajectory behaves when accumulating to the singular locus of
a vector field in question.  Since we proved that in the complex
case the first option is impossible, \ie a complex trajectory {\slsf
always} accumulates to the singular locus of the vector field, it is
natural to provide a step towards the description of a way a complex
trajectory approaches the singular locus.

\smallskip In the first approximation we describe the limiting behavior
of a {\slsf minimal} leaf of a holomorphic foliation in terms of its
invariant neighborhoods - Poincar\'e domains $\calp_D$ associated to
transversal discs cutting $\call_m$ at $m$ (or at any other point)
and their {\slsf leafwise universal coverings}.  Our goal is to
specify the behavior of the Poincar\'e domains themselves and/or the
universal covering Poincar\'e domains (at least), and, even more, of
the leaf $\call_m$ itself. Let us see that the case of hyperbolic
holonomy is practically done. Let $\hln : \pi_1(\call_m,m)\to \dif
(D,m)$ be the holonomy representation. Its image is the holonomy
group $\hln (\call_m,m)$ of the leaf $\call_m$. For $\gamma\in
\pi_1(\call_m,m)$ one can take the derivative $\hln (\gamma)'(m)$ of
its holonomy representative at $m$ and thus obtain a homomorphism
$\hln' : \pi_1(\call_m,m) \to \cc^*$.

\begin{defi}
\label{hyp-par} We say that the holonomy of $\call_m$ is {\slsf
parabolic} if $\hln' \big(\pi_1(\call_m,m)\big)\subset \ss^1$. If it
is not the case we say that the holonomy of $\call_m$ is {\slsf
hyperbolic}.
\end{defi}
The latter case means that there exists a loop $\gamma \subset
\call_m$, starting and ending at $m$, such that $|\hln'(\gamma
)(m)|<1$. One says also that such $\gamma$ has the {\slsf
contracting} or {\slsf hyperbolic} holonomy.

\begin{corol}
\label{approach1} Let $\call_m$ be a minimal leaf of a holomorphic
foliation on $\pp^2$ with hyperbolic holonomy. Then:

\smallskip\sli either $\overline{\call_m}\setminus \Sing\call$ posed
a Stein invariant neighborhood;

\smallskip\slii or, $\overline{\call_m}$ is a rational curve, which cuts
$\Sing\call$ by at least two points;

\smallskip\sliii or, $\overline{\call_m}$ is an elliptic curve, which cuts
$\Sing\call$ by at least one point.
\end{corol}
Options (\slii and (\sliii mean, in another words, that $\call_m =
C\setminus \{p_1,...,p_k\}$, where $p_j\in \Sing\call$ for every $j$ and $C$
is a rational or elliptic curve
in $\pp^2$ such that $C\setminus \{p_1,...,p_k\}$ is imbedded as a
closed analytic subset in $\pp^2\setminus \{p_1,...,p_k\}$.
Moreover, for every $p_j$ there exists a closed subset $E_j$ in
$\call_m$ such that $E_j$ is biholomorphic to the punctured disc
$\check\Delta \deff \{z\in\cc : 0<|z|\le 1\}$ (the so called
{\slsf vanishing end}). Finally, $\overline{
E_j} \cap \Sing\call = \{p_j\}$. Note that the possibility that
$p_i=p_j$ for some $i\not = j$ is not excluded (but $E_i\cap E_j=\emptyset$
for all $i\not= j$).

\smallskip The proof of this corollary repeats step to step the
proof of Theorem \ref{p2-inv-set} and is given in Subsection
\ref{prgPROOF.approach1}. In general, \ie without the assumption of
hyperbolicity we can prove a much weaker result.

\begin{thm}
\label{approach2} Let $\call$ be a holomorphic foliation on $\pp^2$
and let $\call_m$ be a minimal leaf of $\call$ through some point
$m\in X^{\reg}\deff \pp^2\setminus \Sing\call$. Then:

\smallskip\sli either for a sufficiently small transversal disc
$D$ through $m$ the universal covering

\quad Poincar\'e domain $\tilde\calp_D$ is a Stein;

\smallskip\slii or, the closure $\overline{\call_m}$ of $\call_m$
is a rational curve, cutting $\Sing\call$ by exactly one point.
\end{thm}

The second case represents a satisfactory description of a limiting
behavior of a minimal leaf. At the same time it gives a precise
obstruction for $\tilde\calp_D$ to be Stein. The proof of this
theorem requires much more then that of the preceding corollary
(modulo Theorem \ref{p2-inv-set} of course) and in given in the last
Section \ref{sect.RAT} of this paper.

\smallskip Theorem \ref{approach2} aloud  to make more precise the
case (\sli of Corollary \ref{approach1}. The point is that the Stein
invariant neighborhood, which occurs there is equal to the
Poincar\'e domain $\calp_D$ for a sufficiently small disc through
$m$ plus a finite number $\{s_1,...,s_d\}$ of isolated boundary
points of $\calp_D$. The nature of this points is given by the
following:

\begin{corol}
\label{approach3} In the case (\sli of Corollary \ref{approach1} all
$\{s_1,...,s_d\}$ are dicritical for $\call$.
\end{corol}

The proof is given in Subsection \ref{prgRAT.approach}.

\newprg[INT.notes]{Notes. } We are certainly even not trying here to
give any sort of review  on an unobservable amount of literature on
general, \ie not only holomorphic Poincar\'e - Bendixson theory.
Issues concerning Hilbert's $16$th  problem the interested reader
may consult in \cite{Iy3}. We only want to track a bit the {\slsf
incomplete in all senses} history of the {\slsf holomorphic} case.

\smallskip\noindent{\slsf 1.} The nature of ``periodic'' solutions of
complex polynomial equations is studied since the time when these
equation are studied themselves. For example, \cite{BB} stays
continuously to be a paper of reference, see ex. \cite{ELN} and
references there..

\smallskip\noindent{\slsf 2.} In his celebrated paper \cite{P}
Poincar\'e, influenced by \cite{BB} and preceding works of Cauchy
(see his recours on p. 385 of \cite{P}), proved his famous theorem
for {\slsf polynomial} vector fields in real dimension two.

\smallskip\noindent{\slsf 3.} The paper \cite{Be} of Bendixson, where
he proved a $\calc^1$-generalization  of the Poincar\'e's theorem
is, in fact, mostly devoted to the equations with holomorphic
coefficients: Ch. I - to, what is now known as Poincar\'e -
Bendixson theorem,  Ch-s II - VII - to holomorphic differential
equations. At the end of the Introduction Bendixson shows that he is
well aware that the situation with the Poincar\'e's famous theorem
is unclear in the complex case, saying that: ``Nous nous bornons ici
\`a cette remarque, voulant dans ce m\'emoire traiter seulement les
{\slsf courbes int\'egrales  r\'eelles } des \'equations
diff\'erentielles''. Non-integral, \ie local trajectories, which he
calls characteristics, are studied in \cite{Be} for the case of
holomorphic vector fields.

\smallskip\noindent{\slsf 4.} In \cite{AGH} the authors discovered
that flows on nilmanifolds provide a rich source of examples of
minimal sets. In   \cite{So} it was proved that if a connected
solvable Lie group acts holomorphically on a compact K\"ahler
manifold $X$ with $H^1(X,\cc)=0$ then this action has a fixed point
in every invariant complex subspace. In \cite{Lb} further study of
minimal sets of holomorphic actions of Lie groups on compact
K\"ahler manifolds was undertaken.

\smallskip\noindent{\slsf 5.} After the appearance of \cite{CLS1}
the problem of the existence of minimal invariant sets, which do not
intersect the singular locus of the foliation, became very popular.
Let us mention only \cite{BLM}, \cite{C} and, finally, for $n\ge 3$
the statement of Theorem \ref{p2-inv-set} was proved in \cite{L}.

\smallskip\noindent{\slsf 6.} After the dichotomy of Cerveau was
discovered, see \cite{C}, the attention was switched  mainly to the
question of (non)existence of Levi flats in $\pp^2$. But it should
be said that the Levi flat hypersurfaces first appeared in complex
analysis long before, to my best knowledge in the example of Grauert
of pseudoconvex manifold without non constant holomorphic functions,
see \cite{Na1} and \cite{Na2}.

\smallskip\noindent{\slsf 7.} Levi flats as common boundaries of two
pseudoconvex domains were studied in \cite{Shc} and other papers. As
natural boundaries of envelopes of holomorphy of bounded tube
domains they appeared \cite{Iv1}. Based on these papers and several
others (see \cite{Vi} for more details on this activity) the
question of existence of Levi flats in $\pp^2$ was actively
discussed on the seminar on Complex Analysis in Moscow University in
mid 80-s. \cite{Mi} is a track of this activity. This list is very
far from being complete, see ex. \cite{Bd,R}.

\smallskip\noindent{\slsf 8.} The case of smooth, \ie non real analytic
hypersurfaces in dimension $\ge 3$ was excluded by Siu in
\cite{Si3}. The case of dimension $2$ remained open despite of
several very clever attempts made in \cite{Oh1, Si4} and several
other papers dependant on these ones, see historical sketch in
\cite{IM}.

\smallskip\noindent{\slsf 9.} Reduction to nef models was developed
in the series of prominent works, among which I bound myself with
mentioning Seidenberg \cite{Sei}, Miyaoka \cite{Miy},  McQuillan
\cite{McQ1}, \cite{McQ2} and Brunella \cite{Br5}. Poincar\'e domains
(under the name ``skew cylinders'') where introduced in the case of
foliations on Stein manifolds by Ilyashenko in \cite{Iy1} and
studied by him and his students. In the non Stein setting Poincar\'e
domains (under the name of "covering tubes") were studied by
Brunella. The study of a way how a nearby trajectory cuts a
transversal interval on the real plane - is the central idea of the
paper of Poincar\'e \cite{P}. Therefore I think it is right to
associate the analogous objects, also in complex setting, with his
name.

\newprg[prgINT.ack]{Acknowledgement} I am grateful to Dmitri Markushevich
for numerous helpful discussions around complex algebraic aspects of
this paper, as well as to Stefan Nemirovski for sending me his
unpublished preprint \cite{Ne2} and, especially, for pointing out to
me the gap in the first version of this paper. From Frank Loray and
Bertrand Deroin I learned interesting examples of foliations on
$\pp^2$. I am especially grateful to Michael McQuillan for giving to
me explanations on his Noncommmutative Mori theory. Also to Vsevolod
Shevchishin and Alexandre Sukhov, who were the first to listen my
{\slsf expos\'es} on the subject of this paper, for their patience
and criticism.

\newsect[sect.PPD]{Pseudoconvexity of Poincar\'e domains}

\newprg[prgPPD.min]{Minimal sets of holomorphic foliations.}

We consider foliated manifolds, \ie pairs $(X,\call)$, where $X$ is
a complex manifold and $\call$ is a (singular) codimension one
holomorphic foliation on $X$. Later we shall aloud $X$ to have the
so called {\slsf cyclic quotient} singularities, but this will not
introduce any complications into the our exposition. There are
several equivalent ways to define a holomorphic foliation on a
complex manifold. We shall use two of them. First reads as follows.
A codimension one holomorphic foliation $\call$ on a complex
manifold $X$ is given by:

\smallskip\sli an open covering $\{U_j\}$ of $X$ and non identically
zero holomorphic $1$-forms $\omega_j$ on $U_j$;

\smallskip\slii  forms $w_j$ satisfy {\slsf integrability condition}
$w_j\wedge dw_j = 0$ (no condition if $\dim X = 2$);

\smallskip\sliii on intersections $U_j\cap U_k$ the defining forms $w_j$
are proportional, \ie $\omega_j = f_{jk}\omega_k$ for some
$f_{jk}\in \calo^*(U_j\cap U_k)$.

\medskip The set $\{\omega_j =0\}$ is the singularity set $\Sing\call$
of $\call$. Up to reducing common factors of coefficients of
$\omega_j$ the singularity set can be supposed to be of complex
codimension at least two. In regular part $X^{\reg}\deff X\setminus
\Sing\call$ of $X$ the leaves $\call_z$ of $\call$ are defined by
the equations $\omega_j|_{T\call_z} = 0$. The cocycle $\{f_{jk}\}\in
H^1(X,\calo^*)$ defines the {\slsf conormal} bundle
$\caln^*_{\call}$ of $\call$. Second definition (in dimension two)
reads as follows:

\smallskip\sli given an open covering $\{U_j\}$ of $X$ and non
identically zero holomorphic vector fields $\v_j$ on $U_j$;

\smallskip\slii on intersections $U_j\cap U_k$ the vector fields
$\v_j$ are proportional, \ie $\v_j = g_{jk}\v_k$ for some $g_{jk}\in
\calo^*(U_j\cap U_k)$.

\smallskip Leaves of $\call$ are locally the complex curves tangent
to $\v_j$, $\Sing\call = \{\v_j = 0\}$. The cocycle $\{g_{jk}\}\in
H^1(X,\calo^*)$ defines a holomorphic line bundle, which is called
the {\slsf canonical} bundle of $\call$. It will be denoted as
$\calk_{\call}$ or, simply as $\calk$ if no misunderstanding could
occur.

\begin{defi}
\label{inv-set-def} A subset $I\subset X$ of a foliated manifold is
called $\call$-invariant or, simply invariant, if $\call$ is clear
from the context, if:

\smallskip\sli $I\setminus \Sing\call$ is non empty;

\smallskip\slii for every point $z\in I\setminus\Sing\call$ the leaf
$\call_z$ is entirely contained in $I$.
\end{defi}
Would we ask only the condition (\slii in this definition then,
since {\slsf no} leaves are defined through singular points,  any
subset $S\subset \Sing\call$ would formally satisfy such definition.
But sets in $\Sing\call$ are not the objects of our study and they
are not considered to be the invariant sets.

\begin{defi}
\label{minimal-def} A closed invariant set, that doesn't contains
any proper closed invariant subset is called {\slsf a minimal set}
of $\call$.
\end{defi}

A minimal set is, obviously, the closure of any of its leaves. A
leaf $\call_m$ such that $\overline{\call_m}$ is minimal we shall
call {\slsf a minimal leaf}. Let us remark the following:

\begin{nnprop}
\label{ex-min} Let $I$ be a compact invariant set of a codimension
one holomorphic foliation $\call$ on a  complex manifold $X$. Then
there exists a minimal $\call$-invariant set $M$, which is contained
in $I$.
\end{nnprop}
\proof Consider the partially ordered by inclusion set $\cali$ of
closed invariant subsets of $I$. Let $A$ be its linearly ordered
subset. Let us see that there exists a smallest element in $A$.
Indeed, take $I_0\deff \bigcap_{i\in A}I_i$. It is a non empty
compact in $X$ and is obviously invariant unless $I_0\subset
\Sing\call$. We need to prove that $I_0\not\subset\Sing\call$.

\smallskip Let $U$ be an $(n-1)$-complete neighborhood  of $\Sing\call$ in
$X$, see \cite{Ba}. Any one of $I_i$ cannot be contained entirely
in $U$ by the maximum principle. Therefore $I_0\cap \big(X\setminus
U)$ is non empty. By Zorn's Lemma $\cali$ contains a minimal element
and this element is our minimal set.

\smallskip\qed

\newprg[prgPPD.pd]{Poincar\'e domains}
Let $X$ be a complex surface and let $\call$ be a (singular)
holomorphic foliation by curves on $X$.  We closely follow
constructions from \cite{Iy1,Iy2, Sz, Br5}. For a point $m\in
X^{\reg}\deff X\setminus \Sing\call$ denote by $\call_m$ the leaf of
$\call^{\reg}\deff \call|_{X^{\reg}}$ through $m$. Again, let us
underline that throughout this paper we are not considering any sort
of leaves through the singular points of $\call$. Take a smooth,
locally closed disc $D$ in $X^{\reg}$, transversal to $\call^{\reg}$
and cutting the leaf $\call_m$ at $m$. Therefore we require that $D$
is closed in some subdomain $Y$ of $X^{\reg}$ and that for every
point $z\in D$ the intersection $D\cap\call_z$ is transversal. We
shall call such discs - the {\slsf Poincar\'e discs}.

\smallskip Denote by $\hat\call_m$ the holonomy covering of
$\call_m$, \ie the covering with respect to the kernel of the
holonomy representation  $\hln : \pi_1 (\call_m,m)\to \dif (D,m)$.
The image $\hln (\call_m,m) \deff \hln \big(\pi_1(\call_m,m)\big)$
of the fundamental group of $\call_m$ under this representation is
called the {\slsf holonomy group} of $\call$ along the leaf
$\call_m$. It is a subgroup of the group $\dif (D,m)$ of germs of
biholomorphisms of the transversal $D$ fixing the point $m$. Up to a
conjugation  $\hln (\call_m,m)$ doesn't depends on the choice of the
transversal through $m$. We refer to \cite{Gd} for generalities on
holonomy groups of foliations. But in a moment we shall recall some
features of the holonomy representation, which will be crucial for
us along this paper. By $\tilde\call_m$ we denote the universal
covering of $\call_m$. Consider the following sets:

\begin{equation}
\eqqno(cylin-s) \calp_D \deff \bigcup_{z\in D}\call_z \qquad
\hat\calp_D
\deff \bigcup_{z\in D} \hat\call_z \qquad
\tilde\calp_D\deff\bigcup_{z\in D}\tilde\call_z.
\end{equation}
The first set we shall name the {\slsf Poincar\'e domain} of $\call$
over the transversal $D$.  The second is a {\slsf  holonomy covering
Poincar\'e domain} (``tube normaux'' in the terminology of
\cite{Sz}) and the third - the {\slsf universal covering Poincar\'e
domain} ("skew cylinder" in the terminology of \cite{Iy1,Iy2}). We
shall often call them simply holonomy or covering domains for short.
If no misunderstanding can occur, we shall also call $\calp_D ,
\hat\calp_D , \tilde\calp_D$ - the {\slsf Poincar\'e domains} of the
Poincar\'e disc $D$.

\smallskip Poincar\'e domain $\calp_D$ is an open connected subset
of $X^{\reg}$ and therefore of $X$. $D$ admits a tautological
imbedding $i:D\to \calp_D$, namely $i : z\to z$. Note that $i(D)$ is
not closed in $\calp_D$ in most interesting cases. Poincar\'e
domains $\hat\calp_D$ and $\tilde\calp_D$ come together with the
natural topologies and foliations on them. Let us briefly recall
this. Consider the topological space $\calc (D,\call)$ of pathes
$\gamma_{z,w}$ starting from points $z\in D$, ending at
$w\in\call_z$ and contained in $\call_z$, \ie we consider pathes
inside leaves of $\call$ only. The topology in $\calc (D,\call)$ is
the topology of uniform convergence on the space of continuous maps
from the unit interval $[0,1]$ to $X$. Pathes $\gamma_{z,w}$ and
$\beta_{z,w}$ (with the same ends) are equivalent if the holonomy
along $\beta_{z,w}^{-1}\circ \gamma_{z,w}$ is trivial. The holonomy
covering Poincar\'e domain is the quotient of $\calc (D,\call)$
under this equivalence relation.

\smallskip To check that this quotient is Hausdorff let us recall
what is the holonomy representation. Let $\gamma$ be a closed path
in $\call_m$, which starts and ends at $m$. If one ``displaces''
$\gamma$ to a nearby leaf $\call_z$, \ie if one takes a point $z\in
D$ close to $m$ and draws a path $\beta$ starting from $z$ in
$\call_z$ close to $\gamma$,  then $\beta$ certainly hits $D$, but
in general by a point $z'$ different from $z$. This way one obtains
a mapping $z\to z'$, which is called the holonomy representation of
$\gamma$. It depends only on the homotopy class of $\gamma$ and is
{\slsf holomorphic}. To see this one covers $\gamma$ by foliated
charts  and realizes that the holonomy representative map of
$\gamma$ is a composition of obvious holomorphic maps in these local
foliated charts, see \cite{Gd} for more details. The representation
obtained we denote as $\hln : \pi_1(\call_m,m)\to \dif (D,m)$. This
is a formalization of {\slsf d'application de premier retour} of
Poincar\'e.

\smallskip Now it is easy to see that $\hat\calp_D$ is Hausdorff.
Suppose that there exists a sequence $\gamma_{z_n,w_n}$, $z_n\to
z_0$ in $D$, $\call_{z_n}\ni w_n\to w_0\in \call_{z_0}$ in $\calc
(D,\call)$ which, after factorization, converges to the two limit
points $\gamma_{z_0,w_0}$ and $\beta_{z_0,w_0}$.  That means that
$\gamma_{z_n,w_n}$ uniformly converge to $\gamma_{z_0,w_0}$ and
there is another sequence $\beta_{z_n,w_n}$ which uniformly converge
to $\beta_{z_0,w_0}$, such that:

\smallskip\sli $\gamma_{z_n,w_n}$ and $\beta_{z_n,w_n}$ are equivalent
for all $n$,

\smallskip\slii while $\gamma_{z_0,w_0}$ and $\beta_{z_0,w_0}$ are not.

\smallskip In another words the holonomy along the closed path
$\beta_{z_0,w_0}^{-1}\circ \gamma_{z_0,w_0}$ is non trivial, but at the
same time $\hln (\beta_{z_0,w_0}^{-1}\circ \gamma_{z_0,w_0})(z_n)=z_n$.
We got a contradiction with the uniqueness theorem for holomorphic functions.

The natural map $\hat p(\gamma_{z,w}) = w$ is locally homeomorphic
and therefore the pair $(\hat\calp_D , \hat p)$ is a Riemann domain
over $X$. The map $i:D\to\hat\calp_D$, defined as $i : z \to
\gamma_{z,z}$, is a holomorphic imbedding and its image is a {\slsf
closed} disc in $\hat\calp_D$ - the {\slsf base} of the holonomy Poincar\'e
domain. $\hat\calp_D$ admits also a natural projection $\hat\pi$
onto $D$ defined as $\hat\pi (\hat\call_z) = z$. Holonomy Poincar\'e
domain $\hat\calp_D$ inherits a natural foliation $\hat\call$ with
leaves $\hat\call_z$ (the {\slsf holonomy foliation}) and the
locally biholomorphic map $\hat p:(\hat\calp_D,\hat\call)\to
(X,\call)$ is foliated, \ie sends leaves to leaves. Foliation
$\hat\call$ on $\hat\calp_D$ has no holonomy by construction.

\smallskip The same construction can be repeated with the following
equivalence relation: pathes $\gamma_{z,w}$ and $\beta_{z,w}$ are
equivalent if $\beta_{z,w}^{-1}\circ \gamma_{z,w}$ is homotopic to
the constant path $\gamma_{z,z}$ inside of the leaf $\call_z$. The
quotient (if it is Hausdorff!) is the universal covering Poincar\'e
domain $\tilde\calp_D$ in question. The corresponding objects are
marked as $\tilde p$, $\tilde\pi$, $i$ and $\tilde\call$ - the last
we shall call the {\slsf universal foliation} on $\tilde\calp_D$. We
shall need to see this construction starting from pathes in
$\hat\calp_D$. Consider the space of pathes $\calc
([0,1),\hat\call)$ starting from points in $D$ inside of leaves of
the holonomy foliation on the holonomy covering Poincar\'e domain
$\hat\call_D$. Two pathes $\gamma_{z,w}$ and $\beta_{z,w}$ are
equivalent if they are homotopic inside $\call_z$ to a constant path
$\gamma_{z,z}$. The quotient is the universal covering Poincar\'e
domain $\tilde\calp_D$. It possesses a natural projection
$p:\tilde\calp_D\to \hat\calp_D$ sending $\gamma_{z,w}$ to $w$.
Composition $\tilde p\deff \hat p\circ p :\tilde\calp_D\to \calp_D$
is the mapping, which one naturally obtains when constructing
$\tilde\calp_D$ starting from pathes in $\calp_D$.

\medskip   In general, it is useful to point out the following items:

\begin{itemize}
\item The universal covering Poincar\'e domain might  be non
Hausdorff.

\item The Poincar\'e domain  $\calp_D$ in most cases {\slsf cannot be
projected} to $D$, simply because the same leaf $\call_z$ may
intersect $D$ in several (even, in infinite number of) points.

\item  Both $\hat p : \hat\calp_D\to \calp_D$ and $\tilde p :
\tilde\calp_D\to \calp_D$ {\slsf are not regular coverings} in most
interesting cases (even if $\tilde\calp_D$ exists). They are regular
coverings only along the leaves.

\item  Both $(\hat\calp_D, \hat p)$ and $(\tilde\calp_D , \tilde p)$
are usually {\slsf not locally pseudoconvex} over $X$.

\end{itemize}

\newprg[prgPPD.d-g-f]{Docquier-Grauert criterion and Fujita's Theorems}
Recall that a domain $R$ in a complex manifold $X$ is called {\slsf
pseudoconvex} at its boundary point $z_0\in\d R$ if there exists a
Stein neighborhood $U\ni z$ such that $R\cap U$ is Stein in the
sense that each connected component of this intersection is a Stein
domain. $R$  is called {\slsf locally pseudoconvex} if it is
pseudoconvex at each of its boundary points.

\smallskip  Recall  that a Riemann domain over a complex manifold
$X$ is a pair $(R,p)$, which consists from a topological space $R$
and a locally homeomorphic map $p$ of this space into $X$. This
local homeomorphism induces an obvious complex structure on $R$, see
\cite{GR2}, and then $p$ becomes a  {\slsf local biholomorphism} of
complex manifolds $R$ and $X$. A Riemann domain $(R, p )$ over a
complex manifold $X$ is called {\slsf locally pseudoconvex over a
point } $z\in X$ if there exists a Stein neighborhood $U\ni z$ such
that all connected components of $p^{-1}(U)$ are Stein. If there
exists one such neighborhood, say $U\ni z$, then for every Stein
subdomain $V\subset U$ all connected components of $p^{-1}(V)$ will
be again Stein. Indeed, each component $V_1$ of $p^{-1}(V)$ is
pseudoconvex in some connected component $U_1$ of $p^{-1}(U)$.
Finally, $(R,p)$ is called locally pseudoconvex over $X$ if it is
locally pseudoconvex over every point of $X$.

\smallskip Recall that the Hartogs figure in $\cc^n$ is the following domain
\begin{equation}
\eqqno(hart-f) H^n_{\eps} = \Big(\Delta^{n-1}_{\eps}\times
\Delta\Big)\bigcup \Big(\Delta^{n-1}\times A_{1-\eps , 1}\Big).
\end{equation}
Here $\Delta_{\eps}$ stands for the disc of radius $\eps>0$ in
$\cc$, $\Delta\deff\Delta_1$, and $A_{1-\eps , 1}
\deff \Delta\setminus\bar\Delta_{1-\eps}$ is an annulus.

\smallskip Throughout this paper we shall repeatedly use the following
remarkable Docquier-Grauert criterion:

\begin{nnthm} {\slsf (Docquier-Grauert)}
\label{doc-gra} Let $(R,p)$ be a Riemann domain over a Stein
manifold $X$. If every holomorphic imbedding $h:H^n_{\eps}\to R$
extends to a locally biholomorphic mapping $\hat h:\Delta^n\to R$
then $R$ is a Stein manifold.
\end{nnthm}

In \cite{DG} this type of convexity of a domain over a Stein
manifold was called {\slsf $p_7$-convexity}. As an obvious corollary
from this criterion one gets one theorem of K. Stein: {\sl a regular
cover of a Stein manifold is Stein itself.} Remark that the inverse
is not true: think about $\cc^2$ covering a torus $\ttt^2$.

\smallskip We shall crucially use in our proofs the following results
of R. Fujita:

\begin{nnthm} {\slsf (Fujita)}
\label{fujita} \sli Let $(R,p)$ be a locally pseudoconvex Riemann
domain over $\pp^n$. If $p : R\to \pp^2$ is not a homeomorphism then
$R$ is a Stein manifold.

\smallskip\noindent\slii  Let $(R,p)$ be a locally pseudoconvex
Riemann domain over $\pp = \pp^{n_1}\times ...\times\pp^{n_k}$. If
$p$ is not a homeomorphism of $R$ onto a domain, which up to a
permutation, contains $\pp^{n_1} \times \{pt\}$, then $R$ is a Stein
manifold.
\end{nnthm}

For the proof see \cite{F1} and \cite{F2}. We shall repeatedly use the
following statement:

\begin{lem}
\label{intersect} Let $\call$ be a codimension one holomorphic
foliation in $\Delta^n$. Then every leaf of $\call$ intersects the
Hartogs figure $H^n_{\eps}$ (for any $\eps >0$).
\end{lem}
\proof Suppose that the leaf $\call_m$ is such that
$\overline{\call_m}\cap H^n_{\eps}=\emptyset$. Take any curve
$\gamma : [0,1]\in\Delta^{n-1}$, coming from $0$ to the boundary of
$\Delta^{n-1}$, and such that $\gamma ([0,1])$ doesn't intersects
the projection $S$ of the singularity set of $\call$ onto
$\Delta^{n-1}$ ($S$ is at most countable union of locally closed
hypersurfaces in $\Delta^{n-1}$). Consider the family of analytic
discs
\[
\Delta_{\gamma (t)} \deff \{\gamma (t)\}\times \Delta , \qquad t\in
[0,1].
\]
Remark that boundaries of these discs never intersect
$\overline{\call_m}$. Suppose that this family intersects
$\overline{\call_m}$. Let $t_0$ be the first value of $t$ such that
$\Delta_{\gamma (t)}\cap \overline{\call_m} \not= \emptyset$. $t_0$
exists because $\overline{\call_m}$ is closed and $t_0\not= 0$
because $\Delta_{\gamma (0)}\cap \overline{\call_m} = \emptyset$ by
assumption. We get a contradiction with the positivity of
intersections of complex varieties. Therefore $\overline{\call_m}$
should  be contained in $S\times \Delta$. But in that case it should
be an irreducible component of this set. And as such necessarily
intersects $H^n_{\eps}$. Contradiction.

\smallskip\qed

\newprg[prgPPD.loc-ps]{Local pseudoconvexity of Poincar\'e domains}
In what follows $\dim X =2$ if the opposite is not explicitly
stated. $\call$ is a singular holomorphic foliation by curves on
$X$. We start from the following:

\begin{lem}
\label{loc-ps1} Let $z_0$ be a non-isolated boundary point of
$\calp_D$. Then $\calp_D$ is pseudoconvex at $z_0$.
\end{lem}
\proof For a holomorphic foliation $\call$ on $X$ and an open,
connected subset $U\subset X$ we denote by $\call|_U$ the {\slsf
restriction} of $\call$ to $U$. The meaning is obvious, but let us
point out that the leaves of $\call|_D$ are {\slsf not} the
intersections of the leaves of $\call$ with $U$ in general. They are
the {\slsf connected components} of these intersections. To make a
clear distinction between leaves of $\call|_U$ and intersections
$\call_z\cap U$ we will set  $\calf\deff\call|_U$ and denote the
leaf of $\calf$ passing through $z\in U$ as $\calf_z$. Again, let us
repeat that all our leaves are defined outside of singular points.
Now consider two cases.

\medskip\noindent{\slsf Case 1. $z_0\not\in\Sing\call$.} Take an
$\call$-foliated neighborhood $U$ of $z_0$, \ie $U$ is biholomorphic
to $\Delta\times\Delta$ with  $\calf_{z_1}\deff \{z_1\}\times\Delta$
being the leaves of $\calf=\call|_U$. Since $\calp_D$ is an
$\call$-invariant domain the intersection $\calp_D\cap U$ is a union
of leaves of $\calf$, \ie has the form $U_1\times \Delta$ for some
open subset $U_1\subset \Delta$. Therefore $\calp_D\cap U$ is
pseudoconvex.

\begin{rema} \rm
\label{sm-n-sing} Remark that from the considerations, made above,
it follows that $z_0\in \d\calp_D\setminus \Sing\call$ cannot be an
isolated point of $\d\calp_D$. Indeed, if $z_0\in \d\calp_D$ then
$\call_{z_0}\subset \d\calp_D$.
\end{rema}

\medskip\noindent{\slsf Case 2. $z_0\in \Sing\call$.}  Let $U$ be a
neighborhood of $z_0$ biholomorphic to a ball, which doesn't
contains any other then $z_0$ singular point of $\call$. Take some
point $z\in \d\calp_D\cap U$ ($z$ may be well $z_0$), but in our
case, it is not an isolated point of $\d\calp_D$. Therefore we can
find a sequence $z_n\in \d\calp_D\cap U$, $z_n\not\in \Sing\call$,
converging to $z$. After going to a subsequence, $\calf_{z_n}$
converge in the Hausdorff metric to a closed in  $U$ set $L$. $L$
clearly contains $z$, is connected and $L\subset \d\calp_D\cap U$
because every point in  $L$ is a limit of boundary points, namely
points of $\calf_{z_n}$. Moreover, $L$ is $\calf$-invariant. Indeed,
a Hausdorff limit of closed invariant sets is, obviously, a closed
invariant set.

\smallskip We proved that $(\d\calp_D\cap U)\setminus
\{z_0\}$ is the union of leaves of $\calf=\call|_U$. Let
$h:H^2_{\eps}\to \calp_D\cap U$ be a holomorphic imbedding. Then $h$
induces a {\slsf smooth} holomorphic foliation $\cale\deff h^*\calf$
on $H^2_{\eps}$.

\begin{rema} \rm
\label{pull-fol} By pulling back a foliation $\calf$ from a manifold
$U$ by a (locally) biholomorphic map $h:\Omega\to U$ from a complex
manifold $\Omega$ into $U$ we mean the following. $\calf$ is defined
locally, on open sets $U_j$ by holomorphic forms $\omega_j$
satisfying the usual compatibility conditions. Take  pull backs
$h^*\omega_j$ of $\omega_j$. Then these holomorphic $1$-forms will
define a foliation in $\Omega$. It is this foliation we mean when
writing $h^*(\calf)$. The leaf of $\cale \deff h^*(\calf)$ through a
point $z\in\Omega\setminus \Sing\cale$ will be denoted by $\cale_z$.
\end{rema}
$\cale$ extends to a (singular) holomorphic foliation $\hat\cale$ on
$\Delta^2$. $h$ extends onto $\Delta^2$ and $h(\Delta^2)\subset U$
(by the usual Hartogs theorem, because it is a mapping into
$\cc^2$). Suppose that $h(\Delta^2)\cap\d\calp_D\not=\emptyset $.
Then $h(\Delta^2)$ contains some $z\in\d\calp_D\setminus \{z_0\}$,
again because $z_0$ is supposed to be a non isolated boundary point.
Let $w\in \Delta^2$ be some $h$-preimage of $z$. Then
$h(\hat\cale_w)\subset \calf_z\subset \d\calp_D$. But $\hat\cale_w$
intersects $H^2_{\eps}$, see Lemma \ref{intersect}. Therefore
$h(H^2_{\eps})$ should also intersect $\d\calp_D$. Contradiction. By
the Docquier-Grauert criterion we conclude that $\calp_D\cap U$ is
Stein.

\smallskip\qed

Since at its isolated boundary point a domain cannot be locally
pseudoconvex, we obtain the following:

\begin{nncorol}
\label{cor-ps1} The Poincar\'e domain $\calp_D$ is not locally
pseudoconvex at $z_0\in\d\calp_D$ if and only if $z_0$ is an
isolated point of $\d\calp_D$ and $z_0\in \Sing\call$.
\end{nncorol}

This implies immediately the following:

\begin{nncorol}
\label{cor-ps2}
Suppose that $\calm = \overline{\call_m}$ is an exceptional minimal
for a holomorphic foliation $\call$ in $\pp^2$. Then for every
Poincar\'e disc $D\ni m$ one has $\calp_D=\pp^2\setminus \Sing\call$.
\end{nncorol}
\proof If not then there exists at least one boundary point $z_1$ of
$\calp_D$, which is not an isolated point of $\d \calp_D$. Adding to
$\calp_D$ all its {\slsf isolated} boundary points we get a new
domain $\bar\calp_D$, which is pseudoconvex by Corollary
\ref{cor-ps1} and different from $\pp^2$ - it doesn't contains
$z_1$. Therefore $\bar\calp_D$ is Stein by Fujita's theorem and at
the same time it contains an invariant compact $\calm$.
Contradiction with the maximum principle.

\smallskip\qed

\smallskip To analyze the situation with the failure of pseudoconvexity
of Poincar\'e domains we shall employ the universal covering
Poincar\'e domains. However the main problem with $\tilde\calp_D$ is
that the natural topology on it might not be separable in general.
In \cite{Iy1,Iy2} and \cite{IS} the following statement was proved:

\begin{nnprop}
\label{cov-cyl-st} Let $\call$ be a holomorphic foliation by curves
on a Stein manifold $X$, let $\call_m$ be a leaf of $\call$ and $D$
a transversal hypersurface through $m$. Then $\tilde\calp_D$ is
Hausdorff. Moreover, if $D$ is Stein then $\tilde\calp_D$ is also
Stein.
\end{nnprop}

In the forthcoming Subsections we shall describe the obstructions to
the existence (\ie separability of the topology) of universal
coverings of Poincar\'e domains in compact complex surfaces. We
shall also replace the steiness of universal covering Poincar\'e
domains by another property appropriate for our needs - we call it
{\slsf rothsteiness}. But before doing that we shall need to aloud
our complex surfaces to have some mild singularities.

\newprg[prgPPD.models]{Nef models of holomorphic foliations}

Recall that a two-dimensional complex space $X$ has a {\slsf cyclic
quotient} singularity at point $a\in X$ if there exists a
neighborhood $U\ni a$ which is biholomorphic to the quotient
$\Chi^{l,d}\deff\Delta^2/\Gamma_{l,d}$. Here, for the relatively
prime $1\le l < d$ the group $\Gamma_{l,d}$ is defined by acting on
$\Delta^2$ as follows: $(z_1,z_2)\to (e^{\frac{2\pi i}{d}}z_1,
e^{\frac{2\pi il}{d}}z_2)$. Such neighborhood $\Chi^{l,d}$ carries a
natural foliation $\call^{\v}$ - we call it the {\slsf vertical
foliation} - it is defined as such that lifts to the {\slsf standard
vertical foliation} $\call^{\v}_{z_1} = \{z_1\}\times \Delta$ on
$\Delta^2$ under the natural {\slsf cyclic covering map}
$\pi_{l,d}:\Delta^2 \to \Chi^{l,d}$.

\smallskip Let $\call$ be a holomorphic foliation on a projective
surface $X$ with at most cyclic singularities. Our {\slsf standing
assumption} on $\call$ is that at every singular point of $X$ our
foliation $\call$ is biholomorphic to the vertical one. A leaf of
such $\call$ is still the leaf of $\call^{\reg}
\deff \call|_{X\setminus \Sing\call}$. In more colloquial terms that
means that we do not consider the cyclic points of $X$ as singular
points of $\call$ because $\call$, by our standing assumption, is as
good as a smooth foliation in a neighborhood of a cyclic point. That
means in its turn that a cyclic point $a$ does belongs to a certain
leaf of $\call$, one may note this leaf as $\call_a$ as we always
do. But we shall never consider a Poincar\'e disc through a cyclic
point. Finally, let us underline once more that singular points of
$\call$ do not belong to any of leaves of $\call$.

\smallskip Let us briefly discuss some specific features of leaves
passing through the cyclic points. Take the loop $\gamma =
\pi_{l,d}\big(\{0\}\times \{\frac{1}{2}e^{i\theta}: \theta\in
[0,2\pi)\}\big) \subset \Chi^{l,d}$ and take the Poincar\'e disc
$D=\pi_{l,d}(\Delta\times \{1\})$ at
$m=\pi_{l,d}((0,-\frac{1}{2}))$. Then the holonomy of $\call^{\v}$
along $\gamma$ is periodic with period $d$. Remark also that
$\Delta^2= \hat\calp_D^{\v}$ is the holonomy covering Poincar\'e
domain for $\calp^{\v}_D=\Chi^{l,d}$ and $\hat p^{\v}
:\hat\calp_D^{\v}\to\calp_D^{\v}$ is nothing but the natural cyclic
covering map $\pi_{l,d}$. Therefore we construct the holonomy
covering $\hat\calp_D$ for an arbitrary Poincar\'e disc $D$ situated
away from both $\Sing\call$ and $\Sing X$ literally in the same way
as we deed in Subsection \ref{prgPPD.pd}. Every cyclic point $a$
will contribute in a way that $\hat p :\hat\call_a\to \call_a$ will
be a ramified covering of order $d$ over $a$. More precisely, if the
leaf $\call_z$ is such that it passes through a cyclic point $a$
then $\hat p|_{\hat\call_z}:\hat\call_z \to \call_z$ is ramified
with the same order $d$ at {\slsf every} point over $a$. The mapping
$\hat p:\hat\calp_D\to\calp_D$ as a whole will be locally
biholomorphic over non-cyclic points and will behave as the standard
cyclic covering $\pi_{l,d}$ over the cyclic points.

\smallskip The universal covering Poincar\'e domain $\tilde\calp_D$
is now constructed from the holonomy covering domain $\hat\calp_D$
again as in Subsection \ref{prgPPD.pd}. The natural projection
$\tilde p$ will behave in the same way as $\hat p$ does and its
restriction $\tilde p|_z : \tilde\call_z \to \call_z$ to a leaf over
a cyclic point $z$ will have ramifications of the same
(corresponding to $z$) order $d$ at all points over $z$. I.e.,
$\tilde p|_{\tilde\call_z}: \tilde\call_z \to \call_z$ will be an
{\slsf orbifold covering}.

\begin{defi}
\label{hyp-leaf} A leaf $\call_z$ is called hyperbolic if
$\tilde\call_z$ is a disc. It is called parabolic if $\tilde\call_z$
is $\cc$ or $\pp^1$.
\end{defi}

From now on the hyperbolicity/parabolicity of leaves will be
understood in the sense defined above. If at least one leaf of
$\call$ has $\pp^1$ as its (orbifold) universal covering then
$\call$ is a rational quasi-fibration. We are usually excluding this
exceptional (and trivial) case from our statements.

\begin{rema} \rm
\label{non-rat1} A holomorphic foliation $\call$ on a projective
surface $X$ with at most cyclic singularities is called a rational
quasi-fibration if the closure of its generic fiber is a rational
curve. For such foliations our theorem is obvious, therefore in the
sequel we shall suppose that $\call$ is not a rational
quasi-fibration whenever it will be convenient for us.
\end{rema}

\begin{rema} \rm
\label{ps-conv-rem} Results on pseudoconvexity of Poincar\'e domains
were proved in this Section for foliations on {\slsf smooth}
surfaces only, and only in that special case they will be used along
this paper.
\end{rema}

Recall furthermore the following:

\begin{defi}
\label{nef-fol} Let $\call$ be a holomorphic foliation on a
projective surface $X$ with at most cyclic singularities such that
$\Sing\call \cap \Sing X = \emptyset$, and let $\calk_{\call}$
denotes its canonical bundle. $\call$ is called {\slsf nef}
(numerically effective) if $\calk_{\call}$ is nef, \ie if for any
irreducible algebraic curve $C\subset X$ one has $\calk_{\call}\cdot
C \ge 0$.
\end{defi}

Our guideline in the proof of the main theorem of this paper will be
the following statement, which is one of the main results of
McQuillan's ''Noncommutative Mori Theory'' of holomorphic
foliations:

\begin{nnthm} {\slsf (McQuillan)}
\label{mcq} Let $\call$ be a holomorphic foliation by curves on an
algebraic surface $X$ with at most cyclic singularities, which is
not a rational quasi-fibration. Then there exists a bimeromorphic
transformation of $X$ onto a projective surface $Y$ with at most
cyclic singularities, such that the transformed foliation $\calf$ is
{\slsf nef} and therefore enjoys the following alternative:

\smallskip\sli either all leaves of $\calf$ are parabolic;

\smallskip\slii or, the set $\Sing\call \cup \{\text{ parabolic leaves }\}$
is a proper algebraic subset $\cala$ of $X$ and, moreover, the
Lobatchevski distance is continuous on $X\setminus \cala$.
\end{nnthm}

Of course, the condition $\Sing\calf\cap \Sing X=\emptyset$ is also
preserved. The bimeromorphic transformation from $X$ to $Y$ consists
of two steps: first - blowing up singular points of $\call$ to make
all singularities {\slsf reduced}; second - contracting invariant
rational curves violating the nef condition. In the process of the
proof of Theorem \ref{p2-inv-set} we shall go through these two
steps paying a specific attention to what is happening with
exceptional minimal leaves of $\call$ and with the Poincar\'e
domains containing them.

\begin{rema} \rm
\label{non-rat2} {\bf 1.} It should be remarked that a rational
quasi-fibrations cannot be brought to the nef model.

\smallskip\noindent{\bf 2.} The reference for Theorem \ref{mcq} is
\cite{McQ1} and, more specifically \cite{McQ2}. A more analytically
motivated reader might find helpfull consulting the nice expositions
in \cite{Br1,Br2}.
\end{rema}

\newprg[prgPPD.obstr]{Obstructions to the existence of universal covering
Poincar\'e domains}

We shall need the existence of universal covering Poincar\'e domains
in more general cases then foliations on Stein manifolds.

A foliated {\slsf holomorphic} immersion between foliated pairs
$(X,\call)$ and $(Y,\calf)$ is a holomorphic immersion $f:X\to Y$
which sends leaves to leaves. A foliated {\slsf meromorphic}
immersion is a meromorphic map which is a holomorphic foliated
immersion outside of its indeterminacy set. All immersions,
considered in this paper will be mappings between manifolds of the
same dimension, \ie locally biholomorphic over the smooth points,
and they will be aloud to behave as standard cyclic covering over
the cyclic points.

\smallskip Let $\Delta^2=\Delta \times\Delta$ be a bidisc in
$\cc^2$. Recall that the foliation $\call^{\v}$ in $\Delta^2$ with
leaves $\call^{\v}_{\lambda}\deff \{\lambda\}\times \Delta$,
$\lambda\in \Delta $ we called {\slsf vertical}.

\begin{defi}
\label{str-fol} A foliated pair $(X,\call)$, with $\dim X = 2$, will
be called {\slsf straight} if any foliated meromorphic immersion
$h:(\Delta^2,\call^{\v})\to (X,\call)$  is, in fact, holomorphic.
\end{defi}

The class of straight foliated pairs we shall denote  as
$\cals_{\calt}$. An important for us observation is due to
Brunella, see Lemma 1 in \cite{Br5}.

\begin{lem}
\label{brun} Let $(X,\call)$ be a {\slsf nef} foliated pair and let
$p:(\Delta^2,\call^{\v}) \to (X,\call)$ be a foliated meromorphic
immersion. Then $p$ is holomorphic.
\end{lem}

In another words a nef foliated pair is straight. Let us emphasis
that in both Definition \ref{str-fol} and Lemma \ref{brun} the
foliated pair $(X,\call)$ is aloud to have cyclic singularities.

\smallskip Another class of straight foliated pairs on {\slsf smooth}
manifolds appears as follows (and will be needed later in Section
\ref{sect.RAT}). Let us denote by $\calb_{\calh}$ the class of
connected, Hausdorff, countable at infinity complex manifolds $X$
such that every locally biholomorphic map $h:\Omega\to X$ from a
Stein domain $\Omega$ of dimension $n=\dim X$ extends to a locally
biholomorphic map of the envelope of holomorphy $\hat\Omega$ into
$X$. Note that vacuously all $1$-dimensional complex manifolds are
in $\calb_{\calh}$ simply because in dimension one all maps are
holomorphic (meromorphic functions are also holomorphic as mappings
with values in $\pp^1$). All Stein manifolds belong to
$\calb_{\calh}$. Also, if $\tilde X$ is a cover of $X$ then $X$ and
$\tilde X$ belong, or not to $\calb_{\calh}$ simultaneously. In
particular, all tori are in $\calb_{\calh}$. By a {\slsf covering}
between {\slsf smooth} manifolds we mean a local
homeomorphism/biholomorphism $\pi : \tilde X\to X$ such that
$\pi^{-1}$ is extendable along all pathes in $X$ for all initial
values. More generally to $\calb_{\calh}$ do belong all complex
manifolds without rational curves etc.

\smallskip In the future we shall need the following, less obvious,
observation:

\begin{nnprop}
\label{lb-ext} \sli Complex projective space $\pp^n$ belongs to
$\calb_{\calh}$ for every $n\ge 1$;

\smallskip\slii If $X_1 , X_2\in \calb_{\calh}$ then $X_1\times X_2\in
\calb_{\calh}$;

\smallskip\sliii in particular, if $S_1,...,S_n$ are Riemann surfaces
then $S_1\times ... \times S_n\in \calb_{\calh}$.
\end{nnprop}
For (\sli see \cite{Iv2}, for (\slii and (\sliii - \cite{Iv3}. Of
course, only the holomorphic extendability is a problem here. If a
map extends holomorphically (and not {\slsf meromorphically (!)}, as
one should expect in general) then this  extension will be
necessarily also locally biholomorphic.

\smallskip We have in our disposal the following three subclasses in
$\cals_{\calt}$.

\begin{itemize}
\item If $(X,\call)$ is nef then $(X,\call)$ is straight.

\smallskip
 \item $X\in\calb_{\calh}$ then $(X,\call)\in\cals_{\calt}$ for
 {\slsf any} $\call$.

\smallskip
\item If a foliated pair $(X,\call)$ doesn't contains invariant
rational curves then $(X,\call)\in\cals_{\calt}$.
\end{itemize}
About the second item  it should be said that there exist compact
complex manifolds $X\not\in\calb_{\calh}$ carrying holomorphic
foliations by curves $\call$ such that $(X,\call)\in\cals_{\calt}$,
see Example 5.1 in \cite{Iv6}.

\begin{nncorol}
\label{nef-cov} Let $(X,\call)$ be a straight foliated pair and
$\calp_D$ a Poincar\'e domain. Then the universal covering
Poincar\'e domain is Hausdorff. In particular, this is always true:

\smallskip\slii for a {\slsf nef} foliated pair;

\smallskip\slii for {\slsf any} $(X,\call)$ with $X=\pp^2$,
$\pp^1\times \pp^1$ and all from Proposition \ref{lb-ext}.
\end{nncorol}
\proof We follow \cite{Br3}. Obstructions to the separability of the
topology on $\tilde\calp_D$ are {\slsf vanishing cycles}, see
Section 3 in \cite{Iv6}. A vanishing cycle, which appears this way,
can be supposed to be {\slsf imbedded}, see Lemma 3.4 in \cite{Iv6}.
It means that there exists an imbedded loop $\gamma$ in
$\hat\call_m$ (starting and ending at $m$) such that:

\begin{itemize}

\item $\gamma$ is not bounding a disc in $\hat\call_m$;

\smallskip\item But there exist a sequence $m_k\in D$, $m_k\to m$
and imbedded loops $\gamma_k \subset \hat\call_{m_k}$ which
uniformly converge to $\gamma$ and such that every $\gamma_k$ bounds
a disk $L_k$ in $\hat\call_{m_k}$.
\end{itemize}

Perturbing slightly, if necessary, we can include this sequence of
loops into a continuous family $\Gamma = \{\gamma_z: z\in D\}$ (one
might need to take a smaller $D$). Then one constructs a
``generalized Hartogs figure'' $W$ around $\Gamma$, see Section 2 of
\cite{Iv6} for more details. $W$ is a special foliated subdomain in
$\hat\calp_D$.  Then one replaces $W$ by another foliated domain
$(V,\pi)$ over the same disc $D$ and transfers the map $\hat
p|_W:W\to X$ to a foliated map $q : V\to X$. The domain $V$, in
contrast to $W$, is foliated over the disc $D$ with all fibers being
discs. One more feature of $V$ is that the fiber $V_z$ of $V$ over
$z\in D$ is mapped by $q$ into the leaf $\call_z$ (with the same
$z$).

\smallskip The map $q$ extends to $V$, but the extended map
$\hat q$ might become {\slsf meromorphic}, \ie it might have a
discrete set $R$ of points of indeterminacy in $V$. This clearly
doesn't happens in the straight case, because $\hat q$ is foliated
and locally biholomorphic outside of its (eventual) points of
indeterminacy.  Since $V_m$ is mapped into $\call_m$ by $\hat q$ we
see that $\gamma$ also bounds a disc, namely $\hat p^{-1}(\hat
q(V_m))$, \ie $\gamma$ is not a vanishing cycle.

\smallskip\qed

\newprg[prgPPD.roth]{A Rothstein type extension theorem}

We shall need a non standard version of the Rothstein's extension
theorem. Recall that the classical Rothstein's theorem  states the
following, see \cite{Si2}:

\begin{nnthm}{\slsf (Rothstein)}
\label{rothstein1} Let $f$ be a holomorphic/meromorphic function in
the unit bidisc $\Delta^2 = \Delta\times \Delta$. Suppose that for
every $z_1\in \Delta$ the restriction $f_{z_1}\deff f(z_1,\cdot)$
extends as a holomorphic/meromorphic function of the variable $z_2$
onto the disc $\Delta_R$ with $R>1$. Then $f$
holomorphically/meromorphically extends onto the bidisc $\Delta\times
\Delta_R$ as a function of both variables $(z_1,z_2)$.
\end{nnthm}

Let us call a complex manifold (or, a normal complex  space) $X$ a
{\slsf Rothstein manifold} (space) if the statement of Rothstein's
theorem is valid for {\slsf holomorphic} mappings with values in
$X$. Stein manifolds are obviously Rothstein. If $\tilde X\to X$ is
a covering then $\tilde X$ and $X$ are Rothstein, or not
simultaneously. We send the interested reader to Subsection 2.5 in
\cite{Iv5} for more information when the Rothstein-type theorem is
valid. Then compare with Lemma 6 from \cite{Iv4} to derive that this
property is invariant under the regular coverings (the property of
being Stein is not invariant).

\smallskip In this paper we are motivated by a more general statement,
 which basically says that all complex manifolds are ``almost Rothstein''.

\begin{nnprop}
\label{rothstein2} Let $X$ be a complex manifold (or, a normal
complex space). Then $X$ is "almost Rothstein" in the sense that
every holomorphic/meromorphic mapping $f:\Delta^2\to X$, such that
for all $z_1\in \Delta$ the restriction $f_{z_1}\deff f(z_1,\cdot)$
holomorphically/meromorphically extends onto $\Delta_R$, $R>1$, is
holomorphically/meromorphically extendable onto $(\Delta\times
\Delta_R)\setminus (E\times\Delta_R)$, where $E$ is a closed polar
subset of $\Delta$.
\end{nnprop}
For the proof  see Corollary 2.5.1 in \cite{Iv5}.

\begin{defi}
\label{cov-roth-def} A foliated pair $(\tilde\calp,\tilde\call)$,
where $X$ is a (not necessarily compact) {\slsf smooth} complex
surface, we shall call {\slsf Rothstein} if for every foliated
holomorphic immersion $f:(\Delta^2,\call^{\v}) \to (\tilde\calp,
\tilde\call)$ such that for every $z_1\in \Delta$ the restriction
$f_{z_1}$ extends to a holomorphic immersion onto $\Delta_R, R>1$,
the map $f$ extends to a holomorphic immersion onto
$\Delta\times\Delta_R$ as a mapping of both variables $(z_1,z_2)$.
\end{defi}

\smallskip By a Riemann domain over a complex surface $Y$
with at most cyclic singularities we understand a {\slsf smooth}
complex surface $R$ together with a holomorphic mapping $p:R\to Y$
such that:

\smallskip\sli $p$ is locally biholomorphic over non cyclic points;

\smallskip\slii for every cyclic point $b\in Y$ and every $a\in p^{-1}(b)$
there exist neighborhoods $W\ni a$ and $V\ni b$ such that the
restriction $p|_W : W\to V$ is the standard cyclic covering.

\smallskip Let us state now a variation of a Rothstein theorem,
which will be needed in this paper. Let $\tilde\calp_D$ be a universal
covering Poincar\'e domain of a holomorphic foliation on a projective
surface with at most cyclic singularities and let $\tilde p :
(\tilde\calp_D,\tilde\call) \to (X,\call)$ be the canonical foliated
projection, \ie $(\tilde\calp_D,\tilde p)$ is a Riemann domain over
$X$ (provided that $\tilde\calp_D$ is Hausdorff).

\begin{nnprop}
\label{rothstein4} Let $(X,\call)$ be a straight foliated pair,
where $X$ is a projective surface with at most cyclic singularities
and let $\calp_D$ be a Poincar\'e domain in $(X,\call)$. Then the
universal covering foliated pair  $(\tilde\calp_D,\tilde\call)$ is
Rothstein.
\end{nnprop}
\proof First of all $(\tilde\calp_D,\tilde\call)$ exists as a
Hausdorff topological space by Corollary \ref{nef-cov}. Let
$f:(\Delta^2,\call^{\v})\to (\tilde\calp_D,\tilde\call)$ be a
foliated immersion such that for every $z_1\in \Delta$ the
restriction $f_{z_1}$ extends as a holomorphic immersion onto
$\Delta_R$. The composition $g\deff \tilde p \circ f$ is meromorphic
on $\Delta\times \Delta_R$ by the classical theorem of Rothstein
\ref{rothstein1}, because we assumed $X$ to be projective. Let us
check that $g$ is, moreover, a meromorphic immersion. If not then
let $C$ be the critical set of $g$ in $\Delta\times \Delta_R$. $C$
cannot contain a component of the form $\{z_1\}\times \Delta_R$
because $g$ is a holomorphic immersion on $\Delta\times\Delta$.
Therefore the intersection of $C$ with $(\Delta\setminus E)\times
\Delta_R$ it non empty. Here $E$ stands for the polar set of
Proposition \ref{rothstein2} for $f$. Let $c_1\not\in E$ be such
that $(\{c_1\}\times \Delta_R)\cap C$ contains some $c_2$. Since
$f_{z_1}$ for all $z_1$ is supposed to be an immersion then $f$
could fail to be an immersion in a neighborhood of $c=(c_1,c_2)$
only if $f$ contracts to a point some local component $C_1$ of $C$,
which passes through $c$. But this is impossible because the
universal foliation $\tilde\call$ on $\tilde\calp_D$ has no
singularities. Therefore $f$ is an immersion in a neighborhood of
$c$ and so must be also $g=\tilde p\circ f$. Contradiction.

\smallskip Therefore $g$ is a meromorphic immersion. By assumed
straightness of $(X,\call)$ our $g$ is holomorphic everywhere. Take
a point $a=(a_1,a_2)\in \Delta\times \d\Delta$. Let $V$ be a cyclic
neighborhood of $b\deff g(a)$. Taking a sufficiently small
neighborhood $U$ of $a$ of the form $\Delta_r(a_1)\times
\Delta_{\rho}(a_2)$ and an  appropriate coordinates $(v_1,v_2)$  in
$V\ni b$, we can suppose that the mapping $g|_U : U\to V$ has the
standard cyclic form. Let $\{W_j\}$ be the at most countable set of
all connected components of $\tilde p^{-1}(V)$. By connectivity and
the fact that $f$ is a foliated holomorphic (on both variables)
immersion on $\Delta_r(a_1) \times (\Delta_{\rho}(a_2)\cap \Delta )$
we see readily that there exists such $j_0$ that for all $z_1\in
\Delta_r(a_1)$ we have that $f_{z_1}(\Delta_{\rho}(a_2))\subset
W_{j_0}$, \ie that the restriction of $f$ to $\Delta_r(a_1)\times
\Delta_{\rho}(a_2)$ takes its values in $W_{j_0}$. Moreover, this
restriction is jointly holomorphic on $\Delta_r(a_1)\times
(\Delta_{\rho}(a_2)\cap \Delta)$.

\smallskip The disc $f_{a_1}(\Delta_{\rho}(a_2))$ standardly covers
the disc $g_{a_1}(\Delta_{\rho}(a_2))$: like $z\to z^d$ (everywhere
in this text $d=1$ is not excluded). Therefore shrinking both
$W_{j_0}$ and $V$ we obtain a cyclic covering $\tilde p|_W:W\to V$
such that for all $z_1\in \Delta_r(a_1)$ (with some smaller $r>0$)
$f_{z_1}(\Delta_{\rho} (a_2))\subset W$. This is because
$f_{z_1}(\Delta_{\rho} (a_2)\cap \Delta)\subset W$. But $W$ is a
bidisc. Therefore our $f$ extends as a holomorphic map
of two variables onto $\Delta_r(a_1)\times \Delta_{\rho}(a_2)$ by
Rothstein's theorem. The rest is obvious. Remark only that the
``vertical size`` $\rho$ in our construction depends only on
``vertical size'' of the cyclic neighborhood $V$ of $b$ (and on
$g$), but not on $f$.

\smallskip\qed

\newsect[sectAB]{Holomorphic representation of the fundamental group}

We have in mind a certain "unification" of two representations of
the fundamental group $\pi_1(\call_m,m)$ of a leaf of a holomorphic
foliation $\call$. The first is the standard one - the {\slsf
holonomy representation}, it was briefly recalled in Subsection
\ref{prgPPD.pd}, the second is the representation by the deck
transformations of the universal covering $\tilde\call_m\to\call_m$.
Both are one dimensional representations in the sense that the space
of these representations is either the {\slsf one dimensional}
complex disc $D$ - \ie the Poincar\'e disc, or the Riemann surface
$\tilde\call_m$ itself.

\smallskip In this Section we shall construct one more holomorphic
representation of $\pi_1(\call_m,m)$. It will be full dimensional
and will act by foliated biholomorphisms on the universal (and
holonomy) covering Poincar\'e domains. This "unified" representation
will be our principal tool in proving the main results  of this
paper.

\smallskip Moreover, using the hyperbolic feature of the holonomy
group in or setting we shall prove that under these circumstances we
can {\slsf expand} the universal (and holonomy also) covering
Poincar\'e domain to some ``Poincar\'e domain`` $\tilde\calp_{\cc}$
{\slsf foliated over} $\cc$ and it will (in the case of hyperbolic
$\call$) {\slsf regularly} cover $\calp_D$. In the case of $X=\pp^2$
it will lead to a contradiction and this will finish the proof. The
case of parabolic $\call$ will be excluded with a different argument
(but also using the expanded domain $\tilde\calp_{\cc}$).

\newprg[prgAB.rep]{A germ of the  holomorphic representation}
Let $\call$ be a singular  holomorphic foliation by curves on a
compact complex surface $X$ with at most cyclic singularities. Fix
some Riemannian metric $r$ on $X$. Fix some point $m\in X^{0}\deff
X\setminus (\Sing\call \cup \Sing X)$ and let $\call_m$ be the leaf
of $\call^{\reg}\deff \call|_{X^{\reg}}$ through $m$. Here $X^{\reg}
\deff X\setminus \Sing\call$ and may contain cyclic points. Take a
small disc $m\in D\subset X^0$, transversal to $\call^{\reg}$. By
saying ``a small disc'' we mean that $D$ is a disc of a small
geodesic radius with center in $m$. Of course, the transversality to
$\call$ will be always supposed. By saying ``a smaller'' subdisc of
$D$ we mean a subdisc of smaller geodesic radius and with the same
center $m$. In this context writing $D_{k}\subset D$ we mean that
$D_{k}$ has radius $1/k$. Our discs will be {\slsf always} situated
in $X^0$.

\smallskip Take a point $w\in \tilde\call_m$ such that $\tilde p (w)
= m$. Denote by $[\gamma]$ the element of $\pi_1(\call_m,m)$ which
realizes $w$. Take a foliated neighborhood $U\ni m$ in
$\tilde\calp_D$ and let $U_0\ni m$ be its biholomorphic image in $X$
under the canonical foliated projection $\tilde p$. Take a foliated
neighborhood $V\ni w$ in $\tilde\calp_D$ such that $\tilde p|_V:V\to
U_0$ is a foliated biholomorphism (to achieve this one might need to
shrink $U_0$ and therefore $U$).

\smallskip The composition  $\tilde\pi \circ \tilde p|_V^{-1}\circ
\tilde p|_D : (D,m)\to (D,m)$ defines the image $[\gamma ]^{hln}$ of
$[\gamma ] \in \pi_1(\call_m,m)$ in the holonomy group $\hln
(\call_m,m)\subset \dif (D,m)$. Define a germ of a foliated
imbedding as follows:
\begin{equation}
\eqqno(phi-gam1) \phi_{\gamma }\deff \tilde p|_V^{-1}\circ \tilde
p|_U : (\tilde\calp_D,m) \to (\tilde\calp_D,w).
\end{equation}
First we shall extend it to a ''more global`` germ. Namely, denote
by $\bihol (\tilde\calp_D,\tilde\call_m)$ the group of foliated
biholomorphic imbeddings of foliated neighborhoods of
$\tilde\call_m$ into $\tilde\calp_D$, which send $\tilde\call_m$
onto itself. Therefore an element $\phi$ of $\bihol
(\tilde\calp_D,\tilde\call_m)$ is a biholomorphic foliated mapping
$\phi : \tilde\calp_V\to \tilde\calp_W$, where $V$ and $W$ are some
neighborhoods of $m$ in $D$, such that  $\phi (\tilde\call_m) =
\tilde\call_m$. Here for $V\subset D$ we set $\tilde\calp_V =
\tilde\pi^{-1}(V)$ and call it a foliated neighborhood of
$\tilde\call_m$ when $V$ is a neighborhood of $m$.

\begin{lem}
\label{hol-exp1} Suppose that $\tilde\calp_D$ is Rothstein. Then the
germ $\phi_{\gamma }$ extends to a foliated imbedding of a foliated
neighborhood of $\tilde\call_m$ into $\tilde\calp_D$ such that
$\phi_{\gamma }(\tilde\call_m) = \tilde\call_m$. It depends only on
the homotopy class $[\gamma ]$  of $\gamma$. The map
\[
\tilde\Phi:\pi_1 (\call_m,m)\to \bihol (\tilde\calp_D,\tilde\call_m)
\]
\begin{equation}
\eqqno(hat-phi1) \gamma  \mapsto \phi_{\gamma }
\end{equation}
is a monomorphism of groups.
\end{lem}
\proof  Take a path $\tilde\gamma $ in $\tilde\call_m$ from $m$ to
the point $w$, which defines our germ $\phi_{\gamma}$ as in
\eqqref(phi-gam1). Let $U, U_0$ and $V$ be as above. Denote by $D_0$
a sufficiently small subdisc in $D$ such that $U$ is foliated over
$D_0$ by $\tilde\pi|_U:U\to D_0$. Restrictions $\tilde
p|_V^{-1}\circ\tilde p_z$ of our germ onto the leaves, initially
defined on $\tilde\call_z\cap U$ with values in
$\tilde\call_{[\gamma]^{hln}(z)}\cap V$, extend along any path in
$\tilde\call_z$ which starts at $z$. This extensions clearly gave us
a singlevalued biholomorphic maps of $\tilde\call_z$ onto
$\tilde\call_{[\gamma]^{hln}(z)}$ - deck transformations of the universal
coverings of leaves in question. Since $\tilde\calp_D$ is
Rothstein these extensions glue together for $z\in D_0$ to a
holomorphic foliated imbedding $\phi_{\gamma}:\tilde\pi^{-1}(D_0)\to
\tilde\calp_D$. This imbedding send $\tilde\call_m$ onto itself,
because $[\gamma]^{hln}(m) =m$, \ie is a germ of a foliated
biholomorphism in a foliated neighborhood of $\tilde\call_m$.

\smallskip The fact that $\tilde\Phi$ is a homomorphism of groups is
obvious, because $\phi_{\gamma}$ restricted to $\tilde\call_m$ is
nothing but the deck transformation of the universal covering
$\tilde p_z : \tilde\call_m\to \call_m$. From this point of view our
representation is an extension of the deck transformation group to a
neighborhood of $\tilde\call_m$ in $\tilde\calp_D$. The extension
$\phi_{\gamma}$ of every deck transformation is subdued to the
condition that it is $\tilde p$ - equivariant by the construction.
From this remark it becomes obvious that $\phi_{\gamma}$ is uniquely
determined by its restriction $\phi_{\gamma}|_{\tilde\call_m}$. This
proves that our representation is a monomorphism of groups.

\smallskip\qed

\begin{defi}
\label{hol-rep-def} Monomorphism $\tilde\Phi : \pi_1 (\call_m,m)\to
\bihol (\tilde\calp_D, \tilde\call_m)$ we shall call the {\slsf
holomorphic representation} of the fundamental group of the leaf
$\call_m$.
\end{defi}

It is also can be called the {\slsf holomorphic extension of the
deck transformation group} and is more precise then the holonomy
representation (which is not a monomorphism in general).

\begin{rema} \rm
Perhaps the most comprehensive view  on $\tilde\Phi$ is that it is a
"unification" of the two "orthogonal" representations of the
fundamental group of the leaf $\call_m$: one is the holonomy
representation, the second - is the representation by the deck
transformations of the universal covering $\tilde\call_m\to
\call_m$.
\end{rema}

Up to now our exposition was quite general. We newer used any
specific features of the holonomy group that will appear in the
following Subsection.

\newprg[prgAB.expan]{Expansion of the holomorphic representation}

Now we are going to explore  the fact that in our case $\hln
(\call_m,m)$  contains a hyperbolic element. Denote this hyperbolic
element as $\alpha$ and its  holomorphic representative
$\tilde\alpha$, \ie $\tilde\alpha \deff \tilde\Phi (\alpha)$. Now
let us fix a coordinate $t\in \Delta$ such that $\alpha$ becomes to
be a multiplication by the complex number $\alpha$, $0<|\alpha|<1$
in this coordinate. Rescaling $t$, if necessary, we can suppose that
$D = D_{|\alpha|^{-2}}$, where by $D_r$ we denote the subdisc $\{ t:
|t| < r\}$ in $\cc$. Set $A_{r_1, r_2}
\deff D_{r_1}\setminus \bar D_{r_2}$ - the annulus of radii
$0<r_1<r_2$. Fix some $0<\eps < 1 - |\alpha|$. For every integer
$n\ge 0$ consider the foliated domain
\begin{equation}
(\tilde\calp_{A_{|\alpha|(1-\eps), |\alpha|^{-2}}}, \alpha^{-2n}
\tilde\pi) \qquad\text{ over the annulus } \qquad
A_{|\alpha|^{-2n+1}(1-\eps),|\alpha|^{-2n-2}}.
\end{equation}
Mappings $\alpha^k\tilde\pi :\tilde\calp_{A_{|\alpha|(1-\eps),
|\alpha|^{-2}}} \to A_{|\alpha|^{k+1}(1-\eps),|\alpha|^{k-2}}$ for
$k\in\zz$  and $w\in\tilde\calp_{A_{|\alpha|(1-\eps),
|\alpha|^{-2}}}$ are defined as $(\alpha^k\tilde\pi)(w) =
\alpha^k\tilde\pi (w)$. The domain
$(\tilde\calp_{A_{|\alpha|(1-\eps), |\alpha|^{-2}}},
\alpha^0\tilde\pi)$ we consider as identical to
$\tilde\calp_{A_{|\alpha|(1-\eps), |\alpha|^{-2}}}\subset
\tilde\calp_{D}$. Glue domain $(\tilde\calp_{A_{|\alpha|(1-\eps),
|\alpha|^{-2}}},\alpha^{-2n}\tilde\pi)$ to
$(\tilde\calp_{A_{|\alpha|(1-\eps),
|\alpha|^{-2}}},\alpha^{-2n-2}\tilde\pi)$ by the biholomorphism
$\tilde\alpha^2 : \tilde\calp_{A_{|\alpha|^{-1}(1-\eps
),|\alpha|^{-2}}} \to \tilde\calp_{A_{|\alpha|(1-\eps), 1}}$, see
diagram \eqqref(glue) below.

\[
\def\normalbaselines{\baselineskip20pt\lineskip3pt \lineskiplimit3pt }
\def\mapright#1{\smash{\mathop{\longrightarrow}\limits^{#1}}}
\def\mapleft#1{\smash{\mathop{\longleftarrow}\limits^{#1}}}
\def\mapdown#1{\Big\downarrow\rlap{$\vcenter{\hbox{$\scriptstyle#1$}}$}}
\begin{matrix}
&\tilde\calp_{A_{|\alpha|(1-\eps), |\alpha|^{-2}}}&\supset&
\tilde\calp_{A_{|\alpha|^{-1}(1-\eps ),|\alpha|^{-2}}}&
\mapright{\tilde\alpha^2}&\tilde\calp_{A_{|\alpha|(1-\eps), 1}}\cr &
\mapdown{\alpha^{-2n}\tilde\pi}& &\mapdown{\alpha^{-2n}\tilde\pi}&
&\mapdown{\alpha^{-2n-2}\tilde\pi}& \cr
&A_{|\alpha|^{-2n+1}(1-\eps),|\alpha|^{-2n-2}}&\supset&
A_{|\alpha|^{-2n-1}(1-\eps),|\alpha|^{-2n-2}}&\mapright{\id}&
A_{|\alpha|^{-2n-1}(1-\eps),|\alpha|^{-2n-2}}\cr
\end{matrix}
\]

\smallskip
\begin{equation}
\eqqno(glue)
\def\normalbaselines{\baselineskip20pt\lineskip3pt \lineskiplimit3pt }
\def\mapright#1{\smash{\mathop{\longrightarrow}\limits^{#1}}}
\def\mapleft#1{\smash{\mathop{\longleftarrow}\limits^{#1}}}
\def\mapdown#1{\Big\downarrow\rlap{$\vcenter{\hbox{$\scriptstyle#1$}}$}}
\begin{matrix}
&\equiv&\tilde\calp_{A_{|\alpha|(1-\eps),
1}}&\subset&\tilde\calp_{A_{|\alpha|(1-\eps), |\alpha|^{-2}}}\cr &
&\mapdown{\alpha^{-2n-2}\tilde\pi}& &
\mapdown{\alpha^{-2n-2}\tilde\pi}\cr &\equiv&
A_{|\alpha|^{-2n-1}(1-\eps),|\alpha|^{-2n-2}}&\subset&A_{|\alpha|^{-2n-1}
(1-\eps),|\alpha|^{-2n-4}}. \cr
\end{matrix}
\end{equation}

The obtained Poincar\'e domain over $\cc$ denote as
$\tilde\calp_{\cc}$ and call it the {\slsf expanded universal
covering Poincar\'e domain}. The projection $\tilde\pi$ obviously
extends to a holomorphic projection $\tilde\pi :\tilde\calp_{\cc}
\to \cc$ by construction: the extended $\tilde\pi$ on each
$(\tilde\calp_{A_{|\alpha|(1-\eps), |\alpha|^{-2}}},
\alpha^{-2n}\tilde\pi)$ is simply equal to $\alpha^{-2n}\tilde\pi$.

\begin{lem}
\label{expand1} The germ $\tilde\alpha$ extends to a global foliated
biholomorphism of the expanded Poincar\'e domain, which commutes
with $\tilde\pi$. The canonical foliated projection $\tilde
p:\tilde\calp_D\to X$ also extends to a foliated immersion $\tilde
p:\tilde\calp_{\cc} \to \calp_D$ and this extension stays to be
$\tilde\alpha$-equivariant.
\end{lem}
\proof The proof of extendability of $\tilde\alpha$ consists in
checking of the correctness of its natural definitions on the
overlapping subsets. Let us do it for $\tilde\alpha^{-1}$ instead of
$\tilde\alpha$. Subdomain $\big(\tilde\calp_{A_{|\alpha|^{-1}(1-\eps
), |\alpha|^{-1}}},\alpha^{-2n}\tilde\pi\big)$ is identified with
the domain $\big(\tilde\calp_{A_{|\alpha|(1-\eps),|\alpha|}},
\alpha^{-2n-2}\tilde\pi\big)$ by $\tilde\alpha^2$, see
\eqqref(glue). Therefore $w \in
\big(\tilde\calp_{A_{|\alpha|^{-1}(1-\eps ), |\alpha|^{-1}}},
\alpha^{-2n}\tilde\pi\big)$ is identified with $\tilde\alpha^2 w\in
\big(\tilde\call_{ A_{|\alpha|(1-\eps),|\alpha|}},
\alpha^{-2n-2}\tilde\pi\big)$, see the upper horizontal line in the
diagram \eqqref(ext-alpha). One can act by $\tilde\alpha^{-1}$ both
on $w$ and on its twin $\tilde\alpha^2 w$, see down arrows in the
diagram below.

\begin{equation}
\def\normalbaselines{\baselineskip20pt\lineskip3pt \lineskiplimit3pt }
\def\mapright#1{\smash{\mathop{\longrightarrow}\limits^{#1}}}
\def\mapleft#1{\smash{\mathop{\longleftarrow}\limits^{#1}}}
\def\mapdown#1{\Big\downarrow\rlap{$\vcenter{\hbox{$\scriptstyle#1$}}$}}
\begin{matrix}
\big(\tilde\calp_{A_{|\alpha|^{-1}(1-\eps ),
|\alpha|^{-1}}},\alpha^{-2n} \tilde\pi\big)\ni
w&\longleftrightarrow& \tilde\alpha^2 w\in
\big(\tilde\calp_{A_{|\alpha|(1-\eps),|\alpha|}},
\alpha^{-2n-2}\tilde\pi\big)\cr \mapdown{\tilde\alpha^{-1}}&
&\mapdown{\tilde\alpha^{-1}}\cr
\big(\tilde\calp_{A_{|\alpha|^{-2}(1-\eps), |\alpha|^{-2}}},
\alpha^{-2n}\tilde\pi\big)\ni
\tilde\alpha^{-1}w&\longleftrightarrow& \tilde\alpha
w\in\big(\tilde\calp_{A_{1-\eps, 1}},\alpha^{-2n-2}\tilde\pi\big).
\cr
\end{matrix}
\eqqno(ext-alpha)
\end{equation}

On the left one gets $\tilde\alpha^{-1}w$, on the right
$\tilde\alpha w$. But they are identified by $\tilde\alpha$.
Therefore $\tilde\alpha^{-1}$ is correctly defined globally. It is
commuting with $\tilde\pi$ by construction.

\smallskip Since the gluing maps involved in the expansion of the universal
covering Poincar\'e domain are $\tilde p$-equivariant the map $\tilde p$
extends obviously to a locally biholomorphic foliated map $\tilde
p:\tilde\calp_{\cc}\to \calp_D$ and in its turn stays to be
$\tilde\alpha$-equivariant. Moreover, we have that for every $t\in\cc$
the restriction $\tilde p|_t$ is the universal covering of
$\call_{\alpha^Nt}$ for $N$ big enough, namely $\alpha^Nt$ should be
in $D$.  As we already remarked the image of $\tilde\calp_{\cc}$
under $ \tilde p $ is nothing but $\calp_D$ due to the periodicity
of the expanded domain.

\smallskip\qed

\newprg[prgAB.hyperb]{Universal covering Poincar\'e domains of hyperbolic
foliations}

As we explained in the Introduction the proof of the Theorem
\ref{p2-inv-set} will be done separately for two different cases:
when $\call$ is parabolic and when it is hyperbolic. In this
Subsection we consider the hyperbolic case.

\smallskip Therefore let $\call$ be a hyperbolic holomorphic foliation on
a projective surface $X$ with at most cyclic singularities, \ie at
least one leaf of $\call$ is hyperbolic. In that case one can
locally define the following hyperbolic norm on  vectors tangent to
$\call$:
\begin{equation}
\eqqno(lob) \norm{\v}_l = \inf \big\{1/r: \exists \text{ holomorphic
} u:\Delta (r) \to \call_z, u(0) = z, u^{'}(0) = \v\big\}.
\end{equation}
If the leaf $\call_z$ passes through a cyclic point then one defines
the hyperbolic norm on the (orbifold) universal covering $\tilde\call_z$ of 
$\call_z$ and pushes it down to a singular (at cyclic point) metric on 
$\call_z$. In fact we shall nod need to push it down and will work on 
$\tilde\call_z$ (more precisely on $\tilde\calp_D$).
One calls (depending on traditions) $\norm{\v}_l$ the  Lobatchevski/
Poincar\'e/Kobayashi/ hyperbolic norm of $\v\in T\call_z$. Function
\eqqref(lob) is well defined on tangent vectors to leaves, $L(\v)
\deff \ln \norm{\v}_l$ is finite if $\v$ is tangent to a hyperbolic
leaf, and is equal to $-\infty$ is it is tangent to a parabolic
leaf. We say that Lobatchevski metric is continuous if for any local
holomorphic vector filed $\v$ tangent to $\call$ the local function
$\norm{\v}_l$ is continuous.

\smallskip If we suppose that a sufficiently small transversal disc
through $m$ doesn't cuts any parabolic leaf of $\call$, \ie that all
leaves $\calp_{z}$ for $z\in D$ are hyperbolic, then, by
construction, the same holds for all leaves of $\tilde\calp_{\cc}$.
As one  see from the McQuillan's alternative we can suppose that the
Lobatchevski norm is {\slsf continuous} on $\tilde\calp_{\cc}$
(provided $(X,\call)$ is {\slsf nef}).

\smallskip Our aim in this Subsection is to prove that $\tilde p:
\tilde\calp_{\cc} \to \calp_D$ is a regular covering in the
hyperbolic case.

\begin{lem}
\label{global-hyp} Suppose that $\call_m$ has a  hyperbolic holonomy
and that all leaves $\call_z$ for $z\in D$ are hyperbolic. Moreover,
suppose that that $\tilde\calp_{\cc}$ is Hausdorff and Rothstein 
and that the Lobatchevski metric is continuous on $\tilde\calp_{\cc}$. 
If for some $a_1\in \tilde\call_{z_1}$ and $b_1\in\tilde\call_{w_1}$ 
one has that $\tilde p(a_1) = \tilde p(b_1)$ then there exists a global 
$\tilde p$-equivariant automorphism $\phi$ of $\tilde\calp_{\cc}$ such 
that $\phi (a_1) = b_1$.
\end{lem}
\proof Applying  $\tilde\alpha^N$ to $a_1$ and $b_1$ with $N$ big
enough we can suppose that corresponding $z_1,w_1$ belong to $D$.
Take foliated neighborhoods $U\ni a_1$, $V\ni b_1$ and $U_0\ni
\tilde p(a_1)$ such that $\tilde p|_U:U\to U_0$ and $\tilde
p|_V:V\to U_0$ are foliated biholomorphisms. As in the proof of
Lemma \ref{hol-exp1} the composition $\tilde p|_V^{-1}\circ \tilde
p|_U$ extends along the leaves to a foliated biholomorphism $\phi$
of some $\tilde\calp_{D_0}$ and $\tilde\calp_{B_0}$, where $D_0\ni
z_1$ and $B_0\ni w_1$. Move $a_1$ to $z_1$ inside
$\tilde\call_{z_1}$ and follow this move by the move of $b_1$ inside
$\tilde\call_{w_1}$ in order to still have $\tilde p(a_1) = \tilde
p(b_1)$ with $a_1=z_1$ this time.

\smallskip Let us prove that $\phi$ extends along any path in $\cc$
in the sense that for any path $\gamma : [0,1]\to \cc$, $\gamma (0)
= z_1$, there exists a continuous family of discs $D_t$ with centers
at $\gamma (t)$ of radii $r(t)$ (continuously depending on $t$), and
there exist foliated $\tilde p$-invariant biholomorphisms $\phi_t :
\tilde\calp_{D_t}\to \tilde\calp_{B_t}$ (for appropriate domains
$B_t$) such that $\phi_{t_1}$ coincide with $\phi_{t_2}$ on
$\tilde\calp_{D_{t_1}\cap D_{t_2}}$ for $t_1-t_2$ small enough. Of
course we mean that $\phi_0 = \phi$.

\smallskip Let $t_0$ be the supremum of $t$-s such that $\phi$ extends
up to $t$. All we need to prove is that $\phi$ extends also to a
neighborhood of $t_0$. Set $\beta (t) \deff \tilde\pi ( \phi (\gamma
(t)))$. Let us prove first that $\beta (t)$ stays in a compact part
of $\cc$ as $t\nearrow t_0$. If not then there exists a sequence
$0<t_1<t_2<...<t_{2n-1}<t_{2n} < ... \to t_0$ such that $|\beta
(t_{2n-1})|= |\alpha|^{-2k_n+1}(1-\eps)$, $|\beta (t_{2n})|=
|\alpha|^{-2k_n-2}$ and $\beta ([t_{2n-1},t_{2n}])\subset
A_{|\alpha|^{-2k_n+1}(1-\eps),|\alpha|^{-2k_n-2}}$. Remark that
applying $\tilde\alpha^N$ with $N$ big enough (from the very
beginning) we can suppose that $\gamma ([t_1,t_0])\subset A_{|\alpha
| (1-\eps),|\alpha|^{-2}}$.

\smallskip Applying $\tilde\alpha^{2k_n}$ to every piece $\phi (\gamma
([t_{2n-1},t_{2n}]))$ and taking a subsequence we obtain that
$\tilde\pi(\tilde\alpha^{2k_n}(\phi (\gamma ([t_{2n-1},t_{2n}]))))$
converge in Hausdorff sense to a continuum $C\subset A_{|\alpha
|(1-\eps),|\alpha|^{-2}}$. Let us prove now the following:

\smallskip\noindent{\slsf Claim. $\call_{\gamma(t_0)}$ coincides
with $\call_c$ for every $c\in C$.} Indeed, take a sequence
$\tau_n\nearrow t_0$ such that $c_n\deff\tilde\pi
(\tilde\alpha^{k_n}(\phi (\gamma (\tau_n))))\to c$. Denote by $r_n$
the geodesic distance from $\tilde p|_{c_n}(c_n)$ to $\tilde
p_{\gamma (\tau_n)}(\gamma (\tau_n))$ in the leaf $\call_{\gamma
(\tau_n)}$. The sequence $\{r_n\}$ is obviously bounded. This
results from the fact that $\tilde p$ is defined and holomorphic in
a neighborhood of both of $t_0$ and $c$, and the fact that the
hyperbolic distance along the leaves is continuous in $\calp_D$. But
that means that we can find a point $b_n\in \tilde\call_{c_n}$ on a
distance not more then $r_n$ such that $\tilde p_{c_n}(b_n) = \tilde
p_{\gamma (t_n)} (\gamma (t_n))$. After taking a subsequence $b_n\to
b_0\in \tilde\call_c$. But then $\tilde p_{c}(b_0) = \tilde
p_{\gamma (t_0)} (\gamma (t_0))$ by continuity of $\tilde p$ and the
Claim is proved.

\smallskip But this (what we get in the Claim) is
impossible because one leaf can cut $D$ only by at most countable
set. Therefore $\beta (t)$ stays bounded when $t\to t_0$.

\smallskip Applying $\tilde\alpha^N$ once more we can suppose that
$\beta (t)$ stays in $D_{1/2}$. Limiting set of $\beta (t)$ when
$t\to t_0$ can be either a point or a continuum. But the latter is
impossible by the reason already explained above. Therefore $\beta
(t)\to w_0$ for some $w_0\in D$ when $t\to t_0$. Mapping $\phi$
along $\gamma$ writes as $\phi = \tilde p^{-1}\circ\tilde p$ for
some choice of $\tilde p^{-1}$.

\smallskip Both $\tilde p|_{\tilde\call_{w_0}}$ and $\tilde
p|_{\tilde\call_{\gamma (t_0)}}$ cover that same leaf $\call_{\gamma
(t_0)}$ (see the Claim). Let $d_l(\cdot)$ be the Lobatchevski
distance along the leaves of $\tilde\call$. Remark that $d_l(\tilde
p|^{-1}_{\beta (t)}\circ \tilde p|_{\gamma (t)},\beta (t))$ is
continuous up to $t_0$. Indeed, it is nothing else but the distance
from $\tilde p(\gamma (t))$ to $\tilde p(\beta (t))$ and the latter
is continuous up to $t_0$. Therefore the limiting set of $\phi
(\gamma (t))$ when $t\nearrow t_0$ is a compact $K$ in
$\tilde\call_{w_0}$. This implies that biholomorphisms
$\phi_{\tilde\call_{\gamma (t)}}: \tilde\call_{\gamma (t)} \to
\tilde\call_{\beta (t)}$ converge to a biholomorphism
$\tilde\call_{\gamma (t_0)}\to \tilde\call_{w_0}$ as $t\nearrow
t_0$. Now we can easily extend $\tilde p^{-1}|_{\tilde\call_{w_0}}
\circ \tilde p|_{\tilde\call_{\gamma (t_0)}}$ to foliated
neighborhoods as it was done at the beginning of the proof of the
Theorem.

\smallskip\qed

Denote by $G$ the group of all $\tilde p$-equivariant foliated
biholomorphisms of $\tilde\calp_{\cc}$. We have the following:

\begin{nnthm}
\label{prop-disc-hyp} If all leaves of $\calp_D$ are hyperbolic and
the hyperbolic distance is continuous on $\calp_D$ then $\tilde
p:\tilde\calp_{\cc}\to \calp_D$ is a regular covering.
\end{nnthm}
\proof Let us underline that by saying that $\tilde
p:\tilde\calp_{\cc}\to \calp_D$ is a regular covering we mean that
$\calp_D = \tilde\calp_{\cc}/G$, in particular, it is the standard
cyclic covering over the cyclic points. Recall, that an action of a
discrete group $G$ on a complex manifold $\tilde\calp_{\cc}$ is
called proper discontinuous if for every compacts $K_1, K_2\comp
\tilde\calp_{\cc}$ the set
\[
\{g\in G: gK_1\cap K_2\not=\emptyset \}
\]
is finite. But in our case we have a $G$-equivariant local
biholomorphism $\tilde p : \tilde\calp_{\cc}\to X$.  Suppose there
exist $w_n\to w_0$ in $K_1$ such that $g_n(w_n) = v_n \to v_0$ in
$K_2$, and $g_n\in G$ are distinct. Take neighborhoods $W\ni w_0$
and $V\ni v_0$ such that $(\tilde p|_V)^{-1}\circ\tilde p|_W: W\to
V$ is a biholomorphism. Since $\tilde p$ is $G$-invariant we get that
for $n>>1$ $g_n|_W = (\tilde p|_V)^{-1}\circ \tilde p|_W$, \ie all
$g_n$ for $n>>1$, are equal to each other. Contradiction.

\smallskip Action of $G$ is cyclic, \ie every point $a\in \tilde\calp_{\cc}$
admits a neighborhood $U$ such that $G_a \deff \{g\in G: gU\cap U
\not= \emptyset\}$ is isomorphic to $\Gamma_{l,d}$ for the
appropriate $1\le l <d$. Indeed, take $N>>1$ in order to have
$b\deff\alpha^Na\in D_{\frac{1}{2}}$. It is obviously sufficient to
find a needed neighborhood $U$ of $b$. Take $U$ such that $\tilde
p|_U:U\to X$ is cyclic.  Since $\tilde p$ is $G$ - equivariant and
the last acts by global biholomorphisms one cannot have any other
behavior of $G$  at $a$, because then $\tilde p$ would not be cyclic
on $U$.

\smallskip\qed

\newprg[prgAB.imbed-c2]{Imbedding of the expanded Poincar\'e
domain into $\cc^2$.}

Before considering the case of parabolic foliations we need few
preparatory lemmas.

\begin{defi}
\label{fol-domain} By a foliated domain we shall mean a triple
$(W,\pi,D)$ where $D$ is a domain in $\cc$, $W$ is a connected
complex surface, and $\pi:W\to\cc$ is a holomorphic submersion with
connected fibers.
\end{defi}

A holomorphic section of a foliated domain $(W,\pi,D)$ is a
holomorphic map $\sigma : D\to W$ such that $\pi\circ\sigma = \id$.
Remark that our covering Poincar\'e domains $\tilde\calp_D$  do
admit holomorphic sections, namely the natural map $i:D\to
\tilde\calp_D$ defined as $i : z \to \gamma_{z,z}$ is such a
section. The complex linear space $\cc^2$ we shall also see as a
foliated domain $(\cc^2,\pi_1)$, where $\pi_1 :(z_1,z_2)\to z_1$ is
the natural ''vertical`` projection. By a foliated holomorphic
imbedding of $(W,\pi,D)$ into $(\cc^2,\pi_1)$ we mean an imbedding
$H:W\to \cc^2$ such that every leaf $W_z$ is mapped into the leaf
$\cc_z\deff \{z\}\times\cc$, \ie  $H$ has the form
\begin{equation}
\eqqno(fol-imb) (z,\cdot) \to (z, h(z,\cdot)),
\end{equation}
where for every fixed $z$ the function $h(z,\cdot)$ realizes a
conformal imbedding of the domain $W_z$ into $\cc_z$.

\smallskip The following result is due to Brunella, see \cite{Br4}:

\begin{lem}
\label{imb-c2-1} Let $\tilde\calp_D$ be a universal covering Poincar\'e
domain over a simply connected transversal $D$ such that:

\smallskip\sli all fibers $\tilde\call_z$ for $z\in D$ are biholomorphic
to $\cc$;

\smallskip\slii the foliated pair $(X,\call)$ is straight.

\smallskip Then there exists a foliated holomorphic imbedding
$H:\tilde\calp_D\to (\cc^2,\pi_1)$ sending a leaf $\tilde\call_z$ to
$\{z\}\times \cc$ for all $z\in D$.
\end{lem}

\begin{rema} \rm
The statement of \cite{Br4} is more general, but $\tilde\calp_D$ is
understood there also in a more general sense. Namely, one should
add to some leaves of $\tilde\call$ certain "vanishing ends" in
order for $\tilde\calp_D$ to become Hausdorff. But in the straight
case one doesn't needs to do that according to Lemma \ref{nef-cov}.
\end{rema}

Remark that if at least one leaf of $\tilde\calp_D$ is $\pp^1$ then
all are such. Therefore we exclude this case in what follows.

\begin{nncorol}
\label{parab-v-c2} If all leaves of $\tilde\calp_{\cc}$ are
parabolic and different form $\pp^1$ then
$(\tilde\calp_{\cc},\tilde\pi)$ is leafwise biholomorphic to
$(\cc^2,\pi_1)$.
\end{nncorol}
\proof First let us check that $(\tilde\calp_{\cc},\pi, \cc)$
satisfies the Gromov's {\slsf spray condition}, see \cite{Gro} or
\cite{Fo} for definitions. For that take a point $z\in \cc$ and such
a small disc $B\ni z$ that $\pi : \tilde\calp_B\to B$ admits a
section over $B$. Then one can apply Lemma \ref{imb-c2-1}  and imbed
the restriction $(\tilde\calp|_B,\pi , B)$ leafwise into $\cc^2$.
After that the spray condition becomes obvious. This condition via
the theorem of Gromov, see the same sources, provides us a {\slsf
global} section of our fibration.

\smallskip After that one can employ the result of Siu that
every Stein submanifold (a section in our case) admits a Stein
neighborhood, which is, in addition biholomorphic to a neighborhood
of the zero section in the normal bundle to this submanifold, see
Corollary 1 in \cite{Si5}.

\smallskip Therefore the proof of Brunella, given in \cite{Br4},
applies again and finishes the proof of our Corollary.

\smallskip\qed

\newprg[prgAB.parab]{Universal covering Poincar\'e domains of parabolic
foliations} We now ready to consider the case when $\call$ is
parabolic, \ie every leaf $\call_z$ of $\call$ has as its (orbifold)
universal covering either $\pp^1$ or $\cc$. The former case means
that $\call$ is a rational quasi-fibration and therefore we don't
need to consider it. The cases $\overline{\call_m} = \pp^1, \ttt^1$
are also obvious, in the setting of Theorem \ref{p2-inv-set},
but we shall include them for the sake of completeness.

\smallskip Since our coverings $\tilde p_z:\tilde\call_z\to \call_z$
are, in fact, orbifold coverings it is useful to recall the formula 
relating their Euler characteristics:

\begin{equation}
\label{euler}
\chi (\tilde\call_z,\nu) = \chi (\call_z) + \sum_j\left( \frac{1}{\nu (r_j)} - 1\right),
\end{equation}
where $\chi (\call_z)$ is the Euler characteristic of the underlying Riemann surface 
$\call_z$ and $\nu (r_j)$ is the value of ramification function $\nu $ at ramification 
point $r_j$, see  \cite{Mil} and references there. Parabolicity of $\tilde\call_z$ means
that $\chi (\tilde\call_z,\nu)\ge 0$. This leaves very few possibilities for $\call_z$
and the ramification function $\nu$. When $\call_z$ is noncompact we have only the
following ones:

\smallskip{\bf r)} $\call_z = \cc$ with one ramification point, or with two ramification
points of order two.

\smallskip{\bf n)}  $\call_z=\cc^*$ without ramifications.

\begin{nnthm}
\label{parab-case} Let $\call_m$ be a parabolic leaf with the
hyperbolic holonomy of a holomorphic foliation $\call$ on a
projective surface $X$ with at most cyclic singularities. Suppose
that for a sufficiently small Poincar\'e disc $D$ all leaves cutting
$D$ are also parabolic and that $\tilde\calp_D$ is Hausdorff and
Rothstein. Then:

\smallskip\sli either $\overline{\call_m}$ is a rational or elliptic curve;

\smallskip\slii or, $\call_m = \cc^*$ and is an imbedded analytic
subset in some open subset of $X\setminus\Sing\call$.
\end{nnthm}
\proof Let $\overline{\call_m}$ be not algebraic and therefore such
are all $\overline{\call_z}$ for $z$ in a sufficiently small
Poincar\'e disc $D$. Note that $\call_m$ cannot be $\cc$ because
$\pi_1(\call_m)$ contains an element $\alpha$ with hyperbolic
holonomy. Therefore $\call_m = \cc^*$ and the universal covering
$\tilde\call_m\to \call_m$ is unramified. 

\smallskip According to Lemma \ref{imb-c2-1} in this case
$(\tilde\calp_{\cc}, \tilde\pi)$ is leafwise biholomorphic to
$(\cc^2,\pi_1)$. Up to an affine change of the coordinate $z_2$ we
can suppose that the holomorphic representative $\tilde\alpha$ of
$\alpha$ acts as follows:
\begin{equation}
\eqqno(action1)
\begin{cases}
\tilde\alpha (z_1,z_2) = (\alpha z_1, h(z_1,z_2)),\cr 
\tilde\alpha (0,z_2) = (0,z_2+1),
\end{cases}
\end{equation}
because $\tilde p_0:\tilde\call_0 \to \call_0$ is an unramified  
covering of $\cc^*$ by $\cc$. Let us prove the following:

\begin{lem}
\label{action-lem} There exist $\eps>0$ and $0<C<\infty$ such that
the automorphism $\tilde\alpha$ satisfies for all $n$ and all
$(z_1,z_2)$ with $|z_1|<\eps$ the relation
\begin{equation}
\eqqno(action2) \tilde\alpha^n (z_1,z_2) = \big(\alpha^n z_1, z_2+n
+ a_n(z_1) + b_n(z_1)z_2\big),
\end{equation}
where
\begin{equation}
\eqqno(action3) |a_n(z_1)|\le \frac{C|z_1|}{(1-|\alpha|)^2}
\qquad\text{ and } \qquad |b_n(z_1)|\le \frac{C|z_1|}{1-|\alpha|}.
\end{equation}
Moreover, $a_n\to a_0$, $b_n\to b_0$ uniformly on $z_1$ in a
neighborhood of zero. Here $a_0$ and $b_0$ are holomorphic in a
neighborhood of zero.
\end{lem}
\proof Since for every fixed $z_1$ the function $h(z_1,z_2)$ should
be an automorphism of $\cc$ and for $z_1=0$ it has the form
as in \eqqref(action1) we obtain that $h(z_1,z_2) = z_2 + 1 +
a_1(z_1) + b_1(z_1)z_2$ with $a_1(0) = 0$ and $b_1(0)=0$. Set
$A=a_1$, $B=b_1$ and write
\begin{equation}
\eqqno(action4)
\begin{cases}
\tilde\alpha (z_1,z_2) = \big(\alpha z_1, z_2 + 1 + A(z_1) +
B(z_1)z_2 \big), \cr \text{where } |A(z_1)|, |B(z_1)|\le C|z_1|
\text{ for } |z_1|<\eps .
\end{cases}
\end{equation}
Here $0<C<\infty$ and  $\eps >0$ are some constants. Let us prove by
induction that for every $n\in \nn$ one has
\begin{equation}
\eqqno(action5)
\begin{cases}
\tilde\alpha^n (z_1,z_2) = \big(\alpha^n z_1, z_2+n +
\sum_{k=1}^n\big[kA_{n,k}(z_1) + z_2B_{n,k}(z_1)\big]\big), \cr
\text{with } |A_{n,k}(z_1)|, |B_{n,k}(z_1)|\le
C|\alpha|^{k-1}\prod_{j=1}^{n-1} \big[1+C|\alpha|^j|z_1|\big]|z_1|
\text{ for } 1\le k\le n.
\end{cases}
\end{equation}
For $n=1$ \eqqref(action5) is nothing but \eqqref(action4). Next
write, using \eqqref(action4) and \eqqref(action5) the second
component $\tilde\alpha^{n+1}_2$ of $\tilde\alpha^{n+1}$ as follows:
\[
\tilde\alpha^{n+1}_2(z_1,z_2) = z_2+n + 1 +
\sum_{k=1}^n\big[kA_{n,k}(z_1) + z_2B_{n,k}(z_1)\big] +
A(\alpha^nz_1) + \big(z_2+n + \sum_{k=1}^n\big[kA_{n,k}(z_1) +
\]
\[
z_2B_{n,k}(z_1)\big]\big)B(\alpha^nz_1)= z_2+n + 1 +
\sum_{k=1}^nkA_{n,k}(z_1) \big[1 + B(\alpha^nz_1)\big] +
A(\alpha^nz_1) + nB(\alpha^nz_1) +
\]
\[
 +  z_2\sum_{k=1}^nB_{n,k}(z_1)\big[1 + B(\alpha^nz_1)\big] +
 z_2B(\alpha^nz_1).
\]
Set $A_{n+1,n+1}(z_1) \deff \frac{1}{n+1}A(\alpha^nz_1) +
\frac{n}{n+1}B(\alpha^nz_1)$. The estimate $|A_{n+1,n+1}(z_1)|\le
C|\alpha|^n|z_1|$ follows. Analogously set $B_{n+1,n+1}(z_1) \deff
B(\alpha^nz_1)$ and get $|B_{n+1,n+1}(z_1)|\le C|\alpha|^n|z_1|$.

\smallskip For $1\le k\le n$ set $A_{n+1,k}(z_1) \deff A_{n,k}(z_1)
\big[1+B(\alpha^nz_1)\big]$. Then by induction we have
\[
|A_{n+1,k}(z_1)| = |A_{n,k}(z_1)\big[1+B(\alpha^nz_1)\big]|\le
C|\alpha|^{k-1}\prod_{j=1}^{n-1}
\big[1+C|\alpha|^j|z_1|\big]|z_1|(1+C|\alpha|^n|z_1|),
\]
which gives us \eqqref(action5). Further set $B_{n+1,k}(z_1) \deff
B_{n,k}(z_1)\big[1 + B(\alpha^nz_1)\big]$ and get the same estimate.
$\eqqref(action5)$ is proved.

\smallskip Since $\big|\prod_{j=1}^{n-1}\big[1+C|\alpha|^j|z_1|\big]\big|\le K$
we get the estimates:
\begin{equation}
\eqqno(action6) |A_{n,k}(z_1)|, |B_{n,k}(z_1)|\le
C|\alpha|^{k-1}|z_1|.
\end{equation}
These estimates plus the usual summation in \eqqref(action5) gave
the proof of the lemma.

\smallskip\qed

\smallskip Now we can get more information about the global behavior
of $\tilde p$. Suppose that for some $z_1\not=0$ the restriction
$\tilde p_{z_1}: \cc_{z_1}\to X$ has a nontrivial period, \ie that
there exists a non-zero complex number $a(z_1)$ such that
$\tilde p|_{z_1}$ is invariant under the translation $z_2\to z_2 +
a(z_1)$ on the complex line $\cc_{z_1}$.

\begin{lem}
\label{scroling} If for some $z_1\not=0$ the restriction $\tilde
p_{z_1}$ has a nontrivial period, then there exists $\eps >0$ and a
non-vanishing holomorphic function $a$ in $\Delta (0,\eps)$ such
that $a(0)=1$ and
\begin{equation}
\eqqno(scrol) \tilde p(z_1,z_2) \equiv \tilde p(z_1, z_2 + a(z_1))
\end{equation}
for all $(z_1,z_2)\in \Delta (0,\eps)\times \cc$.
\end{lem}
\proof Denote by $\call_{z_1}$ the leaf which is covered by $\tilde
p_{z_1}$. Since $\call_{z_1}$ is also covered by every $\tilde
p_{\alpha^nz_1}$ we can suppose that $|z_1|<\eps$ where $\eps>0$ is
from Lemma \ref{action-lem}. From $\tilde\alpha$-invariance of
$\tilde p$ we see that for every $z_2$ the point
\[
\tilde\alpha (z_1,z_2 + a(z_1)) = (\alpha z_1, z_2 + a(z_1) + 1 +
a_1(z_1) + b_1(z_1)a(z_1) + b_1(z_1)z_2)
\]
should be a translation by some $d_1\cdot a(\alpha z_1)$  of the
point
\[
\tilde\alpha (z_1,z_2) = (\alpha z_1, z_2 + 1 + a_1(z_1) + b_1(z_1)
z_2)
\]
on the line $\cc_{\alpha z_1}$. Here $a(\alpha z_1)$ is a notation
for the period of this translation. Therefore
\[
a(z_1) + a_1(z_1) + b_1(z_1)a(z_1) + b_1(z_1)z_2 = a_1(z_1) +
b_1(z_1) z_2 + d_1\cdot a(\alpha z_1).
\]
From here we get that
\begin{equation}
\eqqno(a-alpha1) a(\alpha z_1) = \frac{1}{d_1}[1 + b_1(z_1)]a(z_1).
\end{equation}
Likewise,  using the formula \eqqref(action2), we get that
\begin{equation}
\eqqno(a-alpha2) a(\alpha^n z_1) = \frac{1}{d_n}[1 +
b_n(z_1)]a(z_1).
\end{equation}
Here, again, $a(\alpha^nz_1)$ is a notation for the period of the
corresponding translation in $\cc_{\alpha^nz_1}$.
\begin{rema} \rm
\label{periods} Periods of $\tilde p_{\alpha^nz_1}$ (we always mean
{\slsf minimal} periods) are uniquely defined, because 
$p_{\alpha^nz_1}$ can be supposed to be noncompact.
\end{rema}
Recall that $b_n\to b_0$. Would $d_n$ be non bounded, some
subsequence $a(\alpha^{n_k} z_1)$ would converge to zero. This
contradicts to the local biholomorphicity of $\tilde p$ at the 
origin. Therefore $[1 + b_0(z_1)]a(z_1) = \lim\limits_{k\to \infty}
[1+b_{n_k}(z_1)]a(z_1)$ is a (may be, non minimal) period of $\tilde
p_{0}$, and it is a limit of (may be not minimal) periods of $\tilde
p_{\alpha^{n_k}z_1}$.
 Would be this period different from $1$ (\ie equal to some
$d\ge 2$) this would contradict to the fact that the holonomy along
the loop $\tilde p_0([0,1])$ is contractible. Therefore we get that
\begin{equation}
\eqqno(a-alpha3) [1 + b_0(z_1)]a(z_1) = 1.
\end{equation}
In the same way one gets
\begin{equation}
\eqqno(solutions) [1 + b_0(\alpha^nz_1)]a(\alpha^nz_1) = 1
\end{equation}
for all $n\in \nn$. The relation \eqqref(a-alpha3) means that
$a(z_1)$ extends to a holomorphic function to a neighborhood of
zero, which can defined by
\begin{equation}
\eqqno(a-alpha4) a(z_1) = 1/[1 + b_0(z_1)].
\end{equation}
Relation \eqqref(solutions) means that this extension in all points
$\{\alpha^nz_1\}$ is a period of $\tilde p_{\alpha^nz_1}$.
Therefore for every $z_2$ the holomorphic equation
\[
\tilde p(z_1,z_2) = \tilde p(z_1, z_2 + a(z_1))
\]
has a converging to zero sequence of solutions \eqqref(solutions).
I.e., \eqqref(scrol) is proved.

\smallskip\qed

\begin{lem}
\label{trans-alpha} $\call_m$ is a locally closed analytic subset of
$X\setminus\Sing\call$.
\end{lem}
\proof  Let us prove that $\call_m$ is such in a neighborhood of
$m$. If not, \ie if $\call_m$ cuts $D$ by a sequence of points
$m_k\to m$ then all $\tilde p|_{m_k}: \tilde\call_{m_k}\to
\call_{m_k}$ are the periodic coverings of the same leaf $\call_m$
and Lemma \ref{scroling} applies. But this means that holonomy along
the path $\gamma = \tilde p_0([0,1])$ should be trivial. Indeed, for
every $z_1$ close to zero $\tilde p_{z_1}(0) = \tilde
p_{z_1}(a(z_1))$. Therefore, if we denote by $U$ a foliated
neighborhood of zero in $\cc^2= \tilde\calp_{\cc}$, by $U_0$ its
image under $\tilde p$, and by $V$ the corresponding foliated
neighborhood of the point $(0,1)\in \cc^2$ then $\tilde \pi\circ
\tilde p|_W^{-1} \circ \tilde p_{U\cap D} = \id$. But the latter is
a holonomy map corresponding to $\gamma$. At the same time the
holonomy along $\gamma$ is the multiplication by $\alpha$.
Contradiction.

\smallskip Therefore $\call_m$ cuts our initial disc $D$ in
a neighborhood of $m$ only finitely many times $m$. Therefore
$\call_m$ is an imbedded analytic set in a neighborhood of $m$ in
$X$. But $D$ can be taken through {\slsf any} point of
$\call_m\setminus \Sing X$. The case of a cyclic point obviously
follows. This proves that $\call_m$ is an imbedded curve in some
open subset of $X\setminus\Sing\call$.

\smallskip\qed

\smallskip Theorem \ref{parab-case} is proved.

\newsect[sect.PROOF]{Proofs of the main results}
\newprg[prgPROOF.blm]{BLM-trichotomy}
Observe, first of all, that on a projective surface $X$ every
holomorphic foliation can be defined by a {\slsf global} meromorphic
$1$-form. Indeed, from the exact sequence
\begin{equation}
\eqqno(exact-mer1) 0\to \calo^*\to \calm^*\to \calm^*/\calo^*\to 0,
\end{equation}
were $\calm^*$ is a sheaf of non zero meromorphic functions on $X$,
we get
\begin{equation}
\eqqno(exact-mer2) 0\to \cc^* \to H^0(X,\calm^*) \to
H^0(X,\calm^*/\calo^*) \to H^1(X,\calo^*) \to ...
\end{equation}
Here $H^0(X,\calm^*/\calo^*)$ is the group of divisors on $X$ and
$H^1(X,\calo^*)$ is the group of holomorphic line bundles on $X$.
From the other hand on projective $X$ every holomorphic line bundle
admits a meromorphic section and as such is a bundle canonically
associated with a divisor of zeroes and poles of this section.
Therefore the map $H^0(X,\calm^*/\calo^*) \to H^1(X,\calo^*)$ in
\eqqref(exact-mer2) is onto.

\smallskip Let now $\omega_j$ be defining holomorphic $1$-forms of
$\call$ on an open subsets $U_j$ with $\omega_k = f_{kj}\omega_j$.
The meromorphic section of the normal bundle $\caln_{\call}^*$, \ie
of the bundle defined by the cocycle $\{f_{kj}\}\in H^1(X,\calo^*)$,
is a couple $\{f_j\}$ of meromorphic functions on open sets $U_j$
such that $f_k = f_{kj}f_j$ - \ie $\{f_j\}$ is a section of
$H^0(X,\calm^*/\calo^*)$. But then $\Omega \deff \{\Omega_j \deff
f_j^{-1}\omega_j\}$ is a globally defined meromorphic form on $X$
and it still defines our $\call$. If $\{v_j\}$ is an another set of
holomorphic forms which define the same $\call$ then it is
straightforward  to see that the corresponding global meromorphic
form it equal to $\Omega $ modulo a non-zero complex number (on a
compact $X$). Therefore the notion of $\Omega$ being (algebraically)
closed is correctly defined.

\smallskip Following \cite{BLM} we set
\[
\Pole \Omega = \bigcup_j\Pole (f_j^{-1}) \cup \Sing\call ,
\]
\[
\Zero \Omega = \bigcup_j\Zero (f_j^{-1}) \cup \Sing\call .
\]

\smallskip The precise formulation of BLM-trichotomy is as follows,
see Section IV in \cite{BLM}:

\begin{nnthm} {\slsf (Bonatti-Langevin-Moussu)}
\label{blm-thm} Let $\call$ be a holomorphic foliation on a compact
complex surface $X$, which can be defined by a global meromorphic
form $\Omega$. Let $\call_m$ be a minimal leaf of $\call$ such that
$\overline{\call_m}$ doesn't intersects either $\Zero \Omega$ or
$\Pole \Omega $. Then the following three cases are possible:

\smallskip\sli $\Omega$ is algebraically closed.

\smallskip\slii $\call_m$ is a compact leaf of $\call$.

\smallskip\sliii there exists an another leaf with $\call_n$ with
the same closure, \ie  $\overline{\call_n}=\overline{\call_m}$, such
that the holonomy group $\hln (\call_m,m)$ contains a hyperbolic
element.
\end{nnthm}

Let now $x,y$ coordinates in an affine chart $U_0$ of $X$ (\ie $X$
is quite special therefore, like $\pp^2$) in which $\call$ is
defined by the vector field $\v = P(x,y)\frac{\d}{\d x} +
Q(x,y)\frac{\d}{\d y}$. The corresponding polynomial defining form
is then $\omega_0 = P(x,y)dy - Q(x,y)dx$. Write $\omega_0 = P\big(dx
- Q/Pdy\big) =: P\cdot \Omega_0$. Perform the standard coordinate
change in $\pp^2$: $x=1/u, y = v/u$, and get that
\[
\omega_0 = 1/u^{2+\max \{\deg P,\deg Q\}}\omega_1,
\]
where
\[
\omega_1 = u^{1+\max \{\deg P,\deg Q\}-\deg P}P_1(u,v)\big(dv +
\frac{Q_1(u,v)}{u^{1+\max \{\deg P,\deg Q\}-\deg P}P_1(u,v)}du\big)
\]
is the polynomial form in the chart $(U_1,u,v)$ which defines
$\call$ there. Here the exact form of $Q_1$ plays no role, whereas
$P_1(u,v)= u^{\deg P}\cdot P(1/u,v/u)$. Set
\begin{equation}
\eqqno(omega1) \Omega_1 = \frac{1}{u}\big(dv +
\frac{Q_1(u,v)}{u^{1+\max \{\deg P,\deg Q\}-\deg P}P_1(u,v)}du\big)
\end{equation}
and check, finally, that $\Omega_0 = \Omega_1$ on $U_1\cap U_2$.
I.e., $\Omega = \{\Omega_j\}_{j=0}^2$ is the global meromorphic form
which defines $\call$. The term $\Omega_2$ (in the chart $U_2$)
doesn't needs to be computed, because $\pp^2\setminus (U_0\cup U_1)$
is a point and a meromorphic form, which is globally defined on
$\pp^2\setminus \{\text{ point }\}$ extends to $\pp^2$. All this was
needed for just to say that the defining meromorphic form for
$\call$ in the chart $(U_0,x,y)$ is given by
\begin{equation}
\eqqno(omega0) \Omega_0 = dx - \frac{Q}{P}dy = P^{-1}\omega_0,
\end{equation}
and in $(U_1,u,v)$ by \eqqref(omega1). As a result we see that
\begin{equation}
\eqqno(poles) \Pole \Omega = \Zero (P) \cup \{ u=0\} \text{ and }
\Zero \Omega = \emptyset
\end{equation}
for every holomorphic foliation on $\pp^2$. Even in the case
$P\equiv 0$, \ie when $\call$ is a rational quasi-fibration, the
corresponding meromorphic form $dx$ has no zeroes. Therefore the
BLM-trichotomy {\slsf applies to foliations on $\pp^2$}.

\smallskip Let us consider the case (\sli \ie that the defining $\call$
meromorphic form $\Omega$ is algebraically closed. Supposing that
$\call$ is not a rational fibration take the standard affine chart
with coordinates $(x,y)$ and write our form as $\Omega_0 = dy +
f(x,y)dx$, where $f$ is rational. The closeness of $\Omega$ means
that $f$ is a function of $x$ only: $f(x) = p(x)/q(x)$ with $p$ and
$q$ relatively prime. If $q$ is not constant, \ie has a zero, say
$x_0$, then the projective line $\{x=x_0\}$ is tangent to $\call$.
Now either it cuts $\Sing\call$, or it is a leaf. In the latter case
$\call$ must be a rational fibration. But such in $\pp^2$ doesn't
exists. If $q$ is constant then $\Omega$ is exact: $\Omega = d(y+
\int p ) = dF$ with $F$ being rational. Then $\call$ is a
quasi-fibration by level sets of a rational function. A minimal set
in this case is a lever set of $F$ and these level sets do
intersect, \ie $\calm$ cannot be away from the $\Sing\call$.

\smallskip The case (\slii of the trichotomy is ruled out by the
Camacho-Sad formula: $\call_m\cdot \call_m = 0$.

\smallskip Therefore we are left with the case (\sliii of the trichotomy.
Now we are in the position to apply the results of Sections
\ref{sect.PPD} and \ref{sectAB}  via the reduction to the {\slsf
nef} models.

\newprg[prgNEF.p2]{Proof of Theorem \ref{p2-inv-set}} Let $X$
be now the complex projective plane $\pp^2$. Suppose that $\calm =
\overline{\call_m}$ is a minimal leaf, which doesn't intersects
$\Sing\call$. Our aim is to arrive to a contradiction. It will be
convenient for the future references to present this proof as a
sequence of steps. Since for the case when $\call$ is a rational
quasi-fibration our theorem is obvious, we may suppose along the
forthcoming considerations that this case doesn't happens. One more
remark, if $\overline{\call_m}$ is an algebraic curve in $\pp^2$
then our theorem is again trivial. The same for any other leaf in
$\calm$. Therefore we will suppose that this is also not the case.

\smallskip\noindent{\slsf Step 1. $\calp_D = \pp^2\setminus \Sing\call$
for every Poincar\'e disc through $m$. Consequently $\pp^2$ doesn't
contains $\call$-invariant algebraic curves.} The first assertion
was proved in Corollary \ref{cor-ps2}. Suppose $C$ is an $\call$ -
invariant algebraic curve. Then it is the closure of some leaf of
$\call$. Since $\calp_D$ intersects this leaf it should contain it.
Therefore $C=\overline{\call_z}$ for some $z\in D$. But the set
$C\cap D$ cannot accumulate to $m$. Therefore taking a smaller disc
$D_k$ we can arrange that $C\cap D_k = \emptyset$ and, consequently,
$\call_z\cap \calp_{D_k}=\emptyset$. But $\calp_{D_k}$ is still
$\pp^2\setminus \{\text{finite set}\}$. Contradiction.

\smallskip\noindent{\slsf Step 2. Seidenberg's reduction.} Let us
perform the first step in transforming of the foliated pair
$(X,\call)$ to its nef model: {\slsf reduction} of singularities.
This reduction consists in blowing up singular points of $\call$.
I.e., one performs the first blow up $\pi_1:X^1\to X$ with center at
some non-reduced singular point of $\call$ to get a new foliation
$\call^1$ on $X^1$, and one does this with $(X^1,\call^1)$ and so
on, until one gets a foliated pair $(X^N,\call^N)=:(Z,\cale)$ with
only reduced singularities. The finiteness of this procedure is the
content of the theorem of Seindenberg. Denote by $\pi_{Z} \deff
\pi_1\circ ... \circ \pi_N$ the resulting modification. Now we want
to remark few things.

\smallskip\sli First: the proper preimage $\calm^{\cale}$ of
$\calm$ under $\pi_Z:Z\to X$ doesn't intersects neither $\Sing\cale$
nor $E^Z$, where $E^Z$ is the exceptional divisor of $\pi_Z$ .

\medskip\slii Second: if $D\ni m$ is a Poincar\'e disc for $\call$
through $\call_m$ with $\overline{\call_m}=\calm$, then it lifts
under $\pi_Z$ to a Poincar\'e disc for $\cale$.

\smallskip These observations are obvious and result from the
fact that neither $\calm$ nor $D$ do not intersect points in
$\Sing\call$. And the blowing up process goes only over these
points. Denote by $\calp_D^{\cale}$ the Poincar\'e domain of $\cale$
corresponding to $\pi_Z^{-1}(D)$ (the latter can be identified with
$D$ due to the remark (\slii ). Decompose $E^Z=E^Z_i\cup E^Z_d$,
where $E^Z_i$ is the union of $\cale$-invariant components of $E^Z$
and $E^Z_d$ - of dicritical ones, \ie such that $\cale$ is
generically transverse to $E^Z_d$.

\smallskip\sliii Third: $\calp_D^{\cale} = Z\setminus
(E^Z_i \cup \Sing\cale)$ and $\calp_D= \calp_D^{\cale}\setminus
E^Z_d$.

\smallskip More precisely we mean that $\pi_Z|_{\calp_D^{\cale}\setminus
E^Z_d}: \calp_D^{\cale} \setminus E^Z_d\to \calp_D$ is a foliated
biholomorphism. Indeed, an invariant component of $E^Z$ is a closure
of a leaf of $\cale$, which cannot cut $D$. What concerns the
dicritical components, take a leaf $L$ of $\cale$ (not contained in
$E^Z_i$), which cuts one of them. Then, since $L\setminus E^Z$ is a
leaf of $\call$ it should intersect $D$. Therefore $E^Z_d\setminus
E^Z_i\subset \calp_D^{\cale}$. We conclude this remark by noticing
that $\calp_D$ is equal to $\calp_D^{\cale}$ minus a divisor.

\smallskip\noindent{\slsf Step 3. Contracting invariant rational curves.}
Let $(Z,\cale)$ be a Seidenberg's reduction of our foliated pair
$(X,\call)$. Suppose that there exists an irreducible algebraic
curve $C\subset Z$ such that $\calk_{\cale}\cdot C<0$, \ie that $C$
violates the nefness of $\calk_{\cale}$. Then $C$ can be only a
smooth $\cale$ - {\slsf invariant} rational curve with negative
self-intersection, which contains exactly one point from
$\Sing\cale$. But moreover, in our special case $C$ should be a
component of $E^Z_i$ (in general this is not the case). Otherwise it
should be an $\call$-invariant algebraic curve, and this was already
prohibited in Step 1.

\smallskip Contracting $C$ to a point we get a new surface with at
most cyclic singularity $a$ and a foliation downstairs, which is
{\slsf non-singular at $a$}, see \cite{McQ2} for the proofs of all
these facts, as well as \cite{Br2}. If there still left algebraic
curves which violate the nefness of the canonical sheaf we can
repeat these procedure and contract them all. As a result one gets a
{\slsf nef} foliated pair $(Y,\calf)$. We denote by $\pi_Y:Z\to Y$
the corresponding proper modification. Let us emphasis that each
tree of invariant rational curves which contracts this way to a
cyclic point of $Y$ in our special case entirely consists from
irreducible components of $E^Z_i$. Indeed, by the Step 1 there
simply {\slsf do not exist} other invariant curves then components
of $E^Z_i$. Therefore the resulting minimal compact $\calm^{\calf}$
doesn't intersects neither $\Sing\calf$ nor $\Sing Y$. Neither does
the Poincar\'e disc $D$.

\smallskip What concerns the Poincar\'e domain $\calp_D^{\calf}$ of
the foliation $\calf$ corresponding to $D$ the situation is more
delicate. $\calp_D^{\calf}$ obviously inherits the whole
$\calp_D^{\cale}$ and, may be, becomes some cyclic points. More
precisely $\calp_D^{\cale} =\calp_D^{\calf}\setminus \{c_i\}$, where
$\{c_i\}$ are some of cyclic points of $Y$ (finite in number). The
nature of these $\{c_i\}$ can be understood from the following:

\begin{lem}
\label{tree} Let $T_i$ be a tree of $\cale$-invariant rational
curves which contracts by $\pi_Y:Z \to Y$ to a cyclic point $c_i\in
\calp_D^{\calf}$. Then for an appropriate neighborhood $V_i\supset
T_i$ one has that:

\smallskip{\bf a)} $V_i\setminus T_i\subset \calp_D^{\cale}$.

\smallskip{\bf b)} $U_i\deff \pi_Y(V_i)$ is the standard cyclic
neighborhood of $c_i$.

\smallskip{\bf c)} $\calp_D^{\cale}$ is {\slsf leafwise biholomorphic}
to $\calp_D^{\calf} \setminus \{c_i\}$.
\end{lem}
\proof $T_i$ doesn't intersects any other $\cale$-invariant curve.
Indeed, would $C$ be such, then its image downstairs in $Y$ would be
an $\calf$-invariant curve passing through $c_i$, and therefore it
would be a closure of a leaf $\calf_z$ from $\calp_D^{\calf}$ and as
such would intersect $D$. But this cannot happen, because all
modifications took place away from $D$.

\smallskip Taking $V_i$ small enough we can insure that $\Sing\cale
\cap V_i \subset T_i$ and that no component of $E^Z_i$ other then
some from $T_i$ intersects $V_i$. Now the both (a) and (b) become
clear.

\smallskip  Item (c) readily follows from (a), because $T_i$-s are 
disjoint from $\calp_D^{\cale}$.

\smallskip\qed

\begin{rema} \rm
\label{non-tree} For a cyclic point $c\not\in \calp_D^{\calf}$ the
contracting to it tree $T$ may well intersect some other
$\cale$-invariant curves, which are not subject of contraction. But
these points are out from our process.
\end{rema}

\smallskip The universal covering Poincar\'e domain $\tilde
\calp_D^{\calf}$ is Hausdorff and Rothstein, because $(Y,\calf)$ is 
nef, see Corollary \ref{nef-cov} and Proposition \ref{rothstein4}. 
Therefore from this point we can  split our proof onto parabolic and
hyperbolic cases.

\smallskip\noindent{\slsf Step 4. End of the proof: parabolic case.}
We suppose that $(Y,\calf)$ is parabolic, \ie all its leaves are
parabolic. Therefore all our constructions from Sections
\ref{sect.PPD} and \ref{sectAB} do apply. We obtain that our minimal
leaf $\calf_m=\call_m$ closes to an algebraic curve, or it is
$\cc^*$ and in the latter case it is a locally closed  analytic set
in $Y\setminus \Sing\calf $. The case of a compact curve was
prohibited in the Step 1. What concerns the second possibility let
us recall that at the same time $\overline{\calf_m}=\calm^{\calf}$
is a compact set in $Y\setminus (\Sing\calf \cup \Sing Y)$ by Step
3. The limiting set $\lim \calf_m$ cannot contain $\calf_m$ because
$\calf_m$ is closed in some open subset of $Y\setminus (\Sing\calf
\cup \Sing Y)$. Therefore it should be empty by the minimality of
$\overline{\calf_m}$. Therefore $\calf_m$ is a closed leaf of
$\calf$ and then such is also $\call_m$. But this was again
prohibited by Step 1.

\smallskip\noindent{\slsf Step 5. End of the proof: hyperbolic case.}
Now consider the case when $\calf$ is hyperbolic. If our leaf
$\calf_m$ happen to belong to the exceptional set $\cala$ in
McQuillan's alternative then the closure of $\call_m$ is also an
algebraic curve in $\pp^2$ and we are done.

\smallskip Therefore we can suppose that $\calf_m$ is hyperbolic and
$\calp_D^{\calf}\subset Y\setminus \cala$. Indeed, if
$\calp_D^{\calf}$ still contains a leaf from $\cala$ then we take a
smaller $D$. Therefore the Lobatchevski-Poincar\'e norm is
continuous on $\calp_D^{\calf}$ and the Theorem \ref{prop-disc-hyp}
applies. We have that $\tilde p^{\calf}:\tilde\calp_{\cc}^{\calf}\to
\calp_D^{\calf}$ is a regular (cyclic) covering. But
$\calp_D^{\calf}\supset \calp_D^{\cale} \supset \calp_D$ and
therefore contains a lot of rational curves, also such that do not
pass through the cyclic points. Therefore
$\tilde\calp_{\cc}^{\calf}$ also contains them. But
$\tilde\calp_{\cc}^{\calf}$ posed  a holomorphic surjective
projection onto $\cc$ and therefore these curves should be the
fibers. This is impossible, because the fibers of
$\tilde\calp_{\cc}^{\calf}$ are Kobayashi hyperbolic. Contradiction.

\smallskip Theorem \ref{p2-inv-set} is proved.

\begin{rema} \rm
If $\{c_i\}$ is non empty then $\tilde\calp_{D}^{\cale}$ is {\slsf
not Hausdorff} and we cannot work with it.
\end{rema}

\begin{rema} \rm
\label{proof-proj} What concerns the proof of Corollary \ref{proj}
let us just remark the following items. Convexity of $\pp^2\setminus
M$ is guaranteed by Fujita's theorem. Therefore the Levi foliation
$\call$ of $M$ extends onto $\pp^2$. Now $M$ should contain a
minimal leaf $\calm$ of $\call$. But already $\calm$ intersects
$\Sing\call$ due to the Theorem \ref{p2-inv-set}. Likewise $M$ does
and a fortiori it cannot be smooth. Corollary \ref{proj} is proved.
\end{rema}

\newprg[prgPROOF.approach1]{Limiting behavior of leaves with hyperbolic
holonomy}

We shall prove now the Corollary \ref{approach1} from the
Introduction, \ie we shall detect the reasons for the failure of
steiness of $\calp_D$ of a minimal leaf with hyperbolic holonomy.
This will be done as in the proof of Theorem \ref{p2-inv-set} along
the reduction to the nef model. First of all let us remark that if
$\call$ is a rational quasi-fibration then {\slsf any} of its leaves
cannot have a hyperbolic holonomy.

\smallskip Furthermore, according to Corollary \ref{cor-ps1} $\calp_D$
is not pseudoconvex at some boundary point $z_0$ if and only if
$z_0$ is an isolated boundary point of $\calp_D$ and $z_0\in
\Sing\call$. In the sequel we denote as $\{z_i\}_{i=1}^k$ the set of
all such points. By $\{w_i\}_{i=1}^l$ we denote the set of points of
$\Sing\call$ which belong to $\calm$, \ie $\call_m$ approaches to
them, but which are not isolated points of $\d\calp_D$. I.e.,
$\d\calp_D$ is a sort of a Levi flat "cone" with vertices at these
$\{w_i\}$.

\smallskip\noindent{\slsf Step 1. If the set $\{w_i\}$ is not empty
then $\calm$ admits a Stein invariant neighborhood. If it is empty
then $\calp_D=\pp^2 \setminus \Sing\call$.} If $\{w_i\}$ is non
empty we can add to $\calp_D$ all $\{z_i\}_{i=1}^k$ and obtain a
domain $\bar\calp_D$ which is still different from $\pp^2$, because
it doesn't contains any of $\{w_i\}$. At the same time $\bar\calp_D$
is pseudoconvex by Lemma \ref{loc-ps1} and is obviously invariant.
Therefore it is a Stein invariant neighborhood as claimed. If
$\{w_i\}$ is empty we obtain in the same manner a pseudoconvex
domain $\bar\calp_D$ which contains an invariant compact $\calm$.
Would $\bar\calp_D$ be different from $\pp^2$ it would be Stein.
Contradiction: a Stein domain cannot contain a compact invariant
set. Therefore we proved that $\calp_D = \pp^2\setminus \Sing\call$
in this case.

\smallskip Suppose that there exist an invariant rational curve
$C$ in $\pp^2$. If $C$ is not the closure of $\call_m$ then we get a
contradiction exactly as in the Step 1 of the proof of Theorem
\ref{p2-inv-set}. If $C$ is the closure of $\call_m$ then the case
(\slii of our Corollary occurs.

\smallskip  From now one we can suppose that $\calp_D = \pp^2
\setminus \Sing\call$ for all Poincar\'e discs and therefore
$\pp^2$ doesn't contains invariant rational curves.

\smallskip\noindent{\slsf Step 2. Seidenberg's reduction.}
Let $(Z,\cale)$ be a Seidenberg's reduction of $(X,\call)$ and
$\pi_Z:Z\to X$ the corresponding modification. Then as in the Step 2
before we have that $\calp_D = \calp_D^{\cale}\setminus E^Z_d$. If
$\cale_m$ is the leaf of $\cale$ through $m$ then $\call_m = \cale_m
\setminus \{e_j\}$, where $\{e_j\}$ is a discrete (in the topology
of $\cale_m$) set of vanishing ends of $\call_m$.

\smallskip\noindent{\slsf Step 3. Contraction of invariant rational
curves.} When contracting invariant rational curves, which violate
the nefness of the canonical bundle of $\cale$ let us observe that
any of $\{e_j\}$ cannot arrive to a cyclic point of $(Y,\calf)$ -
the nef model of $(X,\call)$.

\smallskip Moreover, $\calp_D^{\cale} = \calp_D^{\calf}\setminus
\{c_i\}$, where $\{c_i\}$ are some cyclic points of $Y$. We have the
following two possibilities for $\calf_m$.

\smallskip{\bf 1)} $\calf_m$ cuts some of $\{c_i\}$, then we get
some more vanishing ends (finite number).

\smallskip{\bf 2)} $\calf_m$ accumulates to some of $\{c_i\}$ - this
case is irrelevant for us.

\smallskip\noindent{\slsf Step 4. Parabolic case.} If the nef model
$(Y,\calf)$ is parabolic we have two posiibilities for $\calf_m$.

\smallskip{\bf 1)} {\slsf $\calf_m$ is a torus or, a sphere.} Then so is also
$\overline{\call_m}$ and the number of ends occurred along the
previous steps is finite, \ie cases (\slii and (\sliii of our Corollary
occur.

\smallskip{\bf 2)} {\slsf $\calf_m = \cc^*$ and is a locally closed
analytic subset in $Y\setminus \Sing\calf$.} The limiting set again
cannot contain $\calf_m$ and therefore should be contained in
$\Sing\calf$ - a finite set. By the theorem of Remmert-Thullen
$\calf_m$ closes to a rational curve with two vanishing ends.
$\call_m$ might get more of them. I.e., the case (\slii occurs
again.

\smallskip\noindent{\slsf Step 5. Hyperbolic case.} If $\calf_m$ is
contained in an exceptional set $\cala$ of McQuillan'e theorem then
we have again one of the cases (\slii , or (\sliii. If not, then the
same proof as in the corresponding Step 5 before shows that this
option is impossible.

\smallskip Corollary \ref{approach1} is proved.

\newprg[prgPROOF.p1p1]{Minimal leaves in the product of projective lines}

We shall prove now Corollary \ref{p1p1-inv-set}. Let $\call_m$ be a
minimal leaf of a holomorphic foliation $\call$ on $\pp^1\times
\pp^1$. Denote by $\calm$ its closure $\overline{\call_m}$. Suppose
that $\calm$ doesn't intersects $\Sing\call$. Computations (simpler)
as that from Subsection \ref{prgPROOF.blm} show that a defining
meromorphic form of a holomorphic foliation on $\pp^1\times \pp^1$
has no zeroes. Therefore the BLM-trichotomy applies.

\smallskip\noindent{\slsf Step 1. The use of BLM-trichotomy.} According
to this trichotomy we must study three cases. Let $\Omega$ be a
global meromorphic defining form of $\call$.

\smallskip\noindent{\slsf Case 1. $\Omega$ is closed.} Supposing that
$\call$ is not a rational fibration take the standard affine chart
with coordinates $(x,y)$ (see the Subsection \ref{prgPROOF.blm}) and
write our form as $\Omega_0 = dy + f(x,y)dx$, where $f$ is rational.
It means that $f$ is a function of $x$ only: $f(x) = p(x)/q(x)$ with
$p$ and $q$ relatively prime. If $q$ is not constant, \ie has a
zero, say $x_0$, then $\pp_{x_0}\deff \{x_0\}\times \pp^1$ is
tangent to $\call$. Now either it cuts $\Sing\call$, or it is a
leaf. In the latter case $\call$ is an obvious rational fibration.
I.e., the case (\slii of Corollary \ref{p1p1-inv-set} occurs. If $q$
is constant then $\Omega$ is exact: $\Omega = d(y+ \int p ) = dF$
with $F$ being rational. Then $\call$ is a quasi-fibration by level
sets of a rational function. A minimal set in this case is a leaf
and this leaf can not to intersect the singularity locus =
indeterminacy set, if and only if $\call$ is one of two rational
fibrations.

\smallskip We still have two cases.

\smallskip\noindent{\slsf Case 2. $\calm$ is a leaf.} Write $\calm = aE_1 +
bE_2$, where $E_1=\{pt\}\times \pp^1$ and $E_2=\pp^1\times \{pt\}$ -
generators of $H^2(X,\zz)$. Then $\calm^2 = 2ab$ and this should be
zero. Therefore, again $\calm = \call_m$ is $\{pt\}\times \pp^1$ or
vice versa and we find ourselves in the case (\slii of our
Corollary.

\smallskip\noindent{\slsf Case 3. $\hln (\call_m,m)$ contains a
hyperbolic element.} More precisely, this is so for a, may be, some
other leaf $\call_n$ with the same closure. Take the expanded
Poincar\'e domain $\tilde\calp_{\cc}$. If it is hyperbolic the proof
is identical to that of case $\pp^2$, because $\pp^1\times \pp^1$ 
also has sufficiently many rational curves, with only one difference - 
eventual appearance
of domains of the form $D\times \pp^1$ when studying pseudoconvexity
of $\calp_D$ or of $\tilde\calp_D$. But then $\call$ is the obvious
rational fibration and the case (\slii of Corollary occurs.
Therefore hyperbolic case cannot happen. If $\tilde\calp_{\cc}$ is
parabolic then everything is the same as in $\pp^2$.

\smallskip Corollary \ref{p1p1-inv-set} is proved.

\begin{rema} \rm
\label{proof-prod-proj} What concerns Corollary \ref{prod-proj} let
us remark that if $\pp^1\times \pp^1 \setminus M$ is not Stein then
by Fujita's theorem $M$ is already of the form $\gamma \times \pp^1$
(or $\pp^1\times \gamma$). The rest is obvious.
\end{rema}

\newprg[prgPROOF.bb]{Briot-Bouquet foliations} For the
study of Briot-Bouquet foliations we shall need an analog of
Fujita's theorems for pseudoconvex domains in $\Sigma_g\times
\pp^1$, where $\Sigma_g$ is a compact complex curve of genus $g\ge
1$.

\begin{nnprop}
\label{bb-levi} Let $D$ be a locally pseudoconvex domain over
$\Sigma_g\times \pp^1$. Then:

\smallskip\sli either $D= D_1\times \pp^1$;

\smallskip\slii or $D= \Sigma_g\times D_1$;

\smallskip\slii or $D$ is Stein.
\end{nnprop}
\proof Denote by $\check\cc^2$ the punctured $\cc^2$, \ie
$\check\cc^2 \deff \cc^2\setminus \{0\}$. Take the canonical
projection $p : \Sigma_g \times\check\cc^2\to \Sigma_g\times\pp^1$.
Let $\tilde D$ be the preimage of $D$ under $p$. Consider $\tilde D$
as a domain over $\Sigma_g\times\cc^2$. It is locally pseudoconvex
there over all points apart, may be, of the points in
$\Sigma_g\times \{0\}$. It is also $\cc^*$-invariant under the
action
\[
\lambda\cdot [s; w_1,w_2] = [s;\lambda w_1, \lambda w_2].
\]
Denote by $\hat D$ the pseudoconvex envelope of $\tilde D$. By the
Theorem of Grauert-Remmert, see \cite{GR1}, $\hat D$ might be
defferent from $\tilde D$ only over some points over $\Sigma_g\times
\{0\}$. It is also invariant under the action in question. We have
the following three cases.

\smallskip\noindent{\slsf Case 1. $\hat D$ contains a point over the
$\Sigma_g\times \{0\}$.} Let $(s,0)$ be this point. But then we see,
acting  by $\cc^*$, that $\hat D$ contains the fiber $\{s\}\times
\check\cc^2$. That means that $D$ is of the form $V\times \pp^1$ for
some domain $V\subset \Sigma$.

\smallskip\noindent{\slsf Case 2. $\hat D = \tilde D$ but is not Stein.}
By theorem of Brun, see \cite{Brn}, $\tilde D = \Sigma\times B$ for
some $B\subset \cc^2$. But then $D=\Sigma_g\times V$ for the
projection $V$ of $B$.

\smallskip\noindent{\slsf Case 3. $\hat D = \tilde D$ and is Stein.}
Applying  the theorem from \cite{MM} we get that $D$ is Stein
itself.

\smallskip\qed

\smallskip Now let us accomplish the
classification of exceptional minimals in Briot-Bouquet foliations.
Let $\calm$ be an exceptional minimal of a holomorphic foliation
$\call$ on $\Sigma_g\times \pp^1$, $g\ge 1$.

\smallskip Let $\tilde\Sigma_g$ be the universal covering of $\Sigma_g$ (\ie
$\tilde\Sigma_g$ is $\cc$ or $\Delta$). Denote by $y$ the natural
coordinate on $\tilde\Sigma_g$ and by $x$ an affine coordinate on an
appropriate chart of $\pp^1$. Let $\Omega = dy + f(x,y)dx$ be a
meromorphic form defining $\call$. Here $f(x,y)$ is rational with
respect to $x$ and automorphic with respect to $y$. The same
analysis as in Subsection \ref{prgPROOF.blm} shows that $\Omega$ has
no zeroes and BLM-trichotomy applies.

\smallskip Identically to the Step 1 of Subsection \ref{prgPROOF.p1p1}
we get that if $\Omega$ is closed then $\calm$ can only a fiber $E_1
\deff \{\pt\}\times \pp^1$ or $E_2 \deff \Sigma_g\times \{pt\}$. The
same holds also for the second option in BLM-trichotomy.

\smallskip Therefore we are left with the case when $\calm =
\overline{\call_m}$ and the leaf $\call_m$ contains a loop with
contractible hyperbolic holonomy. All results of this paper are
applicable because we have a repalacement of the Fujita theorem -
Proposition \ref{bb-levi}.
Let $\calp_D$ be the Poincar\'e domain of a Poincar\'e disc $D$
through $m$.  If it is locally pseudoconvex but non Stein
then the situation is clear via Proposition \ref{bb-levi} - our
$\call$ is a fibration. If $\calp_D$ is Stein then we have
three following cases to consider. Remark that the case when $\call$
is a rational fibration is obvious, as usual, and will be ignored
(but it produces the possibility: $\call_s = \{s\}\times \pp^1$ for
all $s\in \Sigma_g$).

\smallskip\noindent{\slsf Case 1. The nef model $(Y,\calf)$ of
$(X, \call)$ is parabolic and $g\ge 2$.} Projections of leaves onto
$\Sigma_g$ can be only constants. Therefore $\call$ is again a
rational fibration, the already excluded case.

\smallskip\noindent{\slsf Case 2. The nef model of $\call$ is
parabolic and $\Sigma_1 = \ttt $.} If $\call_m$ is a torus or, a 
sphere then this
case was already considered. Otherwise $\call_m$ is $\cc^*$ imbedded
as a locally closed  analytic subset to $Y\setminus \Sing\call $ and
its closure doesn't intersect $\Sing\call$. This is again a
contradiction unless $\call_m$ is not a fiber.

\smallskip\noindent{\slsf Case 3. The nef model of $\call$ is
hyperbolic.} As in the case of $\pp^2$, or $\pp^1\times\pp^1$ we
have that $\calp_D= X\setminus \Sing\call$, otherwise we would
either get a contradiction with maximum principle, or conclude that
$\call$ is a trivial fibration $\call_p = \Sigma_g\times \{pt\}$ or
$\{pt\}\times \pp^1$. But $\tilde p: \tilde\calp_{\cc}\to \calp_D =
X\setminus \{ \text{points} \}$ is a regular covering by Theorem
\ref{prop-disc-hyp}. But then $\calp_D$ contains a lot of rational
curves, and so the $\tilde\calp_{\cc}$ does. Contradiction.

\smallskip Therefore we proved the Corollary \ref{bb-min}  from the
Introduction.

\newprg[prgPROOF.hirz1]{Levi problem  in Hirzebruch surfaces}

As in the case with Briot-Bouquet foliations  we need to  discuss
the Levi problem first. Recall that for $k\ge 1$ the Hirzebruch
surface $H_k$ is the projectivization of the bundle $E=\calo \oplus
\calo (-k)\to \pp^1$. We need to understand the description of
locally pseudoconvex domains in $H_k$, which are not Stein.

\smallskip Set $\check\cc^2\deff \cc^2\setminus \{0\}$ and consider
the domain $U\deff \check\cc^2\times\check\cc^2 \subset \cc^4$.
Coordinates in $\cc^4$ we denote as $[z_1, z_2; w_1, w_2]$, or as
$(z,w)$. Set furthermore $E_1\deff \check\cc^2\times \{0\}$ and
$E_2\deff \{0\}\times\check\cc^2$.

\smallskip Denote by $G$ the compact Lie
group $\sss^1\times\sss^1$ and by $G_{\cc}=\cc^*\times \cc^*$ its
complexification. Elements of $G_{\cc}$ we denote as $(\lambda
,\mu)$, where $\lambda$ and $\mu$ are non-zero complex numbers.
Consider the following action of $G_{\cc}$ on $U$:
\begin{equation}
\eqqno(action)
\begin{matrix}
\lambda \cdot [z_1, z_2; w_1, w_2] = [\lambda z_1, \lambda z_2; w_1,
\lambda^kw_2],\cr \mu \cdot [z_1,z_2; w_1,w_2] = [z_1, z_2; \mu w_1,
\mu w_2].
\end{matrix}
\end{equation}
It is not difficult to check that $U/G_{\cc}\equiv H_k$. Denote by
$p:U\to H_k$ the natural projection. The image of $\{w_1=0\}$ is the
exceptional curve in $H_k$, it will be denoted as $E$. By $\pi
:H_k\to \pp^1$ we denote the natural projection coming from the
projection $(z,w)\to z$ after a factorization. $\pi$ realizes $H_k$
as a ruled surface over $\pp^1$.

\smallskip Let us gather the information we need about the Levi problem
in Hirzebruch surfaces:

\begin{nnprop}
\label{levi-pr-hirz} Let $D$ be a locally pseudoconvex domain in a
Hirzerbruch surface $H_k$, $k\ge 1$. Then the following cases are
possible:

\smallskip\sli $D=\pi^{-1}(V) $ for some domain $V\subset \pp^1$ (this
includes the case $D=H_k$);

\smallskip\slii $D$ is a $1$-complete neighborhood of the
exceptional curve $E$;

\smallskip\sliii $D=B\setminus C$, where $B$ is a $1$-complete
neighborhood of the exceptional curve $E$;

\smallskip\sliv $D=H_k\setminus E$;

\smallskip\slv $D$ is Stein.
\end{nnprop}
\proof Here by saying that $D$ (or $B$) is $1$-complete we mean that
after the contraction the exceptional curve $E$ it becomes Stein.
Let $D$ be a locally pseudoconvex domain in $H_k$. Denote by $\tilde
D$ its preimage $p^{-1}(D)$. $\tilde D$ is a domain in $\cc^4$,
which is:

\begin{itemize}
 \item invariant under the action of $G_{\cc}$;

\smallskip\item locally pseudoconvex at every boundary point except,
may be, points in $E_1\cup E_2$.
\end{itemize}

\noindent Denote by $\hat D$ the envelope of holomorphy of $\tilde
D$. This is a locally pseudoconvex domain over $\cc^4$, which is
invariant under the action of $G_{\cc}$. Again, due to \cite{GR1}
$\hat D$ might be different form $\tilde D$ only by some points
over over $E_1\cup E_2$.  We have several cases to consider.

\smallskip\noindent{\slsf Case 1. $\hat D$ contains a point over $E_1$.}
Let $(z^0,0)$ be this point. Since a schlicht neighborhood over
$(z^0,0)$ is contained then in $\hat D$ we get, by acting with $\mu$
on this neighborhood, that $\hat D$ contains $\{z^0\}\times \cc^2$.
In fact $\hat D$ contains all fibers $\{z\}\times \cc^2$ for $z$ in
a neighborhood of $z^0$. That means that $D$ has the form
$\pi^{-1}(v)$ for some $V\subset \pp^1$, \ie the case (\sli of our
Proposition occurs.

\smallskip\noindent{\slsf Case 2. $\hat D$ contains a point over $E_2$.}
Let $(0,w^0)$ be this point. Then $\hat D$ contains a schlicht
neighborhood over the polydisc $\Delta^2_{\eps}(0)\times
\Delta^2_{\eps}(w^0)$ for some $\eps >0$. Acting on this set by
$\lambda$ we easily get that $\hat D$ contains a schlicht domain
over a neighborhood of $\Delta^2_{\eps}(0)\times
\Delta_{\eps}(w_1^0) \times \{|w_2| \ge w_2^0|\}$.

\smallskip{\slsf Subcase 2a. $w_1^0=0$.} Then, acting by $\mu$ we
obtain that $\hat D$ contains a schlicht domain over the cone
\begin{equation}
\eqqno(cone1)
\Delta^2_{\eps}\times \{|w_2| \ge \frac{|w_2^0|}{\eps }|w_1|\}.
\end{equation}
In this case $D$ is a $1$-complete neighborhood of the exceptional
curve $E$, \ie the case (\slii occurs.

\smallskip{\slsf Subcase 2b. $w_1^0\not=0$.} Acting again by $\mu$
we obtain that $\hat D$ contains a schlicht domain over a
neighborhood of the cone
\begin{equation}
\eqqno(cone2)
\Delta^2_{\eps}\times \big(\{|w_2| \ge \frac{|w_2^0|}{|w_1^0|}|w_1|
\setminus \{w_1=0\}\big)\}.
\end{equation}
In this case $D=B\setminus E$ for some $1$-complete domain $B$
containing  $E$, \ie the case (\sliii  occurs if $w_2^0\not=0$ and
the case (\sliv if $w_2^0=0$.

\smallskip\noindent{\slsf
Case 3. $\hat D = \tilde D$, \ie $\tilde D$ is Stein.} In this case
$D$ is a factor of a Stein domain by an free holomorphic action of a
complexification of a compact Lie group. Theorem of Matsushima and
Morimoto, see \cite{MM}, assures then that $D$ is Stein itself.

\smallskip\qed

\newprg[prgPROOF.hirz2]{Exceptional minimals and Levi flats in
Hirzebruch surfaces} In the very same spirit one can clarify the
situation with exceptional minimals in Hirzebruch surfaces $H_k$ for
$k\ge 1$, we keep notations of Subsection \ref{prgPROOF.hirz1}.

\smallskip
Let us prove the Corollary \ref{inv-set-hirz} from the Introduction.
BLM-trichotomy is applicable in this case to and gives the following
possibilities. If the defining meromorphic form is closed then
$\call$ can be only the canonical rational fibration. To understand
the second case let $C$ be a smooth, irreducible algebraic curve in
$H_k$, which is a leaf of some holomorphic foliation $\call$ on
$H_k$. If $E$ denotes the exceptional curve and $F$ - the fiber,
then write $C = nE + lF$. Suppose that $C$ is neither $E$ nor $F$.
Then it should intersect $E$ non-negatively:
\begin{equation}
\eqqno(non-neg)
(nE + lF)\cdot E = -kn + l \ge 0 \text{ and this implies } l\ge nk.
\end{equation}
At the same time by Camacho-Sad formula we have that
\begin{equation}
\eqqno(cam_sad)
0 = C^2 = (nE + lF)^2 = - n^2k + 2ln \ge -n^2k + 2n^2k = n^2k>0.
\end{equation}
Therefore an imbedded curve cannot be a leaf of a holomorphic
foliation unless it is a fiber of the canonical rational fibration.

\smallskip The case with the leaf with hyperbolic holonomy is identic
to the already considered cases of $\pp^2$, $\pp^1\times \pp^1$ and
$\Sigma_g\times \pp^1$. We shall not repeat it.

\smallskip\qed

\smallskip Now we can reach the understanding of Levi flats in
Hirzebruch surfaces. Let $M$ be a Levi flat hypersurface in $H_k$.
Set $H_k\setminus M = D^+\cup D^-$. Then both $D^{\pm}$ are locally
pseudoconvex. Therefore the Proposition \ref{levi-pr-hirz} applies.
Case (\sli of the Proposition \ref{levi-pr-hirz} produces the
statement of this Corollary. Cases (\slii ,(\sliii and (\sliv cannot
happen. In the case \slv we  can extend the Levi foliation onto the
whole of $H_k$.

\smallskip Now our $M$ should contain an exceptional minimal $\calm =
\overline{\call_m}$. Therefore by Corollary \ref{inv-set-hirz} the Levi
foliation is the rational fibration and $\calm$ is a fiber. Moreover,
since all leaves are $\pp^1$-s the hypersurface $M$ should be
foliated by them.

\smallskip Corollary \ref{levi-fl-hirz} follows.

\smallskip\qed

\newprg[prgPROOF.pn]{Minimal sets in projective spaces}

Now we shall prove Corollary \ref{pn-inv-set} from Introduction. Let
$\call$ be a codimension one holomorphic foliation in $\pp^n$. First
we shall make a preparatory step.

\smallskip\noindent{\slsf Step 1. Generic sections by hyperplanes.}
We follow \cite{CLS2} and \cite{BLM}. An open subset of $\pp^n$ is
called {\slsf generic} if its complement is {\slsf thin}, \ie is
contained in at most countable union of locally closed proper
analytic subsets. Recall that if $w=[w_0:...:w_n]$ is a point in the
dual $\pp^{n^*}$, then the corresponding hyperplane $E_w$ in $\pp^n$
is given by the equation $ w_0z_0 + ... + w_nz_n = 0$. There is a
generic subset $G^*_1$ in the dual $\pp^{n^*}$ such that for every
$w\in G^*_1$ the plane $E_w$ from this subset the following holds:

\medskip\sli $E_w$ is not contained in any leaf of $\call$.

\smallskip\slii All components of $\Sing\call^w$ have codimension at
least two. Here $\call^w\deff \call|_{E_w}$ - the restriction of
$\call$ to $E_w$.

\smallskip\noindent For the proof see Lemma 10 in \cite{CLS2}. In
\cite{BLM}, Proposition in Section III, it is proved that in
addition to (\sli , (\slii one has the following assertion. There
exists a generic subset $G^*_2\subset \pp^{n^*}$ and an integer
number $d\ge 0$ such that for every $w_0\in G^*_2$ the hyperplane
$E_{w_0}$ satisfies also the following:

\smallskip\sliii
there exist exactly $d$ points $\{z_1,...,z_d\}$ in
$E_{w_0}\setminus \Sing\call$ where $E_{w_0}$ is tangent to $\call$.

\medskip
Moreover, for every $w_0\in G^*_2$ there exists a neighborhood
$w_0\in W\subset G^*_2$ and a neighborhood $Z$ of $\{z_1,...,z_d\}$,
$\bar Z\cap \Sing\call = \emptyset$, such that for any $w\in W$ all
points in the plane $E_w$, in which $E_w$ is tangent to $\call$, are
contained in $Z$. In particular all these points are contained in a
compact away from $\Sing\call$.

\smallskip\noindent{\slsf Step 2. Points of tangency.} Let $\calm$
be a closed invariant set of our foliation $\call$. For $w\in
G^*_1\cap G^*_2$ denote by $T_w$ the finite set of points where
$\call$ is tangent to $E_w$. Let us see that:

\smallskip\noindent\sliv there exists open and dense subset $G^*_3
\subset \pp^{n^*}$ such that for every $w\in G^*$ one has $T_w\cap
\calm = \emptyset$.

\smallskip Indeed, take $w$ as above and let $T_w = \{z_1,...,z_d\}$.
Suppose that $z_d\in \calm$. Since $\calm$ has empty interior (by
minimality) we can take $z_d^{'}$ close to $z_d$, which doesn't
belongs to $\calm$. Take $w'$ close to $w$ such that $\call$ is
tangent to $T_{w'}$ at $z_d^{'}$. We obtained that $T_{w'}\cap
\calm$ has at most $d-1$ points and the same for hyperplanes in a
neighborhood of $T_{w'}$. After $d$ steps we obtain an open dense
set as in (\sliv .

\smallskip\noindent{\slsf Step 3. Measure positivity of boundaries of
invariant sets.} We shall prove the statement of Corollary
\ref{pn-inv-set} by induction. For $n=2$ Corollary \ref{pn-inv-set}
is proved by Theorem \ref{p2-inv-set}. Therefore, from now on $n\ge
3$.

\smallskip
For $w\in G^*=G^*_3$ from the Step 2 denote $\call^w$ the
restriction of $\call$ to the hyperplane $E_w$. Then:

\smallskip\sli $\Sing\call^w \subset (\Sing\call\cap E_w) \cup T_w$,
where $T_w$ is disjoint from $\Sing \call$ and from $\calm$;

\smallskip\slii $\calm^w\deff \calm\cap E_w$ is a closed invariant
set of $\call^w$ (may be not minimal).

\medskip $\calm^w$ cuts $\Sing\call^w$ by a set of positive
$(2n-6)$-measure for every $w$ from the set of full measure. If
$n\ge 4$ we are done, because this set can be contained only in
$\Sing\call\cap E_w$ (and not in $T_w$, which is just finite). If
$n=3$ we have that $\calm^w$ contains at least one point from
$\Sing\call^w$ and this point is not from $T_w$, \ie it can be only
from $\Sing\call \cap E_w$, and we are done again.

\smallskip Corollary \ref{pn-inv-set} is proved.

\newsect[sect.RAT]{Pseudoconvexity vs. rational curves}

\newprg[prgPPD.an-ob]{Analytic objects}

Throughout this paper we used extension properties of some "analytic
objects" like holomorphic/meromorphic mappings, foliations etc.
These analytic objects have the following two decisive properties:

\medskip{\slsf A1)} The Hartogs type extension theorem is valid for
them. I.e., if any of these objects is given on the Hartogs figure
$H^2_{\eps}$ then it extends to the same type of object onto the
bidisc $\Delta^2$.

\medskip{\slsf A2)} They obey the uniqueness theorem. I.e., if two of
them $\sigma_1$, $\sigma_2$ are defined on a connected manifold (or
space) $U$ and for some open $V\subset U$ one has
$\sigma_1|_V=\sigma_2|_V$, then $\sigma_1=\sigma_2$.

\begin{rema}\rm {\slsf 1.} A holomorphic object can became to be a
meromorphic after extending in a non-Stein case. In the Stein case
such thing cannot happen for functions, holomorphic forms, sections
of holomorphic bundles. But it can happen, even in the Stein case,
for mappings with values in K\"ahler manifolds, for example.

\medskip\noindent{\slsf 2.} A holomorphic foliation on $X$ (resp.
singular holomorphic foliation) is defined by a holomorphic section
(resp. meromorphic section) of the projectivized tangent bundle $\pp
(TX)$. The products $U\times \pp (T_xX)\equiv\pp (TX)|_U$, where $U$
is a local chart, are K\"ahler. Therefore the extension works.
Involutibility, being a holomorphic condition, is preserved by
extension. The extended foliation might become singular.

\medskip\noindent{\slsf 3.} Obvious generalizations to $n>2$ will be
not needed for us in this paper.
\end{rema}

\begin{defi}
\label{an-obj} A sheaf on analytic objects on a complex manifold
(or, a normal space) $X$ is a sheaf of sets which sections obey
properties (A1) and (A2).
\end{defi}

\newprg[prgRAT.ext-ob]{Extension of analytic objects}

First, we recall the following:

\begin{defi}
\label{1-conv} A real hypersurface $\Sigma$ in $2$-dimensional
complex manifold $X$ is called strictly $1$-convex if it locally
admits a smooth defining function $\rho : X\cap U\to \rr$, $\Sigma
\cap U = \{\rho =0\}$, such that the eigenvalues of the Levi form of
$\rho$ are strictly positive at each point. Such function is called
strictly $1$-convex, or strictly plurisubharmonic. More precisely,
$\Sigma$ is called to be strictly $1$-convex from the side $U^- = \{
z\in U:\rho (z)<0\}$.
\end{defi}

One has the following:

\begin{nnthm}
\label{ext-thm} Let $D$ be a  domain in a complex surface $\tilde D$
and let $\Sigma_t$ be a continuous family of $1$-convex
hypersurfaces in $\tilde D$, $t\in (t_1,t_0)$. Suppose that for
every $t\in (t_1,t_0)$ the intersection $\Sigma_t \cap (\tilde
D\setminus D)$ is contained in a relatively compact part of $\tilde
D\setminus D$ and suppose that $\Sigma_{t_0} \subset D$. Let $\cale$
be a sheaf of analytic objects on $D$ and let $E$ be an element of
$\cale|_D$. Then $E$ extends along the family $\{\Sigma_t\}_{t\in
(t_0,t_1)}$.
\end{nnthm}

\smallskip The proof is standard for the standard analytic objects
and will be not given. Now let describe a somewhat non standard
situation in which it will be applied. The idea of the following
construction is inspired by \S 3 from \cite{Iv3}. For  $\eps >0$ and
$\alpha\in (0,\infty)$ consider the following smooth functions
\begin{equation}
\eqqno(rho-alf) \rho_{\eps ,\alpha} (z) = |z_1|^2 - \frac{\eps^2}{4}
- \big(1-\frac{\eps^2}{4}\big)|z_2|^{2\alpha},
\end{equation}
domains
\begin{equation}
\eqqno(d-alfa) D_{\eps , \alpha}^+ = \{ (z_1,z_2)\in \Delta^2:
\rho_{\eps ,\alpha} (z)<0\} \quad\text{ and } \quad  D_{\eps ,
\alpha}^-
\deff \bar\Delta^2\setminus D_{\eps,\alpha},
\end{equation}
and hypersurfaces
\begin{equation}
\eqqno(gam-alf) \Gamma_{\eps , \alpha} = \{ (z_1,z_2)\in \Delta^2:
\rho_{\eps ,\alpha} (z)=0\},
\end{equation}
separating $D_{\eps,\alpha}^+$ from $D_{\eps , \alpha}^-$, see
Figure \ref{gamma-fig}.

\begin{figure}[h]
\centering
\includegraphics[width=2.5in]{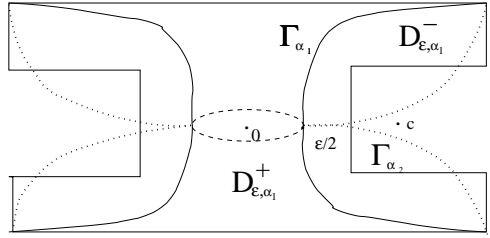}
\caption{For $\alpha_1>>1$ the hypersurface $\Gamma_{\alpha_1}$
(solid line) belongs to $H^2_{\eps}$, for $\alpha_2<1$ (punctured
line) approaches $\{z_2=0\}$ when $\alpha_2\searrow 0$. They are
smooth and strictly pseudoconvex (from the side of $D_{\eps ,
\alpha}^-$) outside of the circle $\{|z_1|=\eps /2, z_2=0\}$.
In the notations of Theorem \ref{ext-thm} $\tilde D = \Delta^2$
and $D=H^2_{\eps}$.}
\label{gamma-fig}
\end{figure}

\begin{lem}
\label{exhaust-l1} \sli For all $\eps > 0$ and $\alpha >0$
hypersurfaces $\Gamma_{\eps ,\alpha}$ are strictly pseudoconvex in
$\Delta^2\setminus (\{|z_1| = \frac{\eps}{2}, z_2=0\})$ from the
side of $D_{\eps , \alpha}^-$.

\smallskip\slii For a fixed $\eps >0$ the domain $D_{\eps , \alpha}$
is contained in $H^2_{\eps}$ for $\alpha$ big enough and
\begin{equation}
\bigcup_{\alpha >0}D_{\eps , \alpha}^+ = \Delta^2\setminus
\big(A_{\frac{\eps}{2},1}\times\{0\}\big).
\end{equation}
\end{lem}
\proof The Levi form of $\rho_{\eps , \alpha}$ at $z$ is

\begin{equation}
\eqqno(levi-form) H_z(\rho_{\eps , \alpha})=
\left(%
\begin{array}{cc}
  1 & 0 \\
  0 & -\alpha^2(1-\eps^2/4)|z_2|^{2\alpha - 2}\\
\end{array}%
\right)
\end{equation}
Since $\d \rho_{\eps , \alpha} = (\bar z_1 , -\alpha
(1-\eps^2/4)\bar z_2|z_2|^{2\alpha - 2})$, a complex tangent vector
to $\Gamma_{\eps ,\alpha}$ at $z$ is $\v_z = (\alpha
(1-\eps^2/4)z_2|z_2|^{2\alpha - 2}, z_1)$.  And therefore
\[
H_z(\rho_{\eps , \alpha})(\v_z,\bar\v_z) =
(1-\eps^2/4)\alpha^2|z_2|^{2\alpha
-2}\big[(1-\eps^2/4)|z_2|^{2\alpha} - |z_1|^2\big] =
-\alpha^2\eps^2/4(1-\eps^2/4)|z_2|^{2\alpha -2}.
\]
The term on the right hand side is negative and therefore the
assertion (\sli of the Lemma is proved. Assertion (\slii is left to
the reader.

\smallskip\qed

Condition $\rho_{\eps,\alpha} >0$ one rewrites as
\begin{equation}
\eqqno(exhaust1) \rho_1(z) \deff \frac{\ln\frac{|z_1|^2 -
\frac{\eps^2}{4}}{1 - \frac{\eps^2}{4}}} {\ln |z_2|^2} < \alpha .
\end{equation}
The preceding calculations mean that the complex Hessian of $\rho_1$
is strictly positive along the complex directions tangent to the
level sets of $\rho_1$. Taking a sufficiently convex function $\psi
:(0,\infty) \to (0,\infty)$ one get a strictly plurisubharmonic
function
\begin{equation}
\eqqno(rho)
\rho \deff \psi (\rho_1)
\end{equation}
in $\Delta \times \Delta^*$ with the same level sets as $\rho_1$.

\smallskip  Let us  formulate the statement we need.

\begin{lem}
\label{exhaust-l2} Let $(X,\call)$ be a straight foliated pair on a
projective surface and let $(\tilde\calp_D,\tilde\call)$ be a
universal covering Poincar\'e domain. Let $(\Delta^2,\cale)$ be a
foliation given by the level sets of a holomorphic function in
$\Delta^2$ with a discrete critical set $\crit (\pi)\subset
\Delta^2\setminus H^2_{\eps}$. Then every foliated immersion $\tilde
h:(H^2_{\eps},\cale) \to (\tilde\calp_D,\tilde\call)$, extends to
$\Delta^2\setminus \crit (\pi)$.
\end{lem}
\proof {\slsf Step 1. Suppose first that $\crit (\pi)$ is empty.}
Take the function $\rho$ as in \eqqref(rho) and let
$D_{\alpha}^{\pm}$ be its upper/lowel level sets. Suppose that
$\tilde h$ is extended up to a point $y_0\in \d D_{\alpha}^+$. Take
the leaf $\cale_{y_0}$ through $y_0$. It cannot be contained in
$D_{\alpha}^-$ in any neighborhood of $y_0$, because $\d
D_{\alpha}^-$ is strictly pseudoconvex at $y_0$. Take a point $y_1
\in D_{\alpha}^+\cap \cale_{y_0}$ close to $y_0$. Trace a small
transversal $L$ to $\cale$ at $y_1$ such that $L\subset
D_{\alpha}^+$. Let $W$ be a foliated bidisc for $\cale$ based on
$L$, \ie for every $y\in L$ one has $\cale_y\cap W$ is the vertical
disc $W_y$. Set $h\deff \tilde p\circ \tilde h$. It extends to a
neighborhood of $y_0$ by the usual Hartogs. Using the fact that
$\tilde p$ is the universal covering we can extend $\tilde
h|_{W_y\cap D_{\alpha}^+} = \tilde p^{-1}\circ h|_{W_y\cap
D_{\alpha}^+}$ with the given initial value $\tilde h(y)$ from
${W_y\cap D_{\alpha}^+}$ onto $W_y$ for every $y$. Using
straightness of $(X\call)$ and therefore rothsteiness of
$(\tilde\calp_D,\tilde\call)$, see Proposition \ref{rothstein4} we
get the extension of $\tilde h$ to a neighborhood of $y_0$. This way
we extend $\tilde h$ to $\Delta^2\setminus A_{\eps,1}\times \{0\}$.
Removal of $A_{\eps ,1}\times \{0\}$ (\ie the Thullen type extension
theorem) doesn't represents any difficulties in this context. I.e.,
$\tilde h$ is extended onto $\Delta^2$.

\smallskip\noindent{\slsf Step 2. Now suppose that $\pi$ has just one
critical point $c_1$, situated away from the Hartogs figure.} I.e.,
such that $c\in \Delta^2\setminus \bar H^2_{\eps}$. Without loss of
generality we may suppose that $c=(1/2,0)$ (with $\eps <<1/2$, see
the Figure \ref{gamma-fig}). Let $\cale$ be the foliation on
$\Delta^2\setminus \{c\}$ by level sets of $\pi$. Take
$(\Delta^2\setminus \{c\},\cale)$ as $(X,\cale)$ and repeat the Step
1 of the proof. Then $\tilde h$ is extended again to
$\Delta^2\setminus A_{\eps, 1}\times \{0\}$.  Removal of $(A_{\eps
,1}\times \{0\})\setminus \{c\}$ is again obvious.

\smallskip\noindent{\slsf Step 3. The $\crit (\pi)$ is descrete.} Then
do the same up to arriving to the first critical point $c_1$ at
value $\alpha_1$. Then  place appropriately a Hartogs figure into
$D_{\alpha_1}^+$ in order to obtain the situation of Step 2 and
therefore extend $\tilde h$ to a punctured neighborhood of $c_1$.
The rest is obvious.

\smallskip\qed

\newprg[prgRAT.proof]{Pseudoconvexity of the universal covering
Poincar\'e domains}

We shall prove first a general statement about appearance of
invariant rational curves as obstructions to local pseudoconvexity
of covering Poincar\'e domains.

\begin{nnthm}
\label{ps-conv-rat} Let  $(X,\call)$ be a straight projective pair
and let $m$ be a point in $X^{\reg}$. Then the following
statements are equivalent:

\smallskip\sli  For every Poincar\'e disc $D$ through $m$ $\tilde\calp_D$
is not locally pseudoconvex over $X$:

\smallskip\slii $\overline{\call_m}$ is a rational curve, cutting
$\Sing\call$ exactly at one point (it is the very same point over which
all $\tilde\calp_D$ are not pseudoconvex).
\end{nnthm}
\proof This will be done in several steps.

\smallskip\noindent{\slsf Step 1. Pseudoconvexity of $\tilde\calp_D$
over non-singular points.} $\tilde\calp_D$  was defined in
Subsection \ref{prgPPD.pd} and $\tilde p : \tilde\calp_D \to
\calp_D\subset X^{\reg}$ denotes the canonical map. The pair
$(\tilde\calp_D , \tilde p)$ is a  Riemann domain over $X$. First of
all let us remark that, as in the case with $\calp_D$ the covering
Poincar\'e domain $\tilde\calp_D$ is always pseudoconvex over non
singular points of $\call$.

\begin{lem}
\label{loc-ps2} If $z_0\not\in \Sing\call$ then $\tilde\calp_D$ is
locally pseudoconvex over $z_0$.
\end{lem}
\proof Pseudoconvexity of $\calp_D$ at such $z_0$ was proved in
Lemma \ref{loc-ps1}. Let us see that an analogous proof goes through
also for $\tilde\calp_D$. Indeed, let $z_0\not\in \Sing\call$. Take
then, as in the proof of Lemma \ref{loc-ps1}, a foliated bidisc
$U\ni z_0$. Set, as above, $\calf\deff \call|_{U}$. Let $U_1$ be a
connected component of $\tilde p^{-1}(U)$ and let $V$ be the image
of $U_1$  under $\tilde p|_{U_1}:U_1\to U$.
 If $z=(z_1,z_2)\in V$ is the
image of some $w\in U_1$ then $\tilde p^{-1}|_{\calf_{z_1}}$ (with
initial value $w$) extends along $\calf_{z_1}\deff \{z_1\}\times
\Delta$ because of simple connectivity of $\Delta$, and $\tilde
p^{-1}|_{\calf_{z_1}}(\calf_{z_1})\subset U_1$ because of the
connectivity of $U_1$. Therefore for every $z=(z_1,z_2)\in V$ we
have that $\calf_{z_1}\subset V$ and $\tilde
p^{-1}|_{\calf_{z_1}}(\calf_{z_1})$ is a disjoint union of discs,
each of them is mapped by $\tilde p$ biholomorphically onto
$\calf_{z_1}$. Since $V$ is connected, we have that $V=V_1\times
\Delta$ for some open, connected $V_1\subset \Delta$. A connected
component $U_1$ is also a connected component of $\tilde p^{-1}(V)$
and it is foliated by discs - preimages of leaves of $\calf$. The
restriction $\tilde p|_{U_1}: U_1\to V$ is a foliated local
biholomorphism, \ie $(U_1, \tilde p|_{U_1})$ is a Riemann domain
over a Stein manifold $V=V_1\times \Delta$. Foliation on $U_1$ we
denote as $\tilde\calf$ - it this the restriction of the universal
foliation $\tilde\call$ to $U_1$.

\smallskip To prove that $U_1$ is $p_7$-convex consider a
holomorphic imbedding $\tilde h:H^2_{\eps}\to U_1$. Let $\cale$ be
the pull back of $\tilde\calf$ to $H^2_{\eps}$ by $\tilde h$. Since
it is the same as pull back of $\call$ by the extended map $
h:\Delta^2 \to V$ (here $h$ stands for the extension of $\tilde
p\circ \tilde h$) we conclude that $\cale$ extends to a {\slsf
smooth} foliation on $\Delta^2$. Let $w$ be a point in $\Delta^2$.
Take a leaf $\cale_w$. It intersects $H^2_{\eps}$ by Lemma
\ref{intersect}. Let $u$ be some point of this intersection and let
$v= h(u)$. Take the lift $\tilde p|_{\calf_v}$ with the initial
value $h(u)$. The image of this lift is a disc - a leaf of
$\tilde\calf$ and it is entirely contained in $U_1$ by preceding
considerations. Since $\tilde h|_{\cale_w\cap H^2_{\eps}} = \tilde
p^{-1}|_{\calf_v}\circ h|_{\cale_w}$ we see that $\tilde
h|_{\cale_w\cap H^2_{\eps}}$ extends onto $\cale_w$ as a mapping
with values in $U_1$ (because $\tilde p^{-1}|_{\calf_v}$ takes its
values in $U_1$). This proves, via Rothsteins-type Proposition
\ref{rothstein2}, the extendability of $\tilde h$ onto $\Delta^2$ as
a mapping with values in $U_1$. Indeed, $(U_1,\tilde\call)$ is
Rothstein, because it is a foliated domain over the Stein pair $(U,
\call)$.

\smallskip We proved that $U_1$ is $p_7$-converx and therefore it
is Stein.

\smallskip\qed

\medskip\noindent{\slsf Step 2. Nearby rational curves.}
Now we shall prove the following:

\begin{lem}
\label{loc-ps3} Let $(X,\call)$ be a straight foliated pair on a
projective surface $X$. Let $m$ be a point in $X^{\reg}$ and
$D\ni m$ a transversal to $\call$, locally closed disc. Suppose that
the universal covering Poincar\'e domain $(\tilde\calp_D, \tilde p)$
 is not locally pseudoconvex over a point $z_0\in X$. Then:

\smallskip\sli  $z_0\in\Sing\call$ and $z_0$ is an isolated point of
$\d\calp_D$;

\smallskip\slii for some $m_1\in D$ the closure $\overline{\call_{m_1}}$
is a rational curve passing through $z_0$.
\end{lem}
\proof We already know that such $z_0$ should belong to $\Sing\call$
and, in particular, $z_0$ should be a boundary  point of $\calp_D$.
Take some small neighborhood $U$ of $z_0$ biholomorphic to a ball
and not containing any other then $z_0$ points of $\Sing\call$.

\smallskip Let $U_1$ be some  connected component of $\tilde p^{-1}(U)$
and let $\tilde h:H^2_{\eps}\to U_1$ be a holomorphic imbedding. Take
$\tilde p\circ \tilde h$ and extend it by Hartogs theorem to a locally
biholomorphic mapping $ h :\Delta^2\to U\subset X$. Let
$\cale\deff  h^*\calf$ be the induced foliation on $\Delta^2$,
where $\calf = \call|_U$. Set $S=\Sing\cale$, it is clear that
$S= h^{-1}(z_0)$ if it is nonempty, \ie if $z_0\in
h(\Delta^2)$. Moreover, it is clear that $\cale|_{H^2_{\eps}} =
\tilde h^*\tilde\calf$, where $\tilde\calf \deff \tilde\call|_{U_1}$.
Remark furthermore that the universal foliation $\tilde\call$ of
$\tilde\calp_D$ possesses a first integral, namely $\tilde\pi$.
Therefore $\cale$ possesses it to, it is nothing but $\tilde\pi\circ
\tilde h$ extended from $H^2_{\eps}$ to $\Delta^2$. Denote this integral as
$\pi$.

\smallskip Shrinking $\Delta^2$ arbitrarily slightly we can suppose that
$S$ is finite. It is clear that $S = \crit (\pi)$. Applying Lemma
\ref{exhaust-l2} to our $\pi$ we  extend $\tilde h$ as a locally
biholomorphic map $\tilde h:\Delta^2\setminus S\to U_1$.

\begin{figure}[h]
\centering
\includegraphics[width=1.7in]{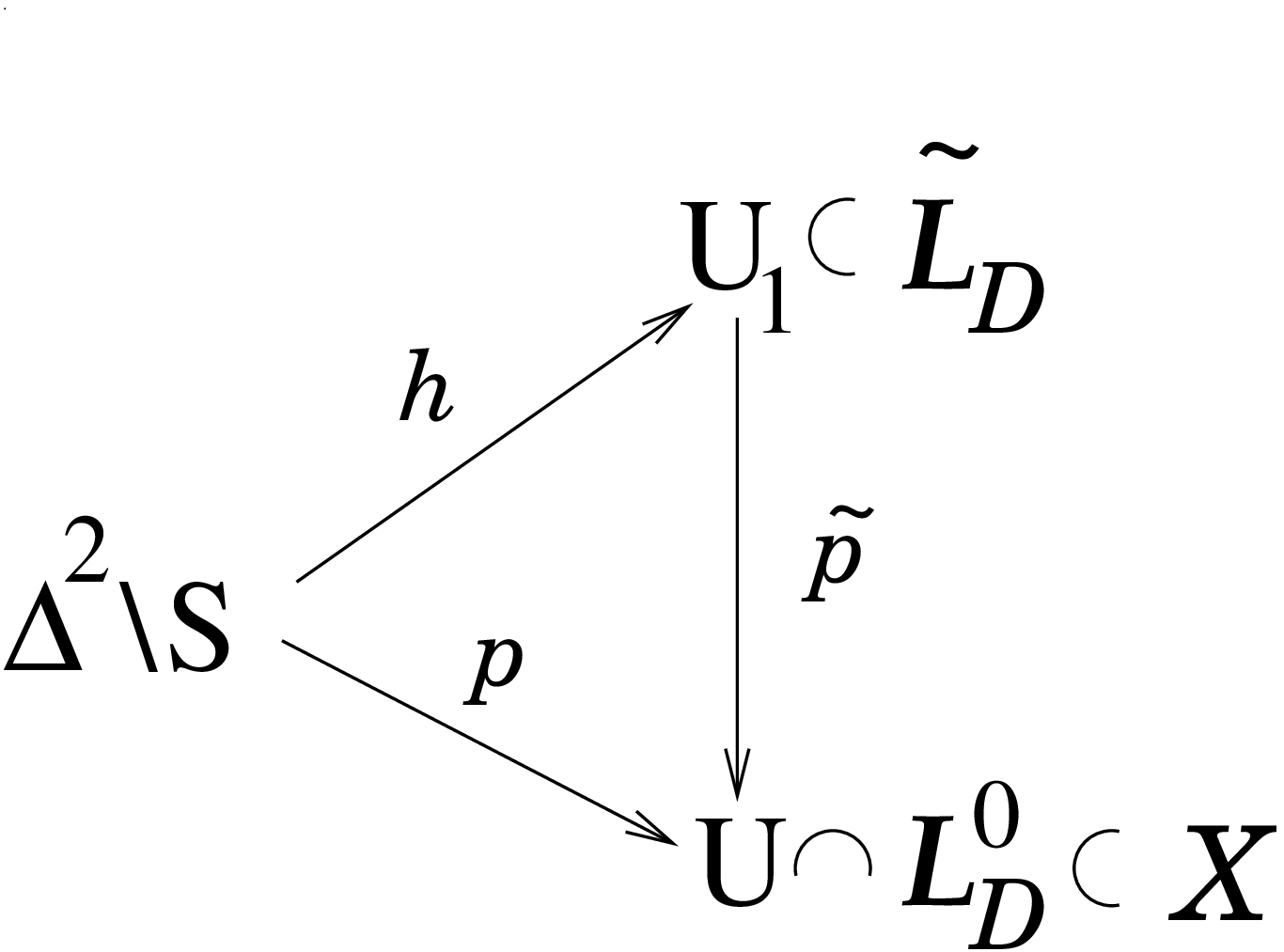}
\caption{A diagram relating $h$, $\tilde p$ and $\tilde h$.}
\label{diagram}
\end{figure}

\smallskip Now we have two cases.

\smallskip\noindent{\slsf Case 1. The point $z_0$ is not an
isolated boundary point of $\calp_D$.}  In that case $U\cap \calp_D$
is Stein, as we already know from Lemma \ref{loc-ps1}, and
consequently $h(\Delta^2)\subset U\cap \calp_D$. Therefore
$S=\emptyset$ and $\tilde h$ extends to a locally biholomorphic mapping
$\tilde h:\Delta^2\to U_1$. I.e., $U_1$ is $p_7$-convex and the
Docquier-Grauert criterion provides the steiness of $U_1$.

\smallskip\noindent{\slsf Case 2. $z_0\in\Sing\call$ and is an
isolated boundary point of $\calp_D$.} Let $U$ be as above. Suppose
that there is a connected component $U_1$ of $\tilde
p^{-1}(\calp_D\cap U)$ which is not Stein. Let $\tilde h:H^2_{\eps}\to U_1$
and  $h \deff \tilde p\circ \tilde h$ be as above and let
$\tilde h:\Delta^2\setminus S\to U_1 \subset \tilde\calp_D$ be the
extension of $\tilde h$ constructed at the beginning of the proof.

\smallskip
This extension is proper near $S$. Remark that $S$ cannot be empty
due to the assumed non Steiness of $U_1$. What concerns properness,
indeed, let $y_n\in\Delta^2$ be a sequence such that $y_n\to s_0\in
S$. We need to prove that $\tilde h(y_n)$ leave every compact in
$\tilde\calp_D$. But $ h =\tilde p\circ \tilde h $ is a locally
biholomorphic map to $X$. Therefore it is biholomorphic in a
neighborhood of $s_0$ and $ h(s_0) = z_0$. Take neighborhoods $W\ni
s_0$ and $V\ni z_0$, biholomorphic to a ball, such that $h|_W:W\to
V$ is a biholomorphism. Set $z_n\deff h(y_n)$ and $x_n\deff \tilde
h(y_n) \in U_1$. Then $\tilde p(x_n) = z_n\to z_0$. Would $x_n$ stay
in a compact part of $\tilde\calp_D$ their images $z_n=\tilde
p(z_n)$ could not approach a singular point $z_0$ of our foliation.
Therefore $\tilde h(y_n)$ leaves every compact in $\tilde\calp_D$.

\smallskip Now we can attach the set $S$ to $U_1$ and get a
completed domain  $\bar U_1$ over $X$, and consequently we can
complete $\tilde\calp_D$. Indeed, for $s_0$ and $z_0= h(s_0)$ take
$W$ and $V$ as above. Then $\tilde p^{-1}\circ h$ lifts over
$W\setminus \{s_0\}$ to an injective map $\tilde p^{-1}\circ
h|_{W\setminus \{s_0\}}: W\setminus \{s_0\} \to U_1$, which is
proper at $s_0$. Now we can attach $s_0$ to $U_1$ by this map.

\smallskip The projection $\tilde p$ will extend as a
locally biholomorphic map to  $\bar U_1$. We can repeat this
procedure with all components of $\tilde p^{-1}(U)$ which are not
Stein. The obtained Riemann domain over $X$ we denote by
$(\bar\calp_D,\bar p)$ ($\bar p$ stands for the extension of $\tilde
p$).

\smallskip $\bar\calp_D\setminus \tilde\calp_D$ is a discrete, non
empty subset of $\bar\calp_D$. We have the holomorphic projection
$\tilde\pi : \tilde\calp_D\to D$ onto the base $D$ of the Poincar\'e
domain $\tilde\calp_D$ and it holomorphically extends to
$\bar\calp_D$. Denote by $\bar\pi$ the extended map. Take some point
$s_0\in \bar\calp_D\setminus \tilde\calp_D$. Set $t_0=\bar\pi (s_0)$
and consider the complex curve $C_{t_0}\deff \bar\pi^{-1}(t_0)$. It
can be nothing else but a $1$-point compactification of a simply
connected leaf $\tilde\call_{t_0}$ by $s_0$. I.e., $C_{t_0}$ can be
only a rational curve. Note that $C_{t_0}$ passes through $s_0$.
Therefore $\overline{\call_{t_0}} = \bar p(C_{t_0})$ is an invariant
rational curve in $X$ passing through $z_0$. All what is left is to
set $m_1=t_0$.

\smallskip\qed

\smallskip\noindent{\slsf Step 3. Rationality of $\overline{\call_m}$.}
Write $D_1\deff D$ and $C_1\deff C_{m_1}$.  Let $D_{k}$ denotes the
subdisc of $D_1$ of geodesics radius $1/k$ with the same center $m$.
$\tilde\calp_{D_k}$ is then a subdomain of $\tilde\calp_{D_1}$.
Suppose that for all $k$ Poincar\'e domains $\tilde\calp_{D_k}$ are
not pseudoconvex over $z_0$, in particular, that $z_0$ is an
isolated boundary point for all $\calp_{D_k}$.

\smallskip Repeating the previous considerations for every $k$
 we get $m_k\in D_k$, a rational curve $C_k$ in
$\bar\calp_{D_k}\subset \bar\calp_D$, which is a one point
compactification of $\tilde\call_{m_k}$. $C_k$ passes through $s_k$
such that $\bar p(s_k) = z_0$ for all $k$. More accurately: $\bar
p_k(s_k)=z_0$. But $\bar p_k$ is the restriction of $\bar p$ to
$\bar\calp_{D_k}$. The images $\bar p(C_k) = \overline{\call_{m_k}}$
are invariant rational curves in $X$ passing through $z_0$. And
$z_0$ is the only point of $\overline{\call_{m_k}}\cap \Sing\call$.

\begin{rema} \rm
All $s_k$ may be well distinct. It is true that they all project to
$z_0$ under $\bar p$, but $\bar\calp_D$ might be infinite sheeted
over $X$.
\end{rema}

\smallskip If some subsequence of $\bar p(C_k)$ stabilizes then it
is equal to $\overline{\call_m}$ in fact, and we are done. If not,
then we got a sequence of distinct invariant rational curves through
$z_0$. But then $\call$ is a rational quasi-fibration and again we
are done. We proved that in this case $\overline{\call_m}$ is an
invariant rational curve passing through $z_0$ as predicts Part
(\slii of Theorem \ref{ps-conv-rat}.

\smallskip\qed

\newprg[prgRAT.approach]{Proof of Corollary \ref{approach3}}
Suppose that in the conditions of Corollary \ref{approach1} the case
(\sli occurs and, moreover, that $\bar\calp_D\setminus \calp_D \not=
\emptyset$. Denote by $\{s_1,..,s_n\}$ the set of all nondicritical
points of this set. We shall prove that this set is empty. If not
let $s\in \Sing\call$ appeared to be an isolated boundary point of
the Poincar\'e domain $\calp_D$ such that $s$ is not dicritical.

\smallskip Consider first the case when the nef model $(Y,\calf)$ of
$(X,\call)$ is parabolic. From the minimality of
$\overline{\call_m}$ we obtain readily that $\lim\call_m\subset
\Sing\call$ and therefore $\overline{\call}$ is a rational curve
cutting $\Sing\call$ by at least two points, \ie the case (\slii
 of Corollary \ref{approach1} occurs.

\smallskip Let us turn to the hyperbolic case. Using contractibility
of the holonomy we can, as before, forbid any invariant rational
curves entering to the Poincar\'e domains appearing in the process
of the nef reduction (by taking smaller $D$). Let $T$ be the tree of
rational curves appeared in the process of Seidenberg's reduction of
singularities over $s$. If $T$ is entirely contracted by the
subsequent modification then all what happens to $s$ is that it is
replaced by a cyclic point. But then it can be only a smooth point.
In this case the universal Poincar\'e domain $\tilde\calp_D$ cannot
be Hausdorff. Contradiction.

\smallskip Therefore $T$ is divided to subtrees $T_i$ of invariant
rational curves, subsequently contracted to cyclic points $c_i$ and
connected by some other rational chains, which are not contracted.
Through each of $c_i$ passes then an $\calf$-invariant rational
curve, which cannot belong to $\calp_D^{\calf}$ and therefore no one
of $c_i$ doesn't belong to $\calp_D$, as well as all these
"connecting" invariant curves. Therefore $\calp_D^{\calf}$ contains
a punctured ball. In the hyperbolic case that means that
$\tilde\calp_{\cc}^{\calf}$ is not Stein. $\tilde\calp_{\cc}$ is the
same as $\tilde\calp_{\cc}^{\calf}$ over this ball, \ie is also not
Stein. But then $\overline{\call_m}$ is a rational curve cutting
$\Sing\call$ by exactly one point (this point is $s$). This is
impossible because $\pi_1(\call_m, m)$ cannot be trivial.
Contradiction.

Corollary \ref{approach3} is proved.

\ifx\undefined\bysame
\newcommand{\bysame}{\leavevmode\hbox to3em{\hrulefill}\,}
\fi

\def\entry#1#2#3#4\par{\bibitem[#1]{#1}
{\textsc{#2 }}{\sl{#3} }#4\par\vskip2pt}
%{ref}{author}{title}ref.

\end{document}